\input amstex
\input amsppt.sty

\magnification1200
\hsize391.46176pt
\vsize536.00175pt

\comment
\catcode`\@=11

\def\remarkheadfont@{\smc}
\let\varindent@\indent
\def\definition{\let\savedef@\definition \let\definition\relax
  \add@missing\endproclaim \add@missing\endroster
  \add@missing\enddefinition \envir@stack\enddefinition
   \def\definition##1{\restoredef@\definition
     \penaltyandskip@{-100}\medskipamount
        {\def\usualspace{{\remarkheadfont@\enspace}}%
        \varindent@\remarkheadfont@\ignorespaces##1\unskip
        \frills@{.\remarkheadfont@\enspace}}%
        \rm \ignorespaces}%
  \nofrillscheck\definition}
\catcode`\@=13
\endcomment

\def\SunaAE{39}
\def\StemAE{38}
\def\StemAD{37}
\def\StanAI{36}
\def\SlatAC{35}
\def\SagaAL{34}
\def\ProcAK{33}
\def\ProcAI{32}
\def\ProcAB{31}
\def\ProcAF{30}
\def\ProcAD{29}
\def\ProcAE{28}
\def\OkadAI{27}
\def\MacdAC{26}
\def\LittAA{25}
\def\LitPAA{24}
\def\LaShAC{23}
\def\LaMSAB{22}
\def\LaksAB{21}
\def\KratAW{20}
\def\KratAP{19}
\def\KratAQ{18}
\def\KratAO{17}
\def\KoTeAB{16}
\def\KoTeAA{15}
\def\KnutAA{14}
\def\KingAC{13}
\def\KingAB{12}
\def\KingAD{11}
\def\KingAA{10}
\def\IsWaAA{9}
\def\HaSuAA{8}
\def\GoulAD{7}
\def\FuHaAA{6}
\def\FuKrAA{5}
\def\DePrAA{4}
\def\DeCoAA{3}
\def\CariAA{2}
\def\BrGrAA{1}

\input amstex
\catcode`\@=11
\font\tenln    = line10
\font\tenlnw   = linew10

\newskip\Einheit \Einheit=0.5cm
\newcount\xcoord \newcount\ycoord
\newdimen\xdim \newdimen\ydim \newdimen\PfadD@cke \newdimen\Pfadd@cke

\newcount\@tempcnta

\newdimen\@tempdima
\newdimen\@tempdimb

\newdimen\@wholewidth
\newdimen\@halfwidth

\newcount\@xarg
\newcount\@yarg
\newcount\@yyarg
\newbox\@linechar
\newdimen\@linelen
\newdimen\@clnwd
\newdimen\@clnht

\newif\if@negarg

\def\@whilenoop#1{}
\def\@whiledim#1\do #2{\ifdim #1\relax#2\@iwhiledim{#1\relax#2}\fi}
\def\@iwhiledim#1{\ifdim #1\let\@nextwhile=\@iwhiledim
        \else\let\@nextwhile=\@whilenoop\fi\@nextwhile{#1}}

\def\@whileswnoop#1\fi{}
\def\@whilesw#1\fi#2{#1#2\@iwhilesw{#1#2}\fi\fi}
\def\@iwhilesw#1\fi{#1\let\@nextwhile=\@iwhilesw
         \else\let\@nextwhile=\@whileswnoop\fi\@nextwhile{#1}\fi}

\def\thinlines{\let\@linefnt\tenln \let\@circlefnt\tencirc
  \@wholewidth\fontdimen8\tenln \@halfwidth .5\@wholewidth}
\def\thicklines{\let\@linefnt\tenlnw \let\@circlefnt\tencircw
  \@wholewidth\fontdimen8\tenlnw \@halfwidth .5\@wholewidth}
\thinlines

\PfadD@cke1pt \Pfadd@cke0.5pt
\def\PfadDicke#1{\PfadD@cke#1 \divide\PfadD@cke by2 \Pfadd@cke\PfadD@cke \multiply\PfadD@cke by2}
\long\def\LOOP#1\REPEAT{\def\BODY{#1}\ITERATE}
\def\ITERATE{\BODY \let\next\ITERATE \else\let\next\relax\fi \next}
\let\REPEAT=\fi
\def\Punkt{\hbox{\raise-2pt\hbox to0pt{\hss$\ssize\bullet$\hss}}}
\def\DuennPunkt(#1,#2){\unskip
  \raise#2 \Einheit\hbox to0pt{\hskip#1 \Einheit
          \raise-2.5pt\hbox to0pt{\hss$\bullet$\hss}\hss}}
\def\NormalPunkt(#1,#2){\unskip
  \raise#2 \Einheit\hbox to0pt{\hskip#1 \Einheit
          \raise-3pt\hbox to0pt{\hss\twelvepoint$\bullet$\hss}\hss}}
\def\DickPunkt(#1,#2){\unskip
  \raise#2 \Einheit\hbox to0pt{\hskip#1 \Einheit
          \raise-4pt\hbox to0pt{\hss\fourteenpoint$\bullet$\hss}\hss}}
\def\Kreis(#1,#2){\unskip
  \raise#2 \Einheit\hbox to0pt{\hskip#1 \Einheit
          \raise-4pt\hbox to0pt{\hss\fourteenpoint$\circ$\hss}\hss}}

\def\Line@(#1,#2)#3{\@xarg #1\relax \@yarg #2\relax
\@linelen=#3\Einheit
\ifnum\@xarg =0 \@vline
  \else \ifnum\@yarg =0 \@hline \else \@sline\fi
\fi}

\def\@sline{\ifnum\@xarg< 0 \@negargtrue \@xarg -\@xarg \@yyarg -\@yarg
  \else \@negargfalse \@yyarg \@yarg \fi
\ifnum \@yyarg >0 \@tempcnta\@yyarg \else \@tempcnta -\@yyarg \fi
\ifnum\@tempcnta>6 \@badlinearg\@tempcnta0 \fi
\ifnum\@xarg>6 \@badlinearg\@xarg 1 \fi
\setbox\@linechar\hbox{\@linefnt\@getlinechar(\@xarg,\@yyarg)}%
\ifnum \@yarg >0 \let\@upordown\raise \@clnht\z@
   \else\let\@upordown\lower \@clnht \ht\@linechar\fi
\@clnwd=\wd\@linechar
\if@negarg \hskip -\wd\@linechar \def\@tempa{\hskip -2\wd\@linechar}\else
     \let\@tempa\relax \fi
\@whiledim \@clnwd <\@linelen \do
  {\@upordown\@clnht\copy\@linechar
   \@tempa
   \advance\@clnht \ht\@linechar
   \advance\@clnwd \wd\@linechar}%
\advance\@clnht -\ht\@linechar
\advance\@clnwd -\wd\@linechar
\@tempdima\@linelen\advance\@tempdima -\@clnwd
\@tempdimb\@tempdima\advance\@tempdimb -\wd\@linechar
\if@negarg \hskip -\@tempdimb \else \hskip \@tempdimb \fi
\multiply\@tempdima \@m
\@tempcnta \@tempdima \@tempdima \wd\@linechar \divide\@tempcnta \@tempdima
\@tempdima \ht\@linechar \multiply\@tempdima \@tempcnta
\divide\@tempdima \@m
\advance\@clnht \@tempdima
\ifdim \@linelen <\wd\@linechar
   \hskip \wd\@linechar
  \else\@upordown\@clnht\copy\@linechar\fi}

\def\@getlinechar(#1,#2){\@tempcnta#1\relax\multiply\@tempcnta 8
\advance\@tempcnta -9 \ifnum #2>0 \advance\@tempcnta #2\relax\else
\advance\@tempcnta -#2\relax\advance\@tempcnta 64 \fi
\char\@tempcnta}

\def\Diagonale(#1,#2)#3{\unskip\leavevmode
  \xcoord#1\relax \ycoord#2\relax
      \raise\ycoord \Einheit\hbox to0pt{\hskip\xcoord \Einheit
         \Line@(1,1){#3}\hss}}
\def\AntiDiagonale(#1,#2)#3{\unskip\leavevmode
  \xcoord#1\relax \ycoord#2\relax \advance\xcoord by -0.05\relax
      \raise\ycoord \Einheit\hbox to0pt{\hskip\xcoord \Einheit
         \Line@(1,-1){#3}\hss}}
\def\Pfad(#1,#2),#3\endPfad{\unskip\leavevmode
  \xcoord#1 \ycoord#2 \thicklines\ZeichnePfad#3\endPfad\thinlines}
\def\ZeichnePfad#1{\ifx#1\endPfad\let\next\relax
  \else\let\next\ZeichnePfad
    \ifnum#1=1
      \raise\ycoord \Einheit\hbox to0pt{\hskip\xcoord \Einheit
         \vrule height\Pfadd@cke width1 \Einheit depth\Pfadd@cke\hss}%
      \advance\xcoord by 1
    \else\ifnum#1=2
      \raise\ycoord \Einheit\hbox to0pt{\hskip\xcoord \Einheit
        \hbox{\hskip-\PfadD@cke\vrule height1 \Einheit width\PfadD@cke depth0pt}\hss}%
      \advance\ycoord by 1
    \else\ifnum#1=3
      \raise\ycoord \Einheit\hbox to0pt{\hskip\xcoord \Einheit
         \Line@(1,1){1}\hss}
      \advance\xcoord by 1
      \advance\ycoord by 1
    \else\ifnum#1=4
      \raise\ycoord \Einheit\hbox to0pt{\hskip\xcoord \Einheit
         \Line@(1,-1){1}\hss}
      \advance\xcoord by 1
      \advance\ycoord by -1
    \fi\fi\fi\fi
  \fi\next}
\def\hSSchritt{\leavevmode\raise-.4pt\hbox to0pt{\hss.\hss}\hskip.2\Einheit
  \raise-.4pt\hbox to0pt{\hss.\hss}\hskip.2\Einheit
  \raise-.4pt\hbox to0pt{\hss.\hss}\hskip.2\Einheit
  \raise-.4pt\hbox to0pt{\hss.\hss}\hskip.2\Einheit
  \raise-.4pt\hbox to0pt{\hss.\hss}\hskip.2\Einheit}
\def\vSSchritt{\vbox{\baselineskip.2\Einheit\lineskiplimit0pt
\hbox{.}\hbox{.}\hbox{.}\hbox{.}\hbox{.}}}
\def\DSSchritt{\leavevmode\raise-.4pt\hbox to0pt{%
  \hbox to0pt{\hss.\hss}\hskip.2\Einheit
  \raise.2\Einheit\hbox to0pt{\hss.\hss}\hskip.2\Einheit
  \raise.4\Einheit\hbox to0pt{\hss.\hss}\hskip.2\Einheit
  \raise.6\Einheit\hbox to0pt{\hss.\hss}\hskip.2\Einheit
  \raise.8\Einheit\hbox to0pt{\hss.\hss}\hss}}
\def\dSSchritt{\leavevmode\raise-.4pt\hbox to0pt{%
  \hbox to0pt{\hss.\hss}\hskip.2\Einheit
  \raise-.2\Einheit\hbox to0pt{\hss.\hss}\hskip.2\Einheit
  \raise-.4\Einheit\hbox to0pt{\hss.\hss}\hskip.2\Einheit
  \raise-.6\Einheit\hbox to0pt{\hss.\hss}\hskip.2\Einheit
  \raise-.8\Einheit\hbox to0pt{\hss.\hss}\hss}}
\def\SPfad(#1,#2),#3\endSPfad{\unskip\leavevmode
  \xcoord#1 \ycoord#2 \ZeichneSPfad#3\endSPfad}
\def\ZeichneSPfad#1{\ifx#1\endSPfad\let\next\relax
  \else\let\next\ZeichneSPfad
    \ifnum#1=1
      \raise\ycoord \Einheit\hbox to0pt{\hskip\xcoord \Einheit
         \hSSchritt\hss}%
      \advance\xcoord by 1
    \else\ifnum#1=2
      \raise\ycoord \Einheit\hbox to0pt{\hskip\xcoord \Einheit
        \hbox{\hskip-2pt \vSSchritt}\hss}%
      \advance\ycoord by 1
    \else\ifnum#1=3
      \raise\ycoord \Einheit\hbox to0pt{\hskip\xcoord \Einheit
         \DSSchritt\hss}
      \advance\xcoord by 1
      \advance\ycoord by 1
    \else\ifnum#1=4
      \raise\ycoord \Einheit\hbox to0pt{\hskip\xcoord \Einheit
         \dSSchritt\hss}
      \advance\xcoord by 1
      \advance\ycoord by -1
    \fi\fi\fi\fi
  \fi\next}
\def\Koordinatenachsen(#1,#2){\unskip
 \hbox to0pt{\hskip-.5pt\vrule height#2 \Einheit width.5pt depth1 \Einheit}%
 \hbox to0pt{\hskip-1 \Einheit \xcoord#1 \advance\xcoord by1
    \vrule height0.25pt width\xcoord \Einheit depth0.25pt\hss}}
\def\Koordinatenachsen(#1,#2)(#3,#4){\unskip
 \hbox to0pt{\hskip-.5pt \ycoord-#4 \advance\ycoord by1
    \vrule height#2 \Einheit width.5pt depth\ycoord \Einheit}%
 \hbox to0pt{\hskip-1 \Einheit \hskip#3\Einheit 
    \xcoord#1 \advance\xcoord by1 \advance\xcoord by-#3 
    \vrule height0.25pt width\xcoord \Einheit depth0.25pt\hss}}
\def\Gitter(#1,#2){\unskip \xcoord0 \ycoord0 \leavevmode
  \LOOP\ifnum\ycoord<#2
    \loop\ifnum\xcoord<#1
      \raise\ycoord \Einheit\hbox to0pt{\hskip\xcoord \Einheit\Punkt\hss}%
      \advance\xcoord by1
    \repeat
    \xcoord0
    \advance\ycoord by1
  \REPEAT}
\def\Gitter(#1,#2)(#3,#4){\unskip \xcoord#3 \ycoord#4 \leavevmode
  \LOOP\ifnum\ycoord<#2
    \loop\ifnum\xcoord<#1
      \raise\ycoord \Einheit\hbox to0pt{\hskip\xcoord \Einheit\Punkt\hss}%
      \advance\xcoord by1
    \repeat
    \xcoord#3
    \advance\ycoord by1
  \REPEAT}
\def\Label#1#2(#3,#4){\unskip \xdim#3 \Einheit \ydim#4 \Einheit
  \def\lo{\advance\xdim by-.5 \Einheit \advance\ydim by.5 \Einheit}%
  \def\llo{\advance\xdim by-.25cm \advance\ydim by.5 \Einheit}%
  \def\loo{\advance\xdim by-.5 \Einheit \advance\ydim by.25cm}%
  \def\o{\advance\ydim by.25cm}%
  \def\ro{\advance\xdim by.5 \Einheit \advance\ydim by.5 \Einheit}%
  \def\rro{\advance\xdim by.25cm \advance\ydim by.5 \Einheit}%
  \def\roo{\advance\xdim by.5 \Einheit \advance\ydim by.25cm}%
  \def\l{\advance\xdim by-.30cm}%
  \def\r{\advance\xdim by.30cm}%
  \def\lu{\advance\xdim by-.5 \Einheit \advance\ydim by-.6 \Einheit}%
  \def\llu{\advance\xdim by-.25cm \advance\ydim by-.6 \Einheit}%
  \def\luu{\advance\xdim by-.5 \Einheit \advance\ydim by-.30cm}%
  \def\u{\advance\ydim by-.30cm}%
  \def\ru{\advance\xdim by.5 \Einheit \advance\ydim by-.6 \Einheit}%
  \def\rru{\advance\xdim by.25cm \advance\ydim by-.6 \Einheit}%
  \def\ruu{\advance\xdim by.5 \Einheit \advance\ydim by-.30cm}%
  #1\raise\ydim\hbox to0pt{\hskip\xdim
     \vbox to0pt{\vss\hbox to0pt{\hss$#2$\hss}\vss}\hss}%
}
\catcode`\@=13

\def\LL{\leavevmode\setbox0=\hbox{L}\hbox to\wd0{\hss\char'40L}}
\def\al{\alpha}

\def\la{\lambda}
\def\rh{\rho}
\def\si{\sigma}

\def\om{\omega}

\def\Th{\Theta}

\def\Ph{\Phi}
\def\Ps{\Psi}
\def\Om{\Omega}
\def\Up{\Upsilon}


\def\C{{\Bbb C}}

\def\today{\ifcase\month\or
 January\or February\or March\or April\or May\or June\or
 July\or August\or September\or October\or November\or December\fi
 \space\number\day, \number\year}
\def\dfrac#1#2{{\displaystyle{#1\over#2}}}

\def\({\left(}
\def\){\right)}
\def\[{\left[}
\def\]{\right]}

\def\Sp{\operatorname{Sp}}

\def\3{\ss}
\catcode`\@=11
\def\dddot#1{\vbox{\ialign{##\crcr
      .\hskip-.5pt.\hskip-.5pt.\crcr\noalign{\kern1.5\p@\nointerlineskip}
      $\hfil\displaystyle{#1}\hfil$\crcr}}}

\newif\iftab@\tab@false
\newif\ifvtab@\vtab@false
\def\tab{\bgroup\tab@true\vtab@false\vst@bfalse\Strich@false%
   \def\\{\global\hline@@false%
     \ifhline@\global\hline@false\global\hline@@true\fi\cr}
   \edef\l@{\the\leftskip}\ialign\bgroup\hskip\l@##\hfil&&##\hfil\cr}
\def\endtab{\cr\egroup\egroup}
\def\vtab{\vtop\bgroup\vst@bfalse\vtab@true\tab@true\Strich@false%
   \bgroup\def\\{\cr}\ialign\bgroup&##\hfil\cr}
\def\endvtab{\cr\egroup\egroup\egroup}
\def\stab{\D@cke0.5pt\null 
 \bgroup\tab@true\vtab@false\vst@bfalse\Strich@true\Let@@\vspace@
 \normalbaselines\offinterlineskip
  \openup\spreadmlines@
 \edef\l@{\the\leftskip}\ialign
 \bgroup\hskip\l@##\hfil&&##\hfil\crcr}
\def\endstab{\crcr\egroup
 \egroup}
\newif\ifvst@b\vst@bfalse
\def\vstab{\D@cke0.5pt\null
 \vtop\bgroup\tab@true\vtab@false\vst@btrue\Strich@true\bgroup\Let@@\vspace@
 \normalbaselines\offinterlineskip
  \openup\spreadmlines@\bgroup}
\def\endvstab{\crcr\egroup\egroup
 \egroup\tab@false\Strich@false}

\newdimen\htstrut@
\htstrut@8.5\p@
\newdimen\htStrut@
\htStrut@12\p@
\newdimen\dpstrut@
\dpstrut@3.5\p@
\newdimen\dpStrut@
\dpStrut@3.5\p@
\def\openup{\afterassignment\@penup\dimen@=}
\def\@penup{\advance\lineskip\dimen@
  \advance\baselineskip\dimen@
  \advance\lineskiplimit\dimen@
  \divide\dimen@ by2
  \advance\htstrut@\dimen@
  \advance\htStrut@\dimen@
  \advance\dpstrut@\dimen@
  \advance\dpStrut@\dimen@}
\def\Let@@{\relax%
    \def\\{\global\hline@@false%
     \ifhline@\global\hline@false\global\hline@@true\fi\cr}%
    \iffalse}\fi}
\def\matrix{\null\,\vcenter\bgroup
 \tab@false\vtab@false\vst@bfalse\Strich@false\Let@@\vspace@
 \normalbaselines\openup\spreadmlines@\ialign
 \bgroup\hfil$\m@th##$\hfil&&\quad\hfil$\m@th##$\hfil\crcr
 \Mathstrut@\crcr\noalign{\kern-\baselineskip}}
\def\endmatrix{\crcr\Mathstrut@\crcr\noalign{\kern-\baselineskip}\egroup
 \egroup\,}
\def\smatrix{\D@cke0.5pt\null\,
 \vcenter\bgroup\tab@false\vtab@false\vst@bfalse\Strich@true\Let@@\vspace@
 \normalbaselines\offinterlineskip
  \openup\spreadmlines@\ialign
 \bgroup\hfil$\m@th##$\hfil&&\quad\hfil$\m@th##$\hfil\crcr}
\def\endsmatrix{\crcr\egroup
 \egroup\,\Strich@false}
\newdimen\D@cke
\def\Dicke#1{\global\D@cke#1}
\newtoks\tabs@\tabs@{&}
\newif\ifStrich@\Strich@false
\newif\iff@rst

\def\Stricherr@{\iftab@\ifvtab@\errmessage{\noexpand\s not allowed
     here. Use \noexpand\vstab!}%
  \else\errmessage{\noexpand\s not allowed here. Use \noexpand\stab!}%
  \fi\else\errmessage{\noexpand\s not allowed
     here. Use \noexpand\smatrix!}\fi}
\def\format{\ifvst@b\else\crcr\fi\egroup\iffalse{\fi\ifnum`}=0 \fi\format@}
\def\format@#1\\{\def\preamble@{#1}%
 \def\Str@chfehlt##1{\ifx##1\s\Stricherr@\fi\ifx##1\\\let\Next\relax%
   \else\let\Next\Str@chfehlt\fi\Next}%
 \def\c{\hfil\noexpand\ifhline@@\hbox{\vrule height\htStrut@%
   depth\dpstrut@ width\z@}\noexpand\fi%
   \ifStrich@\hbox{\vrule height\htstrut@ depth\dpstrut@ width\z@}%
   \fi\iftab@\else$\m@th\fi\the\hashtoks@\iftab@\else$\fi\hfil}%
 \def\r{\hfil\noexpand\ifhline@@\hbox{\vrule height\htStrut@%
   depth\dpstrut@ width\z@}\noexpand\fi%
   \ifStrich@\hbox{\vrule height\htstrut@ depth\dpstrut@ width\z@}%
   \fi\iftab@\else$\m@th\fi\the\hashtoks@\iftab@\else$\fi}%
 \def\l{\noexpand\ifhline@@\hbox{\vrule height\htStrut@%
   depth\dpstrut@ width\z@}\noexpand\fi%
   \ifStrich@\hbox{\vrule height\htstrut@ depth\dpstrut@ width\z@}%
   \fi\iftab@\else$\m@th\fi\the\hashtoks@\iftab@\else$\fi\hfil}%
 \def\s{\ifStrich@\ \the\tabs@\vrule width\D@cke\the\hashtoks@%
          \fi\the\tabs@\ }%
 \def\sa{\ifStrich@\vrule width\D@cke\the\hashtoks@%
            \the\tabs@\ %
            \fi}%
 \def\se{\ifStrich@\ \the\tabs@\vrule width\D@cke\the\hashtoks@\fi}%
 \def\cd{\hfil\noexpand\ifhline@@\hbox{\vrule height\htStrut@%
   depth\dpstrut@ width\z@}\noexpand\fi%
   \ifStrich@\hbox{\vrule height\htstrut@ depth\dpstrut@ width\z@}%
   \fi$\dsize\m@th\the\hashtoks@$\hfil}%
 \def\rd{\hfil\noexpand\ifhline@@\hbox{\vrule height\htStrut@%
   depth\dpstrut@ width\z@}\noexpand\fi%
   \ifStrich@\hbox{\vrule height\htstrut@ depth\dpstrut@ width\z@}%
   \fi$\dsize\m@th\the\hashtoks@$}%
 \def\ld{\noexpand\ifhline@@\hbox{\vrule height\htStrut@%
   depth\dpstrut@ width\z@}\noexpand\fi%
   \ifStrich@\hbox{\vrule height\htstrut@ depth\dpstrut@ width\z@}%
   \fi$\dsize\m@th\the\hashtoks@$\hfil}%
 \ifStrich@\else\Str@chfehlt#1\\\fi%
 \setbox\z@\hbox{\xdef\Preamble@{\preamble@}}\ifnum`{=0 \fi\iffalse}\fi
 \ialign\bgroup\span\Preamble@\crcr}
\newif\ifhline@\hline@false
\newif\ifhline@@\hline@@false
\def\hlinefor#1{\multispan@{\strip@#1 }\leaders\hrule height\D@cke\hfill%
    \global\hline@true\ignorespaces}
\def\Item "#1"{\par\noindent\hangindent2\parindent%
  \hangafter1\setbox0\hbox{\rm#1\enspace}\ifdim\wd0>2\parindent%
  \box0\else\hbox to 2\parindent{\rm#1\hfil}\fi\ignorespaces}
\def\ITEM #1"#2"{\par\noindent\hangafter1\hangindent#1%
  \setbox0\hbox{\rm#2\enspace}\ifdim\wd0>#1%
  \box0\else\hbox to 0pt{\rm#2\hss}\hskip#1\fi\ignorespaces}
\def\item"#1"{\par\noindent\hang%
  \setbox0=\hbox{\rm#1\enspace}\ifdim\wd0>\the\parindent%
  \box0\else\hbox to \parindent{\rm#1\hfil}\enspace\fi\ignorespaces}
\let\plainitem@\item
\catcode`\@=13

\magnification1200

\hsize13cm
\vsize19cm

\TagsOnRight

\def\bT{{\bar T}}
\def\tF{{\tilde F}}
\def\x{{\bold x}}
\def\la{{\lambda}}
\def\LR{\operatorname{LR}}
\def\bLR{\operatorname{\overline {\text {LR}}}}
\def\tLR{\operatorname{\widetilde {LR}}}
\def\con{\operatorname{con}}
\def\GL{\operatorname{GL}}
\def\SSp{\operatorname{SSp}}
\def\Sp{\operatorname{Sp}}
\def\sp{\operatorname{\text {\it sp}}}
\def\O{\operatorname{O}}
\def\o{\operatorname{\text {\it o}}}
\def\SO{\operatorname{SO}}
\def\so{\operatorname{\text {\it so}}}
\def\oddrows{\operatorname{oddrows}}
\def\oddcols{\operatorname{oddcols}}
\def\evencols{\operatorname{evencols}}

\topmatter 
\title Identities for classical group characters of nearly
rectangular shape
\endtitle 
\author C.~Krattenthaler\footnote"$^\dagger$"{Supported in part by EC's Human
Capital and Mobility Program, grant CHRX-CT93-0400 and the\linebreak 
\hbox{Austrian Science Foundation FWF, grant P10191-MAT}}
\endauthor 
\affil 
Institut f\"ur Mathematik der Universit\"at Wien,\\
Strudlhofgasse 4, A-1090 Wien, Austria.\\
e-mail: KRATT\@Pap.Univie.Ac.At\\
WWW: \tt http://radon.mat.univie.ac.at/People/kratt
\endaffil
\address Institut f\"ur Mathematik der Universit\"at Wien,
Strudlhofgasse 4, A-1090 Wien, Austria.
\endaddress
\abstract We derive several identities that feature irreducible
characters of the general linear, the symplectic, the orthogonal, and
the special orthogonal groups. All the identities feature characters
that are indexed by shapes that are ``nearly" rectangular, by which
we mean that the shapes are rectangles except for one row or column
that might be shorter than the others. As applications we prove new
results in plane partitions and tableaux enumeration, including new
refinements of the Bender-Knuth and MacMahon (ex-)conjectures.
\endabstract
\endtopmatter
\leftheadtext{C. Krattenthaler}
\rightheadtext{Characters of nearly rectangular shape}
\document

\subhead 1. Introduction\endsubhead
We prove identities that set into relation irreducible characters of
classical Lie groups. All of them feature a classical group character
indexed by a ``nearly" rectangular shape.
What is remarkable about these identities is that they have 
simple explicit forms (as opposed to many, however more general,
identities in the literature), and that many of them exhibit 
multiplicity-free expansions, typically featuring sums of characters indexed
by shapes that have a fixed number of rows or columns of odd length.

There are three types of identities
that we consider. 

First, in Theorem~1, we express (irreducible)
characters of the general linear group 
$\GL(N)$ (these characters are also known as {\it Schur
functions}) of ``nearly" 
rectangular shape in terms of (irreducible) characters of the
symplectic group $\Sp(N)$, 
and in terms of
(irreducible) characters of the orthogonal group $\O(N)$. 
As may be expected, 
for the proofs of these identities we utilize Littlewood's
\cite{\LittAA} branching rules for restricting representations of
the general linear
groups to representations of the symplectic or orthogonal groups. 

Secondly, in Theorem~2, we
express (irreducible) characters of the symplectic group $\Sp(N)$
(where $N=2n$) 
and (irreducible) characters of the special
orthogonal group $\SO(N)$ (and its spin covering group) 
of ``nearly" rectangular shape in terms of
(irreducible) characters of the general linear group
$\GL(N)$. 
For the proofs of these identities we use the
tableaux descriptions for symplectic characters developed by
DeConcini and Procesi \cite{\DeCoAA, \DePrAA}, and for special orthogonal
characters as developed by Lakshmibai, Musili and Seshadri
\cite{\LaMSAB, \LaShAC, \LaksAB}. 

Finally, 
in Theorem~3, we express the
product of two (irreducible) characters of the symplectic group $\Sp(2n)$,
respectively of the special orthogonal group $\SO(N)$ (and its spin
covering group), of ``nearly"
rectangular shape in terms of (irreducible) characters of the same
type. 
For the proofs of these identitites, we rely on Littelmann's
extension \cite{\LitPAA} of the
Littlewood--Richardson rule for all classical group characters, which
also uses the tableaux by DeConcini, Procesi, Lakshmibai, Musili and
Seshadri. In fact, as a first ``approximation", we convert
Littelmann's rule in the special case of the product of an (irreducible)
symplectic or special orthogonal character of
{\it rectangular shape} by an {\it arbitrary} (irreducible) character
of the same type to a simpler and more
explicit form, the
coefficients in the expansion being expressed in terms of 
modified Littlewood--Richardson coefficients.
The resulting formulas are given in Proposition~1 in Section~6.

We want to emphasize that it is particularly the tableaux
descriptions by DeCon\-ci\-ni, Procesi, Lakshmibai, Musili and
Seshadri that turn out to be very useful and effective here.
The more classical rules (cf\. \cite{\KingAB, \KingAC, \KoTeAB}) involving
Littlewood--Richardson coefficients,
which would also apply to the second and
third problem type that we consider, are nice in theory but
useless in practice,
since it seems to be impossible to keep track of the cancellations
that are caused by application of modification rules. 

As will become apparent, the reason that classical group characters
indexed by ``nearly" rectangular shapes allow particularly nice
identities is because Littlewood--Richardson fillings and the above
mentioned tableaux behave special for ``nearly" rectangular shapes.
This fact was previously observed by Proctor and Stanley
\cite{\ProcAE, \ProcAD, \ProcAB} for rectangular shapes. Aside from these
three papers, the inspiration for the current paper comes from
\cite{\KratAW}. There it was discovered that symplectic characters of
``nearly" rectangular shape have a nice explicit expansion in terms
of general linear characters. This gave the idea to explore what
other identities for classical group characters of ``nearly"
rectangular shapes exist. 

Finally, we remark that Okada's paper \cite{\OkadAI} may be considered
as a precursor to this paper as it contains the specializations 
to rectangular shapes of all our identities. Okada, however, uses a
completely different approach (namely, he uses the minor summation formula of
Ishikawa and Wakayama \cite{\IsWaAA}). It must be pointed out on
the other hand, that there is one type of
identities that Okada adresses, but I am not able to address. This is
restrictions of
representations of $G(N)$ ($G$ can be $\GL$, $\Sp$, or $\SO$) to
a representation of $G(l)\times G(N-l)$. Though computational
evidence exists that the respective identities of Okada can be extended from
rectangular shapes to ``nearly" rectangular shapes as well, so far I was
not able to prove anything because of the lack of an efficient
combinatorial rule. (Again, the rules in \cite{\KingAB, \KoTeAB} are
apparently useless.)

Our paper is organized as follows. In Section~2 we recall the
definitions of irreducible general linear, symplectic, orthogonal,
and special orthogonal characters. Then, in Section~3, we state our
results. These are subsequently proved in Sections~4, 5, and 6. 
Section~7 contains applications of our identities to
plane partition theory, including new refinements of the
Bender--Knuth and MacMahon (ex-)conjectures.
Finally, in
the Appendix we provide the necessary background information, in
particular there we describe all the different types of tableaux that
we use, Littlewood--Richardson fillings, the 
Littlewood--Richardson rule, Littelmann's extension of it to other
classical groups, and Littlewood's 
branching rules for restricting representations of
the general linear
groups to representations of the symplectic or orthogonal groups.

\subhead 2. General linear characters (Schur functions), symplectic,
orthogonal, and special orthogonal characters \endsubhead
Here we recall the classical character formulas. We refer the reader
to \cite{\FuKrAA; \FuHaAA, ch.~24; 
\KingAC; \KoTeAA; \ProcAK, Appendix; \SunaAE} for surveys and more
detailed background information concerning classical group characters.

In what follows $\x=(x_1,x_2,\dots)$ will always be an infinite
sequence of indeterminates. We call a sequence
$\la=(\la_1,\la_2,\dots,\la_r)$ with
$\la_1\ge\la_2\ge\dots\ge\la_r\ge0$ a {\it partition\/} if all the
$\la_i$'s are {\it integers\/}, and a {\it half-partition\/} if all the
$\la_i$'s are {\it half-integers\/}, by which we mean numbers of the
form $k+1/2$, where $k$ is an integer. For (ordinary) partitions we
adopt the convention, that partitions that only differ by trailing
zeros are considered to be the same. The components
$\la_1,\la_2,\dots$ of $\la$ are also
called the {\it parts\/} of $\la$.
We call a sequence
$\la=(\la_1,\la_2,\dots,\la_r)$ with
$\la_1\ge\la_2\ge\dots\ge\la_{r-1}\ge\vert\la_r\vert$ an
{\it $r$-orthogonal partition\/} if all the
$\la_i$'s are {\it integers\/}, and an {\it $r$-orthogonal half-partition\/} 
if all the $\la_i$'s are {\it half-integers\/}.

Given a partition $\la=(\la_1,\la_2,\dots,\la_r)$, we define the
(formal) {\it
general linear character\/}, also called {\it Schur function\/},
$s(\la;\x)$ by (see \cite{\MacdAC, I, (3.4)})
$$s(\la;\x)=\det_{1\le i,j\le r}\big(h_{\la_i-i+j}(\x)\big),
\tag2.1$$
where $h_m(\x)=\sum _{1\le i_1\le \dots\le i_m} ^{}x_{i_1}\cdots
x_{i_m}$ denotes the {\it complete homogeneous symmetric function\/}
of degree $m$. If $r\le n$, 
$$s(\la;x_1,x_2,\dots,x_n,0,0,\dots)$$ 
is the irreducible character for $\GL(n,\C)$ indexed by $\la$ (see e.g\.
\cite{\FuHaAA, (24.10)}). Following
King \cite{\KingAC} we write
$$s_n(\la;\x):=s(\la;x_1,x_2,\dots,x_n,0,0,\dots).\tag2.2$$
Two further similar notations for Schur functions that we use are
$$s_{2n}(\la;\x^{\pm1}):=s(\la;x_1,x_1^{-1},x_2,x_2^{-1},\dots,
x_n,x_n^{-1},0,0,\dots)\tag2.3$$
and
$$s_{2n+1}(\la;\x^{\pm1}):=s(\la;x_1,x_1^{-1},x_2,x_2^{-1},\dots,
x_n,x_n^{-1},1,0,0,\dots).\tag2.4$$

Again, let $\la=(\la_1,\la_2,\dots,\la_r)$ be a partition.
Following Koike and Terada \cite{\KoTeAA, Def.~2.1.1}, we define the (formal) {\it
symplectic character\/} $\sp(\la;\x)$ by 
$$\sp(\la;\x)=\det_{1\le i,j\le r}\big(h_{\la_i-i+1}(\x)\quad \vdots\quad 
h_{\la_i-i+j}(\x)+h_{\la_i-i-j+2}(\x)\big).
\tag2.5$$
Here, the notation of the determinant means
that the first expression gives the entries of the first column and
the second the entries for the remaining columns, $j\ge2$. If $r\le n$,
$$\sp(\la;x_1,x_1^{-1},x_2,x_2^{-1}, \dots,x_n,x_n^{-1},0,0,\dots)$$
is the irreducible character for $\Sp(2n,\C)$ indexed by $\la$ (see
\cite{\FuHaAA, Prop.~24.22}). We write
$$\sp_{2n}(\la;\x^{\pm1}):=
\sp(\la;x_1,x_1^{-1},x_2,x_2^{-1},
\dots,x_n,x_n^{-1},0,0,\dots).\tag2.6$$
If $r\le n+1$,
$$\sp(\la;x_1,x_1^{-1},x_2,x_2^{-1}, \dots,x_n,x_n^{-1},1,0,\dots)$$
is a character for Proctor's odd symplectic group
$\SSp(2n+1,\C)$ (see
\cite{\ProcAF, Prop.~3.1 with $x_N=1$}). We write
$$\sp_{2n+1}(\la;\x^{\pm1}):=
\sp(\la;x_1,x_1^{-1},x_2,x_2^{-1},
\dots,x_n,x_n^{-1},1,0,0,\dots).\tag2.7$$

Following Koike and Terada \cite{\KoTeAA, Def.~2.1.1}, we define the (formal) {\it
orthogonal character\/} $\o(\la;\x)$ by 
$$\o(\la;\x)=\det_{1\le i,j\le r}\big(
h_{\la_i-i+j}(\x)-h_{\la_i-i-j}(\x)\big).
\tag2.8$$
If $r\le n$,
$$\o(\la;x_1,x_1^{-1},x_2,x_2^{-1}, \dots,x_n,x_n^{-1},0,0,\dots)$$
is the irreducible character for $\O(2n,\C)$ indexed by $\la$ (see
\cite{\FuHaAA, Ex.~24.46}), for which we write
$$\o_{2n}(\la;\x^{\pm1}):=
\o(\la;x_1,x_1^{-1},x_2,x_2^{-1},
\dots,x_n,x_n^{-1},0,0,\dots).\tag2.9$$
Similarly, if $r\le n$,
$$\o(\la;x_1,x_1^{-1},x_2,x_2^{-1}, \dots,x_n,x_n^{-1},1,0,0,\dots)$$
is the irreducible character for $\O(2n+1,\C)$ indexed by $\la$ (see
\cite{\FuHaAA, Ex.~24.46}), for which we write
$$\o_{2n+1}(\la;\x^{\pm1}):=
\o(\la;x_1,x_1^{-1},x_2,x_2^{-1},
\dots,x_n,x_n^{-1},1,0,0,\dots).\tag2.10$$

The character $\o_{2n+1}(\la;\x^{\pm1})$ is also the irreducible
character for $\SO(2n+1,\C)$ indexed by $\la$. Therefore we shall
sometimes write $\so_{2n+1}(\la;\x^{\pm1})$ for
$\o_{2n+1}(\la;\x^{\pm1})$. However, for the spin covering group of 
$\SO(2n+1,\C)$ there are
also irreducible characters indexed by half-partitions. Let
$\la=(\la_1,\dots,\la_n)$ be a half-partition, then the irreducible
character $\so_{2n+1}(\la;\x^{\pm1})$ is given by the Weyl formula
(see \cite{\FuHaAA, (24.28)})
$$\so_{2n+1}(\la;\x^{\pm1})=\frac {\det\limits_{1\le i,j\le
n}(x_j^{\la_i+n-i+1/2}-x_j^{-(\la_i+n-i+1/2)})} 
{\det\limits_{1\le i,j\le
n}(x_j^{n-i+1/2}-x_j^{-(n-i+1/2)})}.\tag2.11$$
In fact (see again \cite{\FuHaAA, (24.28)}), if $\la$ is an ordinary
partition, then $\so_{2n+1}(\la;\x^{\pm1})$ is given by the same
formula.

The situation is even more delicate for $\SO(2n,\C)$. 
The irreducible characters for $\SO(2n,\C)$ and its spin covering
group are indexed by
$n$-orthogonal partitions or half-partitions. Let
$\la=(\la_1,\la_2,\dots,\la_n)$ be an $n$-orthogonal partition or
half-partition, i.e., 
$\la_1\ge\la_2\ge\dots\ge\la_{n-1}\ge\vert\la_n\vert$, with
the $\la_i$'s integers, respectively half-integers. Then the
irreducible special orthogonal character $\so_{2n}(\la;\x^{\pm1})$ is
given by the Weyl formula (see \cite{\FuHaAA, (24.40)})
$$\so_{2n}(\la;\x^{\pm1})=\frac {\det\limits_{1\le i,j\le
n}(x_j^{\la_i+n-i}+x_j^{-(\la_i+n-i)})+
\det\limits_{1\le i,j\le
n}(x_j^{\la_i+n-i}-x_j^{-(\la_i+n-i)})} 
{\det\limits_{1\le i,j\le
n}(x_j^{n-i}+x_j^{-(n-i)})}.\tag2.12$$
The irreducible character for $\O(2n,\C)$
indexed by the partition $\la$, $\o_{2n}(\la;\x^{\pm1})$,
equals the irreducible character for $\SO(2n,\C)$
indexed by $\la$, $\so_{2n}(\la;\x^{\pm1})$,
if $\la_n=0$, but splits into two irreducible characters for
$\SO(2n,\C)$ if $\la_n\ne0$, 
one indexed by $\la=(\la_1,\la_2,\dots,\la_n)$, the
other indexed by $\la^-$, which by definition is  
$(\la_1,\la_2,\dots,\la_{n-1},\mathbreak-\la_n)$ (see \cite{\FuHaAA,
first paragraph on p.~411}).

\medskip
The characters $s(\la;\x)$, $\sp(\la;\x)$, $\o(\la;\x)$ 
are also called
{\it universal characters\/}, meaning that by specializing one
obtains the actual (general linear, symplectic, orthogonal,
respectively) characters for any dimension of the corresponding
group. There is no such thing in the even special orthogonal case
and for odd special orthogonal characters indexed by half-partitions.
Of course, it has to be mentioned that there is not only a Weyl
formula for $\so_N(\la;\x^{\pm1})$, but also for $s_n(\la;\x)$ 
(see \cite{\FuHaAA, p.~403, (A.4); \MacdAC, I, (3.1)}) and
$\sp_{2n}(\la;\x)$ (see \cite{\FuHaAA, (24.18)}), the latter reading
$$\sp_{2n}(\la;\x^{\pm1})=\frac {\det\limits_{1\le i,j\le
n}(x_j^{\la_i+n-i+1}-x_j^{-(\la_i+n-i+1)})} 
{\det\limits_{1\le i,j\le
n}(x_j^{n-i+1}-x_j^{-(n-i+1)})}.\tag2.13$$
Actually, the Weyl formulas for all these characters allow a
uniform statement (see \cite{\FuHaAA,
Theorem~24.2; \SunaAE, Theorem~4.1}).

\subhead 3. Character identities\endsubhead
In this section we collect our identities for classical group
characters of ``nearly" rectangular shapes. By ``nearly" rectangular
shape we mean partitions of the form $(c,c,\dots,c,c-p)$ or
$(c,\dots,c,c-1,\dots,c-1)$, i.e., partitions whose Ferrers diagram
(see Section~A1 of the Appendix) is a rectangle, except that
one row or column might be shorter. We write $(c^{r-1}, c-p)$ for the
first partition (given that the partition has $r$ components) and
$(c^{r-p},(c-1)^p)$ for the second (given that the partition has $r$
components, $p$ of which are equal to $c-1$).

\medskip
The first five identities express general linear characters (Schur
functions) in terms of symplectic, respectively orthogonal
characters. Most of the ``partition terminology" that we use here and
henceforth is explained in Section~A1 of the Appendix. Other
terminology concerns the parity of rows or
columns: When we say `an odd (even) row' we mean `a row of odd (even) length.'
We use the same convention with columns. For convenience, 
we denote the number of odd
rows of some (ordinary or skew) shape $\si$ by $\oddrows(\si)$ and
the number of odd columns of $\si$ by $\oddcols(\si)$, and the same
with `even' instead of `odd.'
\proclaim{Theorem 1} Let $r,c,p$ by nonnegative integers.
For $c\ge p$ there holds
$$s\big((c^{r-1},c-p);\x\big)=\underset
\oddcols\!\big((c^r)/\nu\big)=p
\to{\sum _{\nu\subseteq (c^r)} ^{}}\sp(\nu;\x),\tag3.1$$
and for $r\ge p$ there holds
$$s\big((c^{r-p},(c-1)^p);\x\big)=\underset
\oddrows\!\big((c^r)/\nu\big)=p
\to{\sum _{\nu\subseteq (c^r)} ^{}}\o(\nu;\x).\tag3.2$$
In particular, for $r\le N$ and $c\ge p$ there hold
$$s_{N}\big((c^{r-1},c-p);\x^{\pm1}\big)=\cases
\underset
\oddcols\!\big((c^r)/\nu\big)=p
\to{\sum\limits _{\nu\subseteq (c^r)} ^{}}\sp_{N}(\nu;\x^{\pm1})&r\le \lceil
N/2\rceil\\
\underset
\oddcols\!\big((c^{N-r+1})/\nu\big)=c-p
\to{\sum\limits _{\nu\subseteq (c^{N-r+1})} ^{}}\sp_{N}(\nu;\x^{\pm1})&r> \lceil
N/2\rceil,
\endcases\tag3.3$$
and
$$s_{2n+1}\big((c^{r});\x^{\pm1}\big)=\cases
{\sum\limits _{\nu\subseteq (c^r)} ^{}}\sp_{2n}(\nu;\x^{\pm1})&r\le n\\
{\sum\limits _{\nu\subseteq (c^{2n+1-r})} ^{}}\sp_{2n}(\nu;\x^{\pm1})&r> n.
\endcases\tag3.4$$
For $p\le r\le \lfloor N/2\rfloor$ there holds
$$s_N\big((c^{r-p},(c-1)^p);\x^{\pm1}\big)=\underset
\oddrows\!\big((c^r)/\nu\big)=p
\to{\sum\limits _{\nu\subseteq (c^r)} ^{}}\o_N(\nu;\x^{\pm1}).\tag3.5$$
\endproclaim

\remark{Remark} 
The identities (3.3)--(3.5) have obvious interpretations as
branching rules for the
restriction of representation modules of $\GL(N,\C)$ to $\Sp(N,\C)$,
$\Sp(N-1,\C)$, respectively $\O(N,\C)$.

The case $N=2n$, $r\le n$ of (3.3) is implicit in
Proctor's paper
\cite{\ProcAE, Proof of Lemma~4, Claim on p.~558}. In the same paper,
there appear, explicitly, the case $N=2n$, $r\le n$, $p=0$ of (3.3)
\cite{\ProcAE, Lemma~4, equation for $A_{2n-1}(m\om_r)$}
and the case $r\le n$ of (3.4)
\cite{\ProcAE, Lemma~4, equation for $A_{2n}(m\om_r)$}. The 
cases $N=2n$, $r\le n$, $p=0$ of (3.3)--(3.5) appear in \cite{\OkadAI,
Theorem~2.6}. 
\endremark

\medskip
The next four identities express symplectic and special orthogonal
characters in terms of general linear characters (Schur functions).
\proclaim{Theorem 2}Let $n,c,p$ by nonnegative integers. For $n\ge
p$ there holds
$$\sp_{2n}\big((c^{n-p},(c-1)^p);\x^{\pm1}\big)=
(x_1x_2\cdots x_n)^{-c}\cdot\underset \oddrows(\nu)=p
\to{\sum _{\nu\subseteq ((2c)^n)} ^{}}s_n(\nu;\x).\tag3.6$$

If $c$ is a nonnegative integer or half-integer and $p$ a
nonnegative integer, $2c\ge p$, then there holds
$$\so_{2n}\big((c^{n-1},c-p);\x^{\pm1}\big)=
(x_1x_2\cdots x_n)^{-c}\cdot\underset \oddcols\!\big(((2c)^n)/\nu\big)=p
\to{\sum _{\nu\subseteq ((2c)^n)} ^{}}s_n(\nu;\x).\tag3.7$$

Next, if $c$ is a nonnegative integer or half-integer and $p$ a
nonnegative integer, $n\ge p$, then there holds
$$\so_{2n+1}\big((c^{n-p},(c-1)^p);\x^{\pm1}\big)=
(x_1x_2\cdots x_n)^{-c}\cdot
{\sum _{\nu\subseteq ((2c)^{n})} ^{}}
a_{n,p}(\nu)\cdot s_n(\nu;\x),\tag3.8$$
where 
$$\align a_{n,p}(\nu)&=\text {number of vertical strips of length $p$
on the}\\
&\hphantom{{}={}}\text {rim of $\nu$ avoiding the $(2c)$-th column.}
\tag3.9\endalign$$

Finally, if $c$ is a nonnegative integer or half-integer and $p$ a
nonnegative integer, $c\ge p$, then there holds
$$\so_{2n+1}\big((c^{n-1},c-p);\x^{\pm1}\big)=
(x_1x_2\cdots x_n)^{-c}\cdot
{\sum _{\nu\subseteq ((2c)^{n})} ^{}}
b_{n,p}(\nu)\cdot s_n(\nu;\x),\tag3.10$$
where 
$$\align b_{n,p}(\nu)&=\text {number of horizontal strips of length $p$
on the rim of $\nu$}\\
&\hphantom{{}={}}\text {such that the $i$-th cell of the
strip (counted from left to}\\
&\hphantom{{}={}} \text {right) comes before the $(2c-2p+2i)$-th column.}
\tag3.11\endalign$$

\endproclaim
\remark{Remark} 
The identities (3.6)--(3.8), and (3.10) have obvious interpretations as
decomposition formulas
for representation modules of $\Sp(2n,\C)$,
$\SO(2n,\C)$, or $\SO(2n+1,\C)$ as representations of the subgroup
$\GL(n,\C)$. For the interested reader we add that, without much
additional effort, it is also possible to derive decomposition
formulas for $\sp_{2n}\big((c^{n-1},c-p);\x^{\pm1}\big)$ and 
$\so_{2n}\big((c^{n-p},(c-1)^p);\x^{\pm1}\big)$ which are in the
spirit of (3.10), respectively (3.8), 
by slightly modifying the arguments in the
proofs of (3.10) and (3.8).

Formula (3.6) for the first time appeared in
\cite{\KratAW}, although it is implicit already in \cite{\GoulAD,
Theorem~2.6; \KratAQ, (2.2)}. The case $p=0$ of (3.6) appeared
already in a number of papers \cite{\ProcAB, Theorem~3; 
\StemAD, Theorem~4.1; \StemAE, Cor.~7.4.(b); \OkadAI, Theorem~2.3.(2)}. 

The case $p=0$ of (3.8) and (3.10), which is
$$\so_{2n+1}\big((c^{n});\x^{\pm1}\big)=
(x_1x_2\cdots x_n)^{-c}\cdot
{\sum _{\nu\subseteq ((2c)^{n})} ^{}}
 s_n(\nu;\x),\tag3.12$$
appeared also in a number of papers \cite{\ProcAB, Theorem~3; 
\StemAE, Cor.~7.4.(a); \OkadAI, Theorem~2.3.(1)}. 

Finally, the cases $p=0$ and $p=c$ of (3.7) previously appeared in
\cite{\BrGrAA; \OkadAI, Theorem~2.3.(3)}.

\endremark

\medskip
Our final identities of this section display the decomposition of the
product of a rectangularly shaped and a ``nearly" rectangularly shaped
(general linear, symplectic, respectively special orthogonal) character.
As will become apparent from the proofs, 
these identities follow rather easily from decomposition
formulas for the product of a rectangularly shaped symplectic,
respectively special orthogonal, character and an {\it arbitrarily\/}
shaped symplectic, respectively special orthogonal, character,
which we obtain in Proposition~1 in Section~6 from
Littelmann's decomposition formula \cite{\LitPAA}.

\proclaim{Theorem 3} 
Let $n,c,d,p$ be nonnegative integers. 

First let $n\ge p$.
If $c\le d$ then there holds
$$
\sp_{2n}\big((c^n);\x^{\pm1}\big)\cdot
\sp_{2n}\big(((d+1)^p,d^{n-p});\x^{\pm1}\big)=
\underset \oddrows\!\big(\nu/((d-c)^n)\big)=p
\to{\sum _{\((d-c)^n\)\subseteq\nu\subseteq\((c+d+1)^n\)} ^{}}\sp_{2n}(\nu;\x^{\pm1}),
\tag3.13$$
and if $c\ge d$,
$$
\sp_{2n}\big((c^n);\x^{\pm1}\big)\cdot
\sp_{2n}\big((d^{n-p},(d-1)^{p});\x^{\pm1}\big)=
\underset \oddrows\!\big(\nu/((c-d)^n)\big)=p
\to{\sum _{\((c-d)^n\)\subseteq\nu\subseteq\((c+d)^n\)} ^{}}\sp_{2n}(\nu;\x^{\pm1}),
\tag3.14$$

Next, if $c,d$ are nonnegative integers or half-integers and $p$ an
integer with 
$p\le 2d$, then there hold
$$
\so_{2n}\big((c^n);\x^{\pm1}\big)\cdot
\so_{2n}\big((d^{n-1},d-p);\x^{\pm1}\big)=
\underset \oddcols\!\big(((c+d)^n)/\nu\big)=p
\to{\sum _{\(\vert c-d\vert^{n-1},c-d\)\subseteq\nu\subseteq\((c+d)^n\)} ^{}}\so_{2n}(\nu;\x^{\pm1}),
\tag3.15$$
with the understanding that $\nu$ ranges over $n$-orthogonal 
partitions if $c+d$ is
an integer and over $n$-orthogonal half-partitions if $c+d$ is a half-integer, 
and
$$\multline
\so_{2n}\big((c^{n-1},-c);\x^{\pm1}\big)\cdot
\so_{2n}\big((d^{n-1},d-p);\x^{\pm1}\big)
\\=
\underset \evencols\!\big(((c+d)^n)/\nu\big)=p
\to{\sum _{\(\vert c-d\vert^{n-1},c-d\)\subseteq\nu\subseteq\((c+d)^n\)} ^{}}
\so_{2n}((\nu_1,\dots,\nu_{n-1},-\nu_n);\x^{\pm1}),
\endmultline\tag3.16$$
with the same understanding.

For $c,d$ nonnegative integers or half-integers and for
$n>p$ there holds
$$\multline
\so_{2n+1}\big((c^n);\x^{\pm1}\big)\cdot
\so_{2n+1}\big((d^{n-p},(d-1)^p);\x^{\pm1}\big)\\
={\sum _{\big(\vert c-d\vert^{n-p},(\max\{c-d,d-c-1\})^{p}\big)
\subseteq\nu\subseteq\big((c+d)^n\big)} ^{}}
c_{n,p}(\nu)\cdot \so_{2n+1}(\nu;\x^{\pm1}),
\endmultline\tag3.17$$
where, 
with $m_l(\nu)$ denoting the multiplicity of $l$ in $\nu$,
$$\align c_{n,p}(\nu)&=\text {number of vertical strips of length
$p-m_{d-c-1}(\nu)$}\\
&\hphantom{{}={}}\text {that can be added to $\nu$ to obtain another
(half-)partition, avoiding the}\\
&\hphantom{{}={}}\text {$(d-c)$-th, the $(c-d+1)$-st, and the
$(c+d+1)$-st
column of $\nu$.}
\tag3.18\endalign$$
Again, the sum in (3.17) is understood to range over partitions if
$c+d$ is an integer and over half-partitions if $c+d$ is a
half-integer.

Finally, for $c,d$ nonnegative integers or half-integers and for
$d\ge p$ there holds
$$\multline
\so_{2n+1}\big((c^n);\x^{\pm1}\big)\cdot
\so_{2n+1}\big((d^{n-1},d-p);\x^{\pm1}\big)\\
={\sum _{\big(\vert c-d\vert^{n-1},\max\{c-d,d-c-p\}\big)
\subseteq\nu\subseteq\big((c+d)^n\big)} ^{}}
d_{n,p}(\nu)\cdot \so_{2n+1}(\nu;\x^{\pm1}),
\endmultline\tag3.19$$
where $d_{n,p}(\nu)$ is the number of all horizontal strips $\si$, 
that can be added to $\nu$ to obtain another
(half-)partition, avoiding the $(c+d+1)$-st column and the $n$-th
row of $\nu$, and which satisfy the following inequalities: If $\si_i$
denotes the number of cells of the strip $\si$ in the $i$-th row,
$i=1,2,\dots,n-1$, and if $\vert\si\vert$ denotes the total number of
cells in $\si$, then
$$\spreadmatrixlines{2\jot}
\matrix  \big\vert\vert\si\vert - c+d-p\big\vert\le \nu_n,\quad
\nu_n+p-\nu_{n-1}\le \vert\si\vert\le p\\
\text {and}\quad 2\si_1+\dots+2\si_{i-1}+\si_i+\nu_i\ge c-d+2p\quad \text
{for }i=1,2,\dots,n-1.
\endmatrix
\tag3.20$$
Also here, the sum in (3.19) is understood to range over partitions if
$c+d$ is an integer and over half-partitions if $c+d$ is a
half-integer.
\endproclaim

\remark{Remark} All these
formulas, that is (3.13)--(3.20) and (6.1)--(6.5),
have obvious interpretations as decomposition formulas
for the tensor product of two representation modules of $\Sp(2n,\C)$,
$\SO(2n,\C)$, or $\SO(2n+1,\C)$.
For the interested reader we add that, without much
additional effort, it is also possible to derive decomposition
formulas for $\sp_{2n}\big((c^{n});\x^{\pm1}\big)\cdot
\sp_{2n}\big((d^{n-1},d-p);\x^{\pm1}\big)$ and 
$\so_{2n}\big((c^{n});\x^{\pm1}\big)\cdot
\so_{2n}\big((d^{n-p},(d-1)^p);\x^{\pm1}\big)$ which are in the
spirit of (3.19), respectively (3.17), 
by slightly modifying the arguments in the
proofs of (3.19) and (3.17).

Formulas (3.13)--(3.20) generalize the symplectic and special
orthogonal decompositions of Okada
\cite{\OkadAI, Theorem~2.5}, who proved the $p=0$ special cases.
Okada also proves a decomposition formula for the product of two
rectangularly shaped general linear characters (Schur functions), thus
generalizing Stanley's results \cite{\StanAI, Lemma~3.3}. More
generally, Carini \cite{\CariAA, sec.~3.3, p.~105ff} derived decomposition formulas
for the product of two ``nearly" rectangularly shaped general linear characters.
\endremark

\subhead 4. Proof of Theorem~1\endsubhead
Here we use decomposition rules of Littlewood \cite{\LittAA} (see 
Section~A7 of the Appendix). 
All the notions that appear in this section, like
Littlewood--Richardson filling  (LR-filling), Littlewood--Richardson
condition (LR-condition), content, etc., are also explained in the
Appendix, mainly in Section~A6.

\demo {Proof of {\rm (3.1)}} Implicitly, this was already proved in
\cite{\ProcAE, Proof of Lemma~4, Claim on p.~558}. However, since
there is no explicit statement in \cite{\ProcAE} and since one would
have to translate things appropriately, we include a detailed proof of (3.1)
here.

According to (A.13) we have
$$s\big((c^{r-1},c-p);\x\big)=
{\sum _{\nu} ^{}}\sp(\nu;\x)\sum _{\mu,\ \mu'\text { even}}
^{}\LR_{\mu,\nu}^{(c^{r-1},c-p)}.$$
Hence we have to show that 
$$\sum _{\mu,\ \mu'\text { even}}
^{}\LR_{\mu,\nu}^{(c^{r-1},c-p)}=\cases
1&\oddcols\big((c^r)/\nu\big)=p\\0&\text {otherwise}.\endcases\tag4.1$$

To see this, suppose that $\mu$ is a fixed partition whose Ferrers
diagram has only even columns. Consider a LR-filling of
shape $(c^{r-1},c-p)/\mu$. By the LR-condition, there is no choice
for the entries in the first $r-1$ rows. I.e., all the rows except
for the $r$-th row are uniquely determined in the way that is
exemplified in Figure~1. There is only freedom in the $r$-th row.
\vskip10pt
\vbox{
$$
\Einheit0.3cm
\Pfad(0,0),11111111111121111122222222\endPfad
\Pfad(0,0),112111221222211122\endPfad
\Pfad(0,0),22222222211111111111111111\endPfad
\PfadDicke{.5pt}
\Pfad(5,1),1111111\endPfad
\Pfad(5,2),111111111111\endPfad
\Pfad(6,3),11111111111\endPfad
\Pfad(6,4),11111111111\endPfad
\Pfad(6,5),11111111111\endPfad
\Pfad(6,6),11111111111\endPfad
\Pfad(9,7),11111111\endPfad
\Pfad(9,8),11111111\endPfad
\Pfad(3,0),2\endPfad
\Pfad(4,0),2\endPfad
\Pfad(5,0),2\endPfad
\Pfad(6,0),222\endPfad
\Pfad(7,0),2222222\endPfad
\Pfad(8,0),2222222\endPfad
\Pfad(9,0),2222222\endPfad
\Pfad(10,0),222222222\endPfad
\Pfad(11,0),222222222\endPfad
\Pfad(12,1),22222222\endPfad
\Pfad(13,1),22222222\endPfad
\Pfad(14,1),22222222\endPfad
\Pfad(15,1),22222222\endPfad
\Pfad(16,1),22222222\endPfad
\Label\o{\mu}(3,5)
\Label\ro{\ssize 1}(9,8)
\Label\ro{\ssize 1}(10,8)
\Label\ro{\ssize 1}(11,8)
\Label\ro{\ssize 1}(12,8)
\Label\ro{\ssize 1}(13,8)
\Label\ro{\ssize 1}(14,8)
\Label\ro{\ssize 1}(15,8)
\Label\ro{\ssize 1}(16,8)
\Label\ro{\ssize 2}(9,7)
\Label\ro{\ssize 2}(10,7)
\Label\ro{\ssize 2}(11,7)
\Label\ro{\ssize 2}(12,7)
\Label\ro{\ssize 2}(13,7)
\Label\ro{\ssize 2}(14,7)
\Label\ro{\ssize 2}(15,7)
\Label\ro{\ssize 2}(16,7)
\Label\ro{\ssize 1}(6,6)
\Label\ro{\ssize 1}(7,6)
\Label\ro{\ssize 1}(8,6)
\Label\ro{\ssize 3}(9,6)
\Label\ro{\ssize 3}(10,6)
\Label\ro{\ssize 3}(11,6)
\Label\ro{\ssize 3}(12,6)
\Label\ro{\ssize 3}(13,6)
\Label\ro{\ssize 3}(14,6)
\Label\ro{\ssize 3}(15,6)
\Label\ro{\ssize 3}(16,6)
\Label\ro{\ssize 2}(6,5)
\Label\ro{\ssize 2}(7,5)
\Label\ro{\ssize 2}(8,5)
\Label\ro{\ssize 4}(9,5)
\Label\ro{\ssize 4}(10,5)
\Label\ro{\ssize 4}(11,5)
\Label\ro{\ssize 4}(12,5)
\Label\ro{\ssize 4}(13,5)
\Label\ro{\ssize 4}(14,5)
\Label\ro{\ssize 4}(15,5)
\Label\ro{\ssize 4}(16,5)
\Label\ro{\ssize 3}(6,4)
\Label\ro{\ssize 3}(7,4)
\Label\ro{\ssize 3}(8,4)
\Label\ro{\ssize 5}(9,4)
\Label\ro{\ssize 5}(10,4)
\Label\ro{\ssize 5}(11,4)
\Label\ro{\ssize 5}(12,4)
\Label\ro{\ssize 5}(13,4)
\Label\ro{\ssize 5}(14,4)
\Label\ro{\ssize 5}(15,4)
\Label\ro{\ssize 5}(16,4)
\Label\ro{\ssize 4}(6,3)
\Label\ro{\ssize 4}(7,3)
\Label\ro{\ssize 4}(8,3)
\Label\ro{\ssize 6}(9,3)
\Label\ro{\ssize 6}(10,3)
\Label\ro{\ssize 6}(11,3)
\Label\ro{\ssize 6}(12,3)
\Label\ro{\ssize 6}(13,3)
\Label\ro{\ssize 6}(14,3)
\Label\ro{\ssize 6}(15,3)
\Label\ro{\ssize 6}(16,3)
\Label\ro{\ssize 1}(5,2)
\Label\ro{\ssize 5}(6,2)
\Label\ro{\ssize 5}(7,2)
\Label\ro{\ssize 5}(8,2)
\Label\ro{\ssize 7}(9,2)
\Label\ro{\ssize 7}(10,2)
\Label\ro{\ssize 7}(11,2)
\Label\ro{\ssize 7}(12,2)
\Label\ro{\ssize 7}(13,2)
\Label\ro{\ssize 7}(14,2)
\Label\ro{\ssize 7}(15,2)
\Label\ro{\ssize 7}(16,2)
\Label\ro{\ssize 2}(5,1)
\Label\ro{\ssize 6}(6,1)
\Label\ro{\ssize 6}(7,1)
\Label\ro{\ssize 6}(8,1)
\Label\ro{\ssize 8}(9,1)
\Label\ro{\ssize 8}(10,1)
\Label\ro{\ssize 8}(11,1)
\Label\ro{\ssize 8}(12,1)
\Label\ro{\ssize 8}(13,1)
\Label\ro{\ssize 8}(14,1)
\Label\ro{\ssize 8}(15,1)
\Label\ro{\ssize 8}(16,1)
\Label\ro{\ssize *}(2,0)
\Label\ro{\ssize *}(3,0)
\Label\ro{\ssize *}(4,0)
\Label\ro{\ssize *}(5,0)
\Label\ro{\ssize *}(6,0)
\Label\ro{\ssize *}(7,0)
\Label\ro{\ssize *}(8,0)
\Label\ro{\ssize *}(9,0)
\Label\ro{\ssize *}(10,0)
\Label\ro{\ssize *}(11,0)
\Label\o{\hskip0.3cm\oversetbrace \tsize c\to{\hbox{\hskip5cm}}}(8,9)
\Label\l{\hskip0cm\raise1.6cm\hbox{$\sideset r\and \to
     {\left\{\vbox{\vskip1.5cm}\right.}$}}(0,4)
\Label\u{\hskip0.3cm\undersetbrace \tsize p\to{\hbox{\hskip1.5cm}}}(14,1)
\hskip5.1cm
$$
\centerline{\eightpoint Figure 1}
}
\vskip10pt

It should be observed that the content of the uniquely determined
part of the LR-filling equals
$\tilde\mu:=(c-\mu_{r-1},c-\mu_{r-2},\dots,c-\mu_1)$. In particular,
this implies that all the columns in $(c^r)/\tilde\mu$ have odd
length, except for columns $c-\mu_{r-1}+1,c-\mu_{r-1}+2,\dots,c$
(which would have length $r$) in case that $r$ is even. 

Suppose that $\nu$ is the content of the {\it complete\/} LR-filling.
Then the LR-condition is equivalent to saying that $\nu/\tilde\mu$ is
a horizontal strip (see Appendix~A1 for the definition of a 
horizontal strip). Besides, this horizontal strip is of length
$c-p-\mu_r$, the length of the $r$-th row of the LR-filling. Hence,
there are $c-(c-p-\mu_r)=p+\mu_r$ odd columns in $(c^r)/\nu$
if $c$ is odd, and there are
$c-\mu_{r-1}-(c-p-\mu_r)=p+\mu_r-\mu_{r-1}$ odd columns 
in $(c^r)/\nu$ if $r$ is even. The latter is due to the fact that 1
cannot be an entry in the $r$-th row of the LR-filling if $r$ is
even, hence the horizontal strip has to avoid the first row.
Actually, the number of odd
columns is $p$ in both cases, since in case $r$ odd we must have
$\mu_r=0$ because the Ferrers diagram of $\mu$ has only even columns,
and in case $r$ even we have $\mu_{r-1}=\mu_r$ for the same reason.
Thus we have shown that, given
a fixed $\mu$ whose Ferrers diagram has only even columns,
then $\LR_{\mu,\nu}^{(c^{r-1},c-p)}\ne0$ only if $\nu$ is a partition
such that the number of odd columns in $(c^r)/\nu$ equals $p$. In
particular, this establishes the ``otherwise" part of (4.1).

Conversely, given a partition $\nu\subseteq(c^r)$ 
such that the number of odd columns in $(c^r)/\nu$ equals $p$, we
claim that there is exactly one partition $\mu$
whose Ferrers diagram has only
even columns and such that $\LR_{\mu,\nu}^{(c^{r-1},c-p)}\ne0$. And, more
precisely, we have $\LR_{\mu,\nu}^{(c^{r-1},c-p)}=1$, which means that
there is exactly one LR-filling of shape $(c^{r-1},c-p)/\mu$ and content
$\nu$. Altogether, this would establish (4.1).

The claim is established by going through the preceding paragraphs,
backwards. Suppose that $\nu\subseteq(c^r)$ is a partition with 
$\LR_{\mu,\nu}^{(c^{r-1},c-p)}\ne0$, for some $\mu$ whose Ferrers
diagram has only even columns. Then $\tilde\mu$, defined as above as 
$(c-\mu_{r-1},c-\mu_{r-2},\dots,c-\mu_1)$, is a partition in which
(1) all the columns have length $\equiv r$ mod~2, and (2) which is
contained in $\nu$ and differs from
$\nu$ by a horizontal strip of length $c-p-\mu_r$. It is
straight-forward to see that (1) and (2) determine $\tilde\mu$
uniquely. Hence, $\mu_1,\mu_2,\dots,\mu_{r-1}$ are uniquely
determined, thus also $\mu_r$, since the Ferrers diagram of $\mu$ has
to contain only even columns. To be precise, the latter condition implies that 
$\mu_r$ equals $\mu_{r-1}$ if $r$ is even, and 0 if $r$ is odd. 

Let $\mu$ be this uniquely determined partition. To show that
$\LR_{\mu,\nu}^{(c^{r-1},c-p)}=1$ it remains to see that there is
exactly one LR-filling of shape $(c^{r-1},c-p)/\mu$ with content
$\nu$. In fact, as we already observed, the entries of the first $r-1$
rows of a LR-filling of shape $(c^{r-1},c-p)/\mu$ are uniquely
determined. The content of this partial filling of these first $r-1$
rows is $\tilde\mu$. Since $\nu$, which should be the content of
the complete LR-filling, differs from $\tilde\mu$ by a horizontal
strip, also the entries of the $r$-th row are uniquely determined. To
be precise, the length of the $i$-th row of $\nu/\tilde\mu$ gives the
multiplicity of $i$ in the $r$-th row of the LR-filling.

Altogether, this establishes (4.1), and hence (3.1), as desired.\quad \quad
\qed
\enddemo

\demo {Proof of {\rm(3.2)}}
According to (A.14) we have
$$s\big((c^{r-p},(c-1)^p);\x\big)=
{\sum _{\nu} ^{}}\o(\nu;\x)\sum _{\mu,\ \mu\text { even}}
^{}\LR_{\mu,\nu}^{(c^{r-p},(c-1)^p)}.$$
Hence we have to show that 
$$\sum _{\mu,\ \mu\text { even}}
^{}\LR_{\mu,\nu}^{(c^{r-p},(c-1)^p)}=\cases
1&\oddrows\!\big((c^r)/\nu\big)=p\\0&\text {otherwise}.\endcases\tag4.2$$
We could prove this again directly, in a similar style as in the proof
of (4.1). However, once (4.1) is already known, the companion (4.2)
follows straight-forwardly from the well-known identity (see
\cite{\HaSuAA}) $\LR_{\mu,\nu}^{\la}=\LR_{\mu',\nu'}^{\la'}$.\quad \quad
\qed
\enddemo

For the proofs of (3.3) and (3.4) we utilize the following auxiliary result.

\proclaim{Lemma}
Let $N$ be a positive integer. Then, for any partition
$\la=(\la_1,\la_2,\dots,\la_{N})$ we have
$$s_{N}\big((\la_1,\la_2,\dots,\la_{N});\x^{\pm1}\big)=
s_{N}\big((\la_1-\la_{N},\la_1-\la_{N-1},\dots,
\la_{1}-\la_2,0);\x^{\pm1}\big).
\tag4.3$$
\endproclaim

\demo{Proof}
This identity follows upon little manipulation from the dual
Jacobi--Trudi identity (the N\"agelsbach--Kostka identity, see
\cite{\MacdAC, I, (3.5)}) for Schur functions,
$$s(\la;\x)=\det_{1\le i,j\le \la_1}\big(e_{\la'_i-i+j}(\x)\big),
\tag4.4$$
where $e_m(\x)=\sum _{1\le i_1< \dots< i_m} ^{}x_{i_1}\cdots
x_{i_m}$ denotes the {\it elementary symmetric function\/}
of degree $m$. For, by (4.4), and because of
$$e_m(x_1,x_1^{-1},\dots,x_n,x_n^{-1})=e_{2n-m}(x_1,x_1^{-1},\dots,
x_n,x_n^{-1})$$
and a similar identity for $e_m(x_1,x_1^{-1},\dots,x_n,x_n^{-1},1)$,
we have (in the following calculation $e_m(\x^{\pm1})$ is short for 
$e_m(x_1,x_1^{-1},\dots,x_n,x_n^{-1})$, respectively for 
$e_m(x_1,x_1^{-1},\dots,x_n,\mathbreak
x_n^{-1},1)$, depending on whether $N$ is
even or odd)
$$\align s_{N}\big((\la_1,\la_2,\dots,\la_{N});\x^{\pm1}\big)&=
\det_{1\le i,j\le \la_1}\big(e_{\la'_i-i+j}(\x^{\pm1})\big)\\
&=\det_{1\le i,j\le \la_1}\big(e_{N-\la'_i+i-j}(\x^{\pm1})\big)\\
&=\det_{1\le i,j\le \la_1}\big(e_{N-\la'_{\la_1+1-i}+(\la_1+1-i)-
(\la_1+1-j)}(\x^{\pm1})\big)\\
&=s_{N}\big((\la_1-\la_{N},\la_1-\la_{N-1},\dots,
\la_{1}-\la_2,0);\x^{\pm1}\big).
\endalign$$

\enddemo

\demo {Proof of {\rm(3.3)}}In (3.1) we specialize $\x$ to
$(x_1,x_1^{-1},\dots,x_n,x_n^{-1},0,0,\dots)$ if $N=2n$, and
to $(x_1,x_1^{-1},\dots,x_n,x_n^{-1},1,0,0,\dots)$ if $N=2n+1$, which leads to
$$s_{N}\big((c^{r-1},c-p);\x^{\pm1}\big)=\underset
\oddcols\!\big((c^r)/\nu\big)=p
\to{\sum _{\nu\subseteq (c^r)} ^{}}\sp_N(\nu;\x^{\pm1}).\tag4.5$$
This gives (3.3) for $r\le \lceil N/2\rceil$ immediately. In case
$r>\lceil N/2\rceil$, on the right-hand side of (4.5) the partition
$\nu$ could have more than $\lceil N/2\rceil$ parts, in which case
one would have to apply modification rules for the corresponding
symplectic characters $\sp_N(\nu;\x)$ (see e.g\. \cite{\KingAA;
\SunaAE, Theorem~5.4}). However, we circumvent this difficulty by
appealing to the Lemma above.
In fact, by (4.3) and (4.5) we have for $r>\lceil N/2\rceil$,
$$\align
s_{N}\big((c^{r-1},c-p);\x^{\pm1}\big)&=s_{N}\big((c^{N-r},p);
\x^{\pm1}\big)\\
&=\underset
\oddcols\!\big((c^{N-r+1})/\nu\big)=c-p
\to{\sum _{\nu\subseteq (c^{N-r+1})} ^{}}\sp_N(\nu;\x^{\pm1}),
\endalign$$
which is (3.3) for $r>\lceil N/2\rceil$.\quad \quad
\qed
\enddemo

\demo {Proof of {\rm(3.4)}}
In (3.1) we set $N=2n+1$, $p=0$, and 
substitute $(x_1,x_1^{-1},\dots,x_n,\mathbreak
x_n^{-1},1,0,0,\dots)$
for $\x$ to obtain
$$s_{2n+1}\big((c^r);\x^{\pm1}\big)=\underset
\oddcols\!\big((c^r)/\rh\big)=0
\to{\sum _{\rh\subseteq (c^r)} ^{}}\sp_{2n+1}(\rh;\x^{\pm1}).\tag4.6$$

Now let first be $r\le n$. The odd
symplectic characters on the right-hand side of (4.6) are known to
decompose in terms of even symplectic characters as (see
\cite{\ProcAF, Cor.~8.1; \ProcAI, Lemma~9.1 with $z=1$ and $b=0$})
$$\sp_{2n+1}(\rh;\x^{\pm1})=\underset \rh/\nu\text { a horizontal strip}
\to{\sum _{\nu\subseteq \rh} ^{}}\sp_{2n}(\nu;\x^{\pm1}).\tag4.7$$
Combining this with (4.6) yields (3.4) for $r\le n$,
as for any subpartition $\nu$ of the rectangle $(c^r)$ there is
exactly one way of adding a horizontal strip to $\nu$ such that a
partition $\rh$ with $\oddcols\!\big((c^r)/\rh\big)=0$ is obtained.

The case $r>n$ of (3.4) then follows from the $r\le n$ case by use of
(4.3) with $N=2n+1$.
\quad \quad \qed

\enddemo

\demo {Proof of {\rm(3.5)}}Here, in (3.2) we specialize $\x$ to
$(x_1,x_1^{-1},\dots,x_n,x_n^{-1},0,0,\dots)$ if $N=2n$, and
to $(x_1,x_1^{-1},\dots,x_n,x_n^{-1},1,0,0,\dots)$ if $N=2n+1$.
\quad \quad \qed
\enddemo

\subhead 5. Proof of Theorem~2\endsubhead
In this section we combine the tableaux description of symplectic
characters due to DeConcini and Procesi \cite{\DeCoAA, \DePrAA} (see
Section~A3 of the
Appendix) and tableaux descriptions of special orthogonal characters due to
Lakshmibai, Musili and Seshadri \cite{\LaMSAB, \LaShAC, \LaksAB}
(see Sections~A4 and A5 of the Appendix) 
with Robinson--Schensted--Knuth type
algorithms. We remark that while
there are restriction rules that would apply here (they involve ordinary
Littlewood--Richardson coefficients,
see \cite{\KingAB, (4.20)--(4.22); 
\KoTeAB, Theorem~2.1 with $k=n$, Theorem~A1}), these do
not appear to be very helpful for our purposes. This is because they
involve modification rules for characters, these cause alternating signs,
and these in turn cause a lot of cancellations, and all
this is simply not tractable for the applications that we have in mind.

\demo {Proof of {\rm (3.6)}}
By (A.4) we know that the left-hand side of (3.6) equals
$$\underset \text {of shape $(c^{n-p},(c-1)^p)$}\to
{\sum _{S\text { a $(2n)$-symplectic tableau}} ^{}}(\x^{\pm1})^S.\tag5.1$$
On the other hand, by (A.2) we have that
the right-hand side of (3.6) equals
$$(x_1x_2\cdots x_n)^{-c}\cdot\underset \text {with $\oddrows(\nu)=p$}\to
{\underset \text {of shape $\nu\subseteq\((2c)^{n}\)$}\to
{\sum _{T\text { an  $n$-tableau}} ^{}}}\x^T.\tag5.2$$
Comparing (5.1) and (5.2), we see that
(3.6) will be proved if we can find a bijection, $\Ph$ say, between
$(2n)$-symplectic tableaux $S$ of shape $(c^{n-p},(c-1)^p)$ and
$n$-tableaux $T$ with at most $2c$ columns and exactly $p$ odd rows
such that
$$(\x^{\pm1})^S=(x_1x_2\cdots x_n)^{-c}\cdot\x^T,\quad \text {if
}T=\Ph(S).\tag5.3$$

The bijection that we are going to construct proceeds in two steps.
In the first step we map $(2n)$-symplectic tableaux of shape
$(c^{n-p},(c-1)^p)$ to certain pairs (see the paragraph including
(5.4)) by analyzing how these
symplectic tableaux look like. In the second step we map these pairs
to the above described $n$-tableaux by a Robinson--Schensted--Knuth
type correspondence.

\medskip
{\it First step}. Let $S$ be a $(2n)$-symplectic tableaux of shape
$(c^{n-p},(c-1)^p)$. By the definition of $(2n)$-symplectic tableaux
in Section~A3 of the Appendix, 
$S$ is a $(2n)$-tableau of
shape $\((2c)^{n-p},(2c-2)^p\)$ such that columns $2c-1$,~$2c$, columns
$2c-3$,~$2c-2$,
etc., form $(2n)$-symplectic admissible pairs. An example
with $n=6$, $c=4$, $p=2$ is displayed in Figure~2.
\vskip10pt
\vbox{
$$
\Einheit0.5cm
\Pfad(0,0),11111122112222\endPfad
\Pfad(0,0),22222211111111\endPfad
\PfadDicke{1.5pt}
\Pfad(0,1),1122112111122\endPfad
\Pfad(0,1),2222211111111\endPfad
\PfadDicke{.5pt}
\Pfad(2,1),1111\endPfad
\Pfad(0,2),111111\endPfad
\Pfad(0,3),11111111\endPfad
\Pfad(0,4),11111111\endPfad
\Pfad(0,5),11111111\endPfad
\Pfad(1,0),222222\endPfad
\Pfad(2,0),222222\endPfad
\Pfad(3,0),222222\endPfad
\Pfad(4,0),222222\endPfad
\Pfad(5,0),222222\endPfad
\Pfad(6,0),222222\endPfad
\Pfad(7,2),2222\endPfad
\Label\ro{1}(0,5)
\Label\ro{1}(1,5)
\Label\ro{1}(2,5)
\Label\ro{1}(3,5)
\Label\ro{3}(4,5)
\Label\ro{3}(5,5)
\Label\ro{3}(6,5)
\Label\ro{4}(7,5)
\Label\ro{2}(0,4)
\Label\ro{2}(1,4)
\Label\ro{2}(2,4)
\Label\ro{3}(3,4)
\Label\ro{4}(4,4)
\Label\ro{4}(5,4)
\Label\ro{5}(6,4)
\Label\ro{5}(7,4)
\Label\ro{3}(0,3)
\Label\ro{4}(1,3)
\Label\ro{4}(2,3)
\Label\ro{5}(3,3)
\Label\ro{7}(4,3)
\Label\ro{7}(5,3)
\Label\ro{9}(6,3)
\Label\ro{10}(7,3)
\Label\ro{4}(0,2)
\Label\ro{5}(1,2)
\Label\ro{7}(2,2)
\Label\ro{7}(3,2)
\Label\ro{8}(4,2)
\Label\ro{8}(5,2)
\Label\ro{12}(6,2)
\Label\ro{12}(7,2)
\Label\ro{5}(0,1)
\Label\ro{6}(1,1)
\Label\ro{8}(2,1)
\Label\ro{9}(3,1)
\Label\ro{11}(4,1)
\Label\ro{11}(5,1)
\Label\ro{7}(0,0)
\Label\ro{10}(1,0)
\Label\ro{10}(2,0)
\Label\ro{11}(3,0)
\Label\ro{12}(4,0)
\Label\ro{12}(5,0)
\Label\o{\hskip0cm\oversetbrace \tsize 2c\to{\hbox{\hskip4cm}}}(4,6)
\Label\l{\hskip0cm\raise1.4cm\hbox{$\sideset n\and \to
     {\left\{\vbox{\vskip1.6cm}\right.}$}}(0,3)
\Label\r{\hskip0cm\raise.2cm\hbox{$\sideset \and p\to
     {\left.\vbox{\vskip.6cm}\right\}}$}}(8,1)
\hskip4cm
$$
\centerline{\eightpoint Figure 2}
}
\vskip10pt

We claim that $(2n)$-symplectic tableaux $S$ of shape
$(c^{n-p},(c-1)^p)$ are in bijection with pairs
$(\bT,\{e_1,e_2,\dots,e_p\})$, by a bijection $\Ph_1$ say,
where $\bT$ is an $n$-tableau whose
shape has only even rows and is contained in
$\((2c)^{n-p},(2c-2)^p\)$,
and where $\{e_1,e_2,\dots,e_p\}$ is a set of numbers satisfying
$$\spreadmatrixlines{2\jot}\matrix 1\le e_1<e_2<\dots<e_p\le n,\\
\text {and } e_l\notin [i(m),j(m)]\text { for }1\le l\le p,\ 1\le m\le s,
\endmatrix\tag5.4$$
given that 
$$\matrix i(1)&j(1)\\\vdots&\vdots \\ i(s)&j(s)\endmatrix$$
are the $(2c-1)$-st and $(2c)$-th column of $\bT$, such that
$$(\x^{\pm1})^S=(x_1x_2\cdots x_n)^{-c}\cdot x_{e_1}\cdots x_{e_p}\cdot
\x^{\bT},\quad \text {if }(\bT,\{e_1,e_2,\dots,e_p\})=\Ph_1(S).
\tag5.5$$

The construction of the bijection $\Ph_1$ is based on
an analysis of the symplectic tableaux under
consideration. By Observation~2 in Section~A3 
of the Appendix, both columns in a
$(2n)$-symplectic admissible pair have the same number of entries
$\le n$. Hence, the entries $\le n$ in $S$ form an $n$-tableau, $\bT$
say, with only even rows. Of course, the shape of $\bT$ is contained
in $\((2c)^{n-p},(2c-2)^p\)$. In Figure~2 we have marked the area that
is covered by entries $\le n$ by a bold line. 
The resulting tableau is displayed in the left half of
Figure~3.

The next observation is that, given $\bT$, we can recover $S$ almost
completely, only the $(2c-1)$-st and the $(2c)$-th column cannot be
necessarily recovered. For, all the columns of $S$ except
columns $2c-1$ and $2c$ have length $n$. Thus, by Observation~3
in Section~A3
of the Appendix, if the entries $\le n$ in a column of length $n$ are
$\{i(1),i(2),\dots,i(s)\}$ then the entries $>n$ are
$\{n+1,n+2,\dots,2n\}\backslash
\{2n+1-i(1),2n+1-i(2),\dots,2n+1-i(s)\}$. Only for recovering columns
$2c-1$ and $2c$ of $S$ from $\bT$ we need more information than just
the $(2c-1)$-st and $(2c)$-th column of $\bT$.

Let the $(2c-1)$-st and $(2c)$-th column of $S$ be
$$\matrix  
\left.\matrix i(1)&j(1)\\
\hbox to34.03pt{\hss$\vdots$\hss}&
\hbox to35.28pt{\hss$\vdots$\hss}\\ 
i(s)&j(s)\endmatrix\right\}
\text {entries }\le n\\
\left.\matrix i(s+1)&j(s+1)\\\vdots&\vdots \\ i(n-p)&j(n-p)\endmatrix\right\}
\text {entries }> n
\endmatrix\tag5.6$$
Observe that by definition of a $(2n)$-symplectic admissible pair
(see Definition~1 in the Appendix) we have
$$\multline
\{i(1),\dots,i(s),2n+1-i(s+1),\dots,2n+1-i(n-p)\}
\\=
\{j(1),\dots,j(s),2n+1-j(s+1),\dots,2n+1-j(n-p)\}.
\endmultline\tag5.7$$
Let $\{e_1,e_2,\dots,e_p\}$ be the complement of this set in
$\{1,2,\dots,n\}$. With other words, $e_1,e_2,\dots,e_p$ are the
numbers $e$ between 1 and $n$ with the property that neither $e$ nor its
``conjugate" $2n+1-e$ occur in the $(2c-1)$-st or $(2c)$-th column of
$S$. Without loss of generality we may assume $e_1<e_2<\dots<e_p$.
In our running example (recall $p=2$) we have $\{e_1,e_2\}=\{2,6\}$.

Obviously, because of (5.7) and the definition of $\bT$ and
$\{e_1,e_2,\dots,e_p\}$, 
all the information about the $(2c-1)$-st and the $(2c)$-th
column (displayed in (5.6)) is contained in the $(2c-1)$-st and $(2c)$-th
column of $\bT$ (the entries $\le n$ in (5.6)) and $\{e_1,e_2,\dots,e_p\}$.
Moreover, because of the definition of a $(2n)$-symplectic admissible
pair, we have $e_l\notin[i(m),j(m)]$, for all $l$ and $m$. And conversely,
if we have
$$\matrix i(1)&j(1)\\\vdots&\vdots \\ i(s)&j(s)\endmatrix,\quad
\{e_1,e_2,\dots,e_p\}$$
such that (5.4) is satisfied, then (5.6) with
$$\align
\{i(s+1),\dots,i(n-p)\}&:=\{n+1\dots,2n\}\backslash
\{2n+1-i(1),\dots,2n+1-i(s),\\
&\hskip4cm 2n+1-e_1,\dots,2n+1-e_p\}\\
\{j(s+1),\dots,j(n-p)\}&:=\{n+1\dots,2n\}\backslash
\{2n+1-j(1),\dots,2n+1-j(s),\\
&\hskip4cm 2n+1-e_1,\dots,2n+1-e_p\}
\endalign$$
will be a $(2n)$-symplectic admissible pair. 

Hence, we have defined the desired bijection $\Ph_1$. It is easy to
check that the weight property (5.5) holds under this correspondence.

Our running example in Figure~2 is mapped under $\Ph_1$ to the pair
in Figure~3.
\vskip10pt
\vbox{
$$
\(\hbox{\quad }
\Einheit0.5cm
\Pfad(0,-2),2222211111111\endPfad
\Pfad(0,-2),1122112111122\endPfad
\PfadDicke{.5pt}
\Pfad(0,-1),11\endPfad
\Pfad(0,0),11\endPfad
\Pfad(0,1),1111\endPfad
\Pfad(0,2),11111111\endPfad
\Pfad(1,-2),22222\endPfad
\Pfad(2,0),222\endPfad
\Pfad(3,0),222\endPfad
\Pfad(4,1),22\endPfad
\Pfad(5,1),22\endPfad
\Pfad(6,1),22\endPfad
\Pfad(7,1),22\endPfad
\Label\ro{1}(0,2)
\Label\ro{1}(1,2)
\Label\ro{1}(2,2)
\Label\ro{1}(3,2)
\Label\ro{3}(4,2)
\Label\ro{3}(5,2)
\Label\ro{3}(6,2)
\Label\ro{4}(7,2)
\Label\ro{2}(0,1)
\Label\ro{2}(1,1)
\Label\ro{2}(2,1)
\Label\ro{3}(3,1)
\Label\ro{4}(4,1)
\Label\ro{4}(5,1)
\Label\ro{5}(6,1)
\Label\ro{5}(7,1)
\Label\ro{3}(0,0)
\Label\ro{4}(1,0)
\Label\ro{4}(2,0)
\Label\ro{5}(3,0)
\Label\ro{4}(0,-1)
\Label\ro{5}(1,-1)
\Label\ro{5}(0,-2)
\Label\ro{6}(1,-2)
\hskip4cm
\quad ,\quad \{2,6\}
\)
$$
\centerline{\eightpoint Figure 3}
}
\vskip10pt

\medskip
{\it Second step}. In the second step we construct a bijection
$\Ph_2$ between pairs $(\bT,\{e_1,\mathbreak
e_2,\dots,e_p\})$ satisfying the above conditions including (5.4), and
$n$-tableaux $T$ with at most $2c$ columns and exactly $p$ odd rows,
such that
$$\x^T=x_{e_1}\cdots x_{e_p}\cdot \x^{\bT}.\tag5.8$$

Let $(\bT,\{e_1,e_2,\dots,e_p\})$ be such a pair. 
We insert $e_p,e_{p-1},\dots,e_1$, in this order, into $\bT$,
according to the following procedure. Let $\bT_0:=\bT$. Suppose that,
by inserting $e_p,e_{p-1},\dots,e_{p-l+1}$ we already formed
$\bT_l$. Next we insert $e_{p-l}$ into $\bT_l$ in the following way.
Choose the first row (from top to bottom) of $\bT_l$ such that
$e_{p-l}$ is less than the entry in
the $(2c-1)$-st column in that row of $\bT_l$. If there is no
such row choose the first row that does not have an entry in the
$(2c-1)$-st column of that row. Then, starting with that row of
$\bT_l$, ROW-INSERT $e_{p-l}$ into $\bT_l$, i.e., (cf\. \cite{\KnutAA,
p.~712; \KratAP, pp.~87--88}) find the leftmost entry in that row
that is larger than $e_{p-l}$, bump it and replace it by $e_{p-l}$, 
if there is none then place $e_{p-l}$ at the end of that row. If an
entry was bumped then repeat this same procedure with the bumped
entry and the next row, etc. Thus one obtains $\bT_{l+1}$. Finally, set
$T=\Ph_2\big((\bT,\{e_1,e_2,\dots,e_p\})\big):=\bT_p$. 
Our running example from Figure~3 is mapped under $\Ph_2$ to the
tableau in Figure~4.
\vskip10pt
\vbox{
$$
\Einheit0.5cm
\Pfad(0,-2),2222211111111\endPfad
\Pfad(0,-2),1121211211122\endPfad
\PfadDicke{.5pt}
\Pfad(0,-1),11\endPfad
\Pfad(0,0),111\endPfad
\Pfad(0,1),11111\endPfad
\Pfad(0,2),11111111\endPfad
\Pfad(1,-2),22222\endPfad
\Pfad(2,-1),2222\endPfad
\Pfad(3,0),222\endPfad
\Pfad(4,0),222\endPfad
\Pfad(5,1),22\endPfad
\Pfad(6,1),22\endPfad
\Pfad(7,1),22\endPfad
\Label\ro{1}(0,2)
\Label\ro{1}(1,2)
\Label\ro{1}(2,2)
\Label\ro{1}(3,2)
\Label\ro{2}(4,2)
\Label\ro{3}(5,2)
\Label\ro{3}(6,2)
\Label\ro{4}(7,2)
\Label\ro{2}(0,1)
\Label\ro{2}(1,1)
\Label\ro{2}(2,1)
\Label\ro{3}(3,1)
\Label\ro{3}(4,1)
\Label\ro{4}(5,1)
\Label\ro{5}(6,1)
\Label\ro{5}(7,1)
\Label\ro{3}(0,0)
\Label\ro{4}(1,0)
\Label\ro{4}(2,0)
\Label\ro{4}(3,0)
\Label\ro{6}(4,0)
\Label\ro{4}(0,-1)
\Label\ro{5}(1,-1)
\Label\ro{5}(2,-1)
\Label\ro{5}(0,-2)
\Label\ro{6}(1,-2)
\hskip4cm
$$
\centerline{\eightpoint Figure 4}
}
\vskip10pt

Since $\bT$
was an $n$-tableau with only even rows, and since later ``insertion
paths" 
are (weakly) to the left
of previous ones, $T$ is an $n$-tableau with exactly $p$ odd rows.
Trivially, (5.8) is satisfied.

To show that $\Ph_2$ is a bijection, we have to construct the inverse
mapping. Take an $n$-tableau $T$ with at most $2c$ columns and
exactly $p$ odd rows. Choose the last row (from top to bottom) of $T$
that has odd length. Now, starting with that row, perform a, slightly
modified, ROW-DELETE
(cf\. \cite{\KnutAA, p.~713; \KratAP, pp.~88}). Namely, remove the last
entry, $x$ say, from that row and find the rightmost entry, $y$ say,
in the previous row that is less than $x$, replace $y$ by $x$
and repeat this same procedure with $y$ and the row before the row
that contained $y$, etc., until no row is left to be
considered or until an entry in the $(2c)$-th column would be
replaced. In the latter case (this is the modification), 
do not replace the entry in the
$(2c)$-th column, but stop the procedure. Thus we obtain an
$n$-tableau $T_1$ with $p-1$ odd rows and an entry, $e_1$ say, that
was replaced in the last step of the procedure. This procedure is
repeated with $T_1$, thus obtaining $T_2$ and $e_2$, etc. In the end
we obtain $T_p$, which is an $n$-tableau with only even rows, and in
the course of our algorithm we obtained the elements
$e_1,e_2,\dots,e_p$. By standard properties of ROW-INSERT and
ROW-DELETE (cf\. \cite{\SagaAL}), it is not difficult to see 
that this algorithm exactly reverses $\Ph_2$, step by step. In
particular, the start of ROW-INSERT in a row that is possibly
different from the first and the modified ending of ROW-DELETE
complement each other exactly.

\medskip
The composition $\Ph_2\circ\Ph_1$ is by definition the desired
bijection between $(2n)$-sym\-plect\-ic tableaux $S$ of shape $(c^{n-p},(c-1)^p)$ and
$n$-tableaux $T$ with at most $2c$ columns and exactly $p$ odd rows.
From (5.5) and (5.8) the weight property (5.3) follows immediately.
This completes the proof of (3.6).\quad \quad \qed
\enddemo

\demo {Proof of {\rm (3.7)}}
Using (A.6) and (A.2) in (3.7), we see that
(3.7) will be proved if we can find a bijection, $\Ps$ say, between
$(2n)$-orthogonal tableaux $S$ of shape $(c^{n-1},c-p)$ and
$n$-tableaux $T$ with at most $2c$ columns and exactly $p$ columns
with parity different from $n$ (by which we mean that the {\it
lengths\/} of the columns have parity different from $n$)
such that
$$(\x^{\pm1})^S=(x_1x_2\cdots x_n)^{-c}\cdot\x^T,\quad \text {if
}T=\Ps(S).\tag5.9$$

Again, the bijection that we are going to construct proceeds in two steps.
In the first step we map $(2n)$-orthogonal tableaux of shape
$(c^{n-1},c-p)$ to certain pairs (see the paragraph including
(5.10)) by analyzing how these
orthogonal tableaux look like. In the second step we map these pairs
to the above described $n$-tableaux by Robinson--Schensted--Knuth
insertion.

\medskip
{\it First step}. Let $S$ be a $(2n)$-orthogonal tableaux of shape
$(c^{n-1},c-p)$. By Observation~1 in
Section~A5 of the Appendix, $S$ consists of a pair $(S_3,S_2)$
of $(2n)$-tableaux, $S_3$ being of shape $((2c-p)^n)$ and each column
of which containing an even number of entries $>n$,  
$S_2$ being of shape $(p^n)$, each column
of which containing an odd number of entries $>n$, and all the
entries in the first row of $S_2$ being at most $n$, such that the
concatenation $S_3\cup S_2'$ is a $(2n)$-tableau, 
where $S_2'$ is the tableau arising
from $S_2$ by replacing the topmost element, $e_i$ say, in column $i$
of $S_2$ by its ``conjugate" $2n+1-e_i$, for all $i=1,2,\dots,p$, and
by rearranging the columns in increasing order.
An example with $n=6$, $c=7/2$, $p=2$ is displayed in the left half
of Figure~5. The right half shows the concatenation $S_3\cup S_2'$
(note that $e_1=1$, $e_2=2$).
\vskip10pt
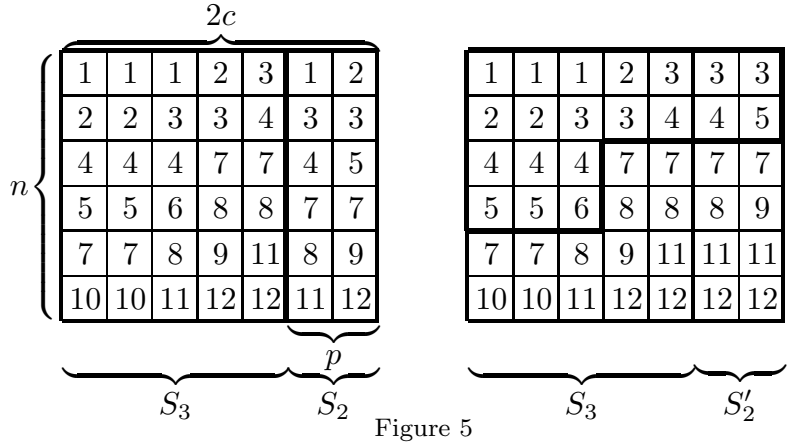
\vbox{
$$
\Einheit0.5cm
\Pfad(0,0),1111111222222\endPfad
\Pfad(0,0),2222221111111\endPfad
\PfadDicke{1.5pt}
\Pfad(5,0),222222\endPfad
\PfadDicke{.5pt}
\Pfad(0,1),1111111\endPfad
\Pfad(0,2),1111111\endPfad
\Pfad(0,3),1111111\endPfad
\Pfad(0,4),1111111\endPfad
\Pfad(0,5),1111111\endPfad
\Pfad(1,0),222222\endPfad
\Pfad(2,0),222222\endPfad
\Pfad(3,0),222222\endPfad
\Pfad(4,0),222222\endPfad
\Pfad(5,0),222222\endPfad
\Pfad(6,0),222222\endPfad
\Label\ro{1}(0,5)
\Label\ro{1}(1,5)
\Label\ro{1}(2,5)
\Label\ro{2}(3,5)
\Label\ro{3}(4,5)
\Label\ro{1}(5,5)
\Label\ro{2}(6,5)
\Label\ro{2}(0,4)
\Label\ro{2}(1,4)
\Label\ro{3}(2,4)
\Label\ro{3}(3,4)
\Label\ro{4}(4,4)
\Label\ro{3}(5,4)
\Label\ro{3}(6,4)
\Label\ro{4}(0,3)
\Label\ro{4}(1,3)
\Label\ro{4}(2,3)
\Label\ro{7}(3,3)
\Label\ro{7}(4,3)
\Label\ro{4}(5,3)
\Label\ro{5}(6,3)
\Label\ro{5}(0,2)
\Label\ro{5}(1,2)
\Label\ro{6}(2,2)
\Label\ro{8}(3,2)
\Label\ro{8}(4,2)
\Label\ro{7}(5,2)
\Label\ro{7}(6,2)
\Label\ro{7}(0,1)
\Label\ro{7}(1,1)
\Label\ro{8}(2,1)
\Label\ro{9}(3,1)
\Label\ro{11}(4,1)
\Label\ro{8}(5,1)
\Label\ro{9}(6,1)
\Label\ro{10}(0,0)
\Label\ro{10}(1,0)
\Label\ro{11}(2,0)
\Label\ro{12}(3,0)
\Label\ro{12}(4,0)
\Label\ro{11}(5,0)
\Label\ro{12}(6,0)
\Label\o{\hskip-0.5cm\oversetbrace \tsize 2c\to{\hbox{\hskip3.5cm}}}(4,6)
\Label\l{\hskip0cm\raise1.4cm\hbox{$\sideset n\and \to
     {\left\{\vbox{\vskip1.6cm}\right.}$}}(0,3)
\Label\u{\hskip.5cm\undersetbrace \tsize S_3\to{\hbox{\hskip2.5cm}}}(2,-1)
\Label\u{\hskip0cm\undersetbrace \tsize p\to{\hbox{\hskip1cm}}}(6,0)
\Label\u{\hskip0cm\undersetbrace \tsize S_2\to{\hbox{\hskip1cm}}}(6,-1)
\hskip3.5cm
\hbox{\hskip1cm}
\PfadDicke{1pt}
\Pfad(0,0),1111111222222\endPfad
\Pfad(0,0),2222221111111\endPfad
\Pfad(5,0),222222\endPfad
\PfadDicke{1.5pt}
\Pfad(0,2),11122111122\endPfad
\Pfad(0,2),22221111111\endPfad
\PfadDicke{.5pt}
\Pfad(0,1),1111111\endPfad
\Pfad(0,2),1111111\endPfad
\Pfad(0,3),1111111\endPfad
\Pfad(0,4),1111111\endPfad
\Pfad(0,5),1111111\endPfad
\Pfad(1,0),222222\endPfad
\Pfad(2,0),222222\endPfad
\Pfad(3,0),222222\endPfad
\Pfad(4,0),222222\endPfad
\Pfad(5,0),222222\endPfad
\Pfad(6,0),222222\endPfad
\Label\ro{1}(0,5)
\Label\ro{1}(1,5)
\Label\ro{1}(2,5)
\Label\ro{2}(3,5)
\Label\ro{3}(4,5)
\Label\ro{3}(5,5)
\Label\ro{3}(6,5)
\Label\ro{2}(0,4)
\Label\ro{2}(1,4)
\Label\ro{3}(2,4)
\Label\ro{3}(3,4)
\Label\ro{4}(4,4)
\Label\ro{4}(5,4)
\Label\ro{5}(6,4)
\Label\ro{4}(0,3)
\Label\ro{4}(1,3)
\Label\ro{4}(2,3)
\Label\ro{7}(3,3)
\Label\ro{7}(4,3)
\Label\ro{7}(5,3)
\Label\ro{7}(6,3)
\Label\ro{5}(0,2)
\Label\ro{5}(1,2)
\Label\ro{6}(2,2)
\Label\ro{8}(3,2)
\Label\ro{8}(4,2)
\Label\ro{8}(5,2)
\Label\ro{9}(6,2)
\Label\ro{7}(0,1)
\Label\ro{7}(1,1)
\Label\ro{8}(2,1)
\Label\ro{9}(3,1)
\Label\ro{11}(4,1)
\Label\ro{11}(5,1)
\Label\ro{11}(6,1)
\Label\ro{10}(0,0)
\Label\ro{10}(1,0)
\Label\ro{11}(2,0)
\Label\ro{12}(3,0)
\Label\ro{12}(4,0)
\Label\ro{12}(5,0)
\Label\ro{12}(6,0)
\Label\u{\hskip.5cm\undersetbrace \tsize S_3\to{\hbox{\hskip2.5cm}}}(2,-1)
\Label\u{\hskip0cm\undersetbrace \tsize S_2'\to{\hbox{\hskip1cm}}}(6,-1)
\hskip3.5cm
$$
\centerline{\eightpoint Figure 5}
}
\vskip10pt

We claim that $(2n)$-orthogonal tableaux $S$ of shape
$(c^{n-1},c-p)$ are in bijection with pairs
$(\bT,\{e_1,e_2,\dots,e_p\})$, by a bijection $\Ps_1$ say,
where $\bT$ is an $n$-tableau whose
shape has only columns with the same parity as $n$ and is contained in
$\((2c)^{n-1},2c-p\)$,
and where $\{e_1,e_2,\dots,e_p\}$ is a set of numbers satisfying
$$\spreadmatrixlines{2\jot}\matrix 1\le e_1\le e_2\le \dots\le e_p\le n,\\
e_i\text { is less than the topmost element of the $(2c-p+i)$-th
column of $\bT$,}
\endmatrix\tag5.10$$
such that
$$(\x^{\pm1})^S=(x_1x_2\cdots x_n)^{-c}\cdot x_{e_1}\cdots x_{e_p}\cdot
\x^{\bT},\quad \text {if }(\bT,\{e_1,e_2,\dots,e_p\})=\Ps_1(S).
\tag5.11$$

The construction of the bijection $\Ps_1$ is based on
an analysis of the orthogonal tableaux under
consideration. By definition, all the columns in $S_3\cup S_2'$
contain an even number of entries $>n$.
Hence, the entries $\le n$ in $S_3\cup S_2'$ form an $n$-tableau, $\bT$
say, with all columns having the same parity as $n$. 
Of course, the shape of $\bT$ is contained
in $\((2c)^{n-1},2c-p\)$. In the right half of
Figure~5 we have marked the area that
is covered by entries $\le n$ by a bold line. 
The resulting tableau is displayed in the left half of
Figure~6. 

Now, let $\Ps_1(S)$ be defined by
$(\bT,\{e_1,e_2,\dots,e_p\})$, where, as before, 
$e_i$ is the topmost element of
the $i$-th column of $S_2$. Our running example in Figure~5 is mapped under
$\Ps_1$ to the pair in Figure~6. It is obvious that 
(5.10) and (5.11) hold under this mapping. 
Besides, it is trivial to recover $S$ from
$(\bT,\{e_1,e_2,\dots,e_p\})$. Hence, $\Ps_1$ is a bijection, as
desired.
\vskip10pt
\vbox{
$$
\(\hbox{\quad }
\Einheit0.5cm
\Pfad(0,-2),22221111111\endPfad
\Pfad(0,-2),11122111122\endPfad
\PfadDicke{.5pt}
\Pfad(0,-1),111\endPfad
\Pfad(0,0),111\endPfad
\Pfad(0,1),1111111\endPfad
\Pfad(1,-2),2222\endPfad
\Pfad(2,-2),2222\endPfad
\Pfad(3,0),22\endPfad
\Pfad(4,0),22\endPfad
\Pfad(5,0),22\endPfad
\Pfad(6,0),22\endPfad
\Label\ro{1}(0,1)
\Label\ro{1}(1,1)
\Label\ro{1}(2,1)
\Label\ro{2}(3,1)
\Label\ro{3}(4,1)
\Label\ro{3}(5,1)
\Label\ro{3}(6,1)
\Label\ro{2}(0,0)
\Label\ro{2}(1,0)
\Label\ro{3}(2,0)
\Label\ro{3}(3,0)
\Label\ro{4}(4,0)
\Label\ro{4}(5,0)
\Label\ro{5}(6,0)
\Label\ro{4}(0,-1)
\Label\ro{4}(1,-1)
\Label\ro{4}(2,-1)
\Label\ro{5}(0,-2)
\Label\ro{5}(1,-2)
\Label\ro{6}(2,-2)
\hskip3.5cm
\quad ,\quad \{1,2\}
\)
$$
\centerline{\eightpoint Figure 6}
}
\vskip10pt

\medskip
{\it Second step}. In the second step we construct a bijection
$\Ps_2$ between pairs $(\bT,\{e_1,\mathbreak
e_2,\dots,e_p\})$ satisfying the above conditions including (5.10), and
$n$-tableaux $T$ with at most $2c$ columns and exactly $p$ columns of
parity different from $n$, such that
$$\x^T=x_{e_1}\cdots x_{e_p}\cdot \x^{\bT}.\tag5.12$$

Let $(\bT,\{e_1,e_2,\dots,e_p\})$ be such a pair. 
We ROW-INSERT $e_1,e_{2},\dots,e_p$, in this order, into $\bT$. (The
reader should observe the differences to the algorithm in the second
step of the proof of (3.6): Here we use {\it genuine} ROW-INSERT,
without any modification, as opposed to the proof of
(3.6). Furthermore, the order in which the elements
$e_1,e_2,\dots,e_p$ are inserted is exactly the reversed order in the
proof of (3.6).) I.e.,
(cf\. \cite{\KnutAA,
p.~712; \KratAP, pp.~87--88})
let $\bT_0:=\bT$. Suppose that,
by inserting $e_1,e_{2},\dots,e_{l}$ we already formed
$\bT_l$. Next we insert $e_{l+1}$ into $\bT_l$. Namely,
we find the leftmost entry in the first row of $\bT_l$
that is larger than $e_{l+1}$, bump it and replace it by $e_{l+1}$, 
if there is none we place $e_{l+1}$ at the end of that row. If an
entry was bumped then we repeat this same procedure with the bumped
entry and the next row, etc. Thus one obtains the
tableau $\bT_{l+1}$. Finally, we set
$T=\Ps_2\big((\bT,\{e_1,e_2,\dots,e_p\})\big):=\bT_p$. 
Our running example from Figure~6 is mapped under $\Ps_2$ to the
tableau in Figure~7.
\vskip10pt
\vbox{
$$
\Einheit0.5cm
\Pfad(0,-2),222221111111\endPfad
\Pfad(0,-2),121121211122\endPfad
\PfadDicke{.5pt}
\Pfad(0,-1),1\endPfad
\Pfad(0,0),111\endPfad
\Pfad(0,1),1111\endPfad
\Pfad(0,2),1111111\endPfad
\Pfad(1,-2),22222\endPfad
\Pfad(2,-1),2222\endPfad
\Pfad(3,0),222\endPfad
\Pfad(4,1),22\endPfad
\Pfad(5,1),22\endPfad
\Pfad(6,1),22\endPfad
\Label\ro{1}(0,2)
\Label\ro{1}(1,2)
\Label\ro{1}(2,2)
\Label\ro{1}(3,2)
\Label\ro{2}(4,2)
\Label\ro{3}(5,2)
\Label\ro{3}(6,2)
\Label\ro{2}(0,1)
\Label\ro{2}(1,1)
\Label\ro{2}(2,1)
\Label\ro{3}(3,1)
\Label\ro{3}(4,1)
\Label\ro{4}(5,1)
\Label\ro{5}(6,1)
\Label\ro{3}(0,0)
\Label\ro{4}(1,0)
\Label\ro{4}(2,0)
\Label\ro{4}(3,0)
\Label\ro{4}(0,-1)
\Label\ro{5}(1,-1)
\Label\ro{6}(2,-1)
\Label\ro{5}(0,-2)
\hskip3.5cm
$$
\centerline{\eightpoint Figure 7}
}
\vskip10pt

Since $\bT$
was an $n$-tableau with all columns having the same parity as $n$, 
and since later ``insertion
paths" are strictly to the right
of previous ones, $T$ is an $n$-tableau with exactly $p$ columns of
parity different from $n$. Moreover, since $e_i$ is 
less than the topmost element of the $(2c-p+i)$-th
column of $\bT$, the insertion process will at no stage produce more
than $2c$ columns. Hence, $T$ has at most $2c$ columns.
Trivially, (5.12) is satisfied.

To show that $\Ps_2$ is a bijection, we have to construct the inverse
mapping. Take an $n$-tableau $T$ with at most $2c$ columns and
exactly $p$ columns of parity different from $n$. 
Choose the rightmost column of $T$
that has parity different from $n$. Suppose that the bottommost
entry in this column is in row $I$.
Now, starting with that row, perform ROW-DELETE,
i.e., (cf\. \cite{\KnutAA, p.~713; \KratAP, pp.~88}) remove the last
entry, $x_I$ say, from the $I$-th row and find the rightmost entry,
$x_{I-1}$ say,
in the $(I-1)$-st row that is less than $x_I$, replace $x_{I-1}$ by
$x_I$
and repeat this same procedure with $x_{I-1}$ and the $(I-2)$-nd row, 
etc., until no row is left to be
considered. Thus we obtain an
$n$-tableau $T_1$ with $p-1$ columns of parity different from $n$ 
and an entry, $e_p$ say, that
was replaced in the last step of the procedure. This procedure is
repeated with $T_1$, thus obtaining $T_2$ and $e_{p-1}$, etc. In the end
we obtain $T_p$, which is an $n$-tableau with all columns having the
same parity as $n$, and in
the course of our algorithm we obtained the elements
$e_p,e_{p-1},\dots,e_1$. Again, by standard properties of ROW-INSERT and
ROW-DELETE (cf\. \cite{\SagaAL}), it is not difficult to see 
that this algorithm exactly reverses $\Ps_2$, step by step.

\medskip
The composition $\Ps_2\circ\Ps_1$ is by definition the desired
bijection between $(2n)$-orth\-og\-on\-al tableaux $S$ of shape $(c^{n-1},c-p)$ and
$n$-tableaux $T$ with at most $2c$ columns and exactly $p$ columns of
parity different from $n$.
From (5.11) and (5.12) the weight property (5.9) follows immediately.
This completes the proof of (3.7).\quad \quad \qed
\enddemo

\demo {Proof of {\rm (3.8)}}
Here we have to deal with $(2n+1)$-orthogonal tableaux. Since these
are not too far from $(2n)$-symplectic tableaux, the arguments here
are very similar to those in the proof of (3.6). In fact, the basic
steps are the same, only the details differ. So we shall be sometimes
sketchy and provide details only if necessary.

Using (A.5) and (A.2) in (3.8), we see that
(3.8) will be proved if we can find a bijection, $\Th$ say, between
$(2n+1)$-orthogonal tableaux $S$ of shape $(c^{n-p},(c-1)^p)$ and
pairs $(T,\si)$, where $T$ is an
$n$-tableau whose shape is contained in $((2c)^{n})$,
and where $\si$ is a vertical strip of length $p$
on the rim of $\nu$ that avoids the $(2c)$-th column,
such that
$$(\x^{\pm1})^S=(x_1x_2\cdots x_n)^{-c}\cdot\x^T,\quad \text {if
}(T,\si)=\Th(S),\text { for some }\si.\tag5.13$$

Again, we proceed in two steps.

\medskip
{\it First step}. Let $S$ be a $(2n+1)$-orthogonal tableaux of shape
$(c^{n-p},(c-1)^p)$. By the definition of $(2n+1)$-orthogonal tableaux
in Section~A4 of the Appendix, $S$ is a $(2n)$-tableau of
shape $\((2c)^{n-p},(2c-2)^p\)$ such that columns $2c-1$,~$2c$, columns
$2c-3$,~$2c-2$,
etc., form $(2n+1)$-orthogonal admissible pairs. An example
with $n=6$, $c=4$, $p=3$ is displayed in Figure~8.
\vskip10pt
\vbox{
$$
\Einheit0.5cm
\Pfad(0,0),11111122211222\endPfad
\Pfad(0,0),22222211111111\endPfad
\PfadDicke{1.5pt}
\Pfad(0,0),12112212111122\endPfad
\Pfad(0,0),22222211111111\endPfad
\PfadDicke{.5pt}
\Pfad(0,1),111111\endPfad
\Pfad(0,2),111111\endPfad
\Pfad(0,3),111111\endPfad
\Pfad(0,4),11111111\endPfad
\Pfad(0,5),11111111\endPfad
\Pfad(1,0),222222\endPfad
\Pfad(2,0),222222\endPfad
\Pfad(3,0),222222\endPfad
\Pfad(4,0),222222\endPfad
\Pfad(5,0),222222\endPfad
\Pfad(6,0),222222\endPfad
\Pfad(7,3),222\endPfad
\Label\ro{1}(0,5)
\Label\ro{1}(1,5)
\Label\ro{1}(2,5)
\Label\ro{1}(3,5)
\Label\ro{3}(4,5)
\Label\ro{3}(5,5)
\Label\ro{4}(6,5)
\Label\ro{4}(7,5)
\Label\ro{2}(0,4)
\Label\ro{2}(1,4)
\Label\ro{2}(2,4)
\Label\ro{3}(3,4)
\Label\ro{4}(4,4)
\Label\ro{5}(5,4)
\Label\ro{5}(6,4)
\Label\ro{5}(7,4)
\Label\ro{3}(0,3)
\Label\ro{3}(1,3)
\Label\ro{4}(2,3)
\Label\ro{6}(3,3)
\Label\ro{7}(4,3)
\Label\ro{7}(5,3)
\Label\ro{12}(6,3)
\Label\ro{12}(7,3)
\Label\ro{4}(0,2)
\Label\ro{5}(1,2)
\Label\ro{5}(2,2)
\Label\ro{8}(3,2)
\Label\ro{8}(4,2)
\Label\ro{9}(5,2)
\Label\ro{5}(0,1)
\Label\ro{6}(1,1)
\Label\ro{6}(2,1)
\Label\ro{9}(3,1)
\Label\ro{11}(4,1)
\Label\ro{11}(5,1)
\Label\ro{6}(0,0)
\Label\ro{9}(1,0)
\Label\ro{10}(2,0)
\Label\ro{11}(3,0)
\Label\ro{12}(4,0)
\Label\ro{12}(5,0)
\Label\o{\hskip0cm\oversetbrace \tsize 2c\to{\hbox{\hskip4cm}}}(4,6)
\Label\l{\hskip0cm\raise1.4cm\hbox{$\sideset n\and \to
     {\left\{\vbox{\vskip1.6cm}\right.}$}}(0,3)
\Label\r{\hskip0cm\raise1.2cm\hbox{$\sideset \and p\to
     {\left.\vbox{\vskip.9cm}\right\}}$}}(8,1)
\hskip4cm
$$
\centerline{\eightpoint Figure 8}
}
\vskip10pt

We claim that $(2n+1)$-orthogonal tableaux $S$ of shape
$(c^{n-p},(c-1)^p)$ are in bijection with pairs
$(\bT,\{e_1,e_2,\dots,e_p\})$, by a bijection $\Th_1$ say,
where $\bT$ is an $n$-tableau contained in
$\((2c)^{n-p},(2c-2)^p\)$,
and where $\{e_1,e_2,\dots,e_p\}$ is a set of numbers satisfying
$$\spreadmatrixlines{2\jot}\matrix 1\le e_1<e_2<\dots<e_p\le n,\\
 e_l\notin [i(m),j(m)]\text { for }1\le l\le p,\ 1\le m\le
s,\\
\text {and } e_l\notin [i(s+1),n]\text { for }1\le l\le p,
\endmatrix\tag5.14$$
given that 
$$\matrix i(1)&j(1)\\\vdots&\vdots \\ i(s)&j(s)\\
i(s+1)\\\vdots\\i(t)\endmatrix$$
are the $(2c-1)$-st and $(2c)$-th column of $\bT$, such that
$$(\x^{\pm1})^S=(x_1x_2\cdots x_n)^{-c}\cdot x_{e_1}\cdots x_{e_p}\cdot
\x^{\bT},\quad \text {if }(\bT,\{e_1,e_2,\dots,e_p\})=\Th_1(S).
\tag5.15$$

The construction of the bijection $\Th_1$ is based on
an analysis of the orthogonal tableaux under
consideration. 
Clearly, the entries $\le n$ in $S$ form an $n$-tableau, $\bT$
say, whose shape is contained
in $\((2c)^{n-p},(2c-2)^p\)$. Note that in difference
to the symplectic case, here $\bT$ is not necessarily a tableau with
only even rows. In Figure~8 we have marked the area that
is covered by entries $\le n$ by a bold line. 
The resulting tableau is displayed in the left half of
Figure~9.

Now, let the $(2c-1)$-st and $(2c)$-th column of $S$ be
$$\matrix  
n\ge\left\{
\matrix \vphantom{j}i(1)\\
\hbox to34.03pt{\hss$\vdots$\hss}
\\\vdots\\\vphantom{j}i(t) \endmatrix\right.
\\
n<\left\{\matrix \vphantom{j}i(t+1)\\
\hbox to34.03pt{\hss$\vdots$\hss}
\\\vphantom{j}i(n-p)\endmatrix\right.
\endmatrix
\matrix  
\left.
\matrix j(1)\\
\hbox to34.03pt{\hss$\vdots$\hss}
\\j(s) \endmatrix\right\}\le n
\\
\left.\matrix j(s+1)\\\vdots\\\vdots\\j(n-p)\endmatrix\right\}>n
\endmatrix
\tag5.16$$
As in the symplectic case, we define $\{e_1,e_2,\dots,e_p\}$ to be the
set of numbers $e$ between 1 and $n$ with the property that neither $e$ nor its
``conjugate" $2n+1-e$ occur in the $(2c-1)$-st or $(2c)$-th column of
$S$. Again, without loss of generality we may assume $e_1<e_2<\dots<e_p$.
In our running example (recall $p=3$) we have
$\{e_1,e_2,e_3\}=\{2,3,6\}$.
From the definition of a $(2n+1)$-orthogonal admissible pair (see
Definition~2 in the Appendix)
it follows that the numbers $e_1,e_2,\dots,e_p$ satisfy (5.14). Note
that (5.14) differs from the ``symplectic analogue" (5.4) by the
additional condition in the third line. 

We define $\Th_1(S)$ to be $(\bT,\{e_1,e_2,\dots,e_p\})$. 
Our running example in Figure~8 is mapped under $\Th_1$ to the pair
in Figure~9.
\vskip10pt
\vbox{
$$
\(\hbox{\quad }
\Einheit0.5cm
\Pfad(0,-3),22222211111111\endPfad
\Pfad(0,-3),12112212111122\endPfad
\PfadDicke{.5pt}
\Pfad(0,-2),1\endPfad
\Pfad(0,-1),111\endPfad
\Pfad(0,0),1111\endPfad
\Pfad(0,1),11111\endPfad
\Pfad(0,2),11111111\endPfad
\Pfad(1,-2),22222\endPfad
\Pfad(2,-2),22222\endPfad
\Pfad(3,-2),22222\endPfad
\Pfad(4,0),222\endPfad
\Pfad(5,1),22\endPfad
\Pfad(6,1),22\endPfad
\Pfad(7,1),22\endPfad
\Label\ro{1}(0,2)
\Label\ro{1}(1,2)
\Label\ro{1}(2,2)
\Label\ro{1}(3,2)
\Label\ro{3}(4,2)
\Label\ro{3}(5,2)
\Label\ro{4}(6,2)
\Label\ro{4}(7,2)
\Label\ro{2}(0,1)
\Label\ro{2}(1,1)
\Label\ro{2}(2,1)
\Label\ro{3}(3,1)
\Label\ro{4}(4,1)
\Label\ro{5}(5,1)
\Label\ro{5}(6,1)
\Label\ro{5}(7,1)
\Label\ro{3}(0,0)
\Label\ro{3}(1,0)
\Label\ro{4}(2,0)
\Label\ro{6}(3,0)
\Label\ro{4}(0,-1)
\Label\ro{5}(1,-1)
\Label\ro{5}(2,-1)
\Label\ro{5}(0,-2)
\Label\ro{6}(1,-2)
\Label\ro{6}(2,-2)
\Label\ro{6}(0,-3)
\hskip4cm
\quad ,\quad \{2,3,6\}
\)
$$
\centerline{\eightpoint Figure 9}
}
\vskip10pt
\noindent
As in the
symplectic case, it can be shown that $S$ can be uniquely recovered
from $(\bT,\{e_1,e_2,\dots,e_p\})$. Also, it is easy to
check that the weight property (5.15) holds under this correspondence.

\medskip
{\it Second step}. In the second step we construct a bijection
$\Th_2$ between pairs $(\bT,\{e_1,\mathbreak
e_2,\dots,e_p\})$ satisfying (5.14) 
as before, and pairs $(T,\si)$, where $T$ is an
$n$-tableau whose shape is contained in $((2c)^{n})$,
and where $\si$ is a vertical strip of length $p$
on the rim of $\nu$ that avoids the $(2c)$-th column,
such that
$$\x^T=x_{e_1}\cdots x_{e_p}\cdot \x^{\bT}.\tag5.17$$

To obtain $T$ from such a pair $(\bT,\{e_1,e_2,\dots,e_p\})$ we use
the mapping $\Ph_2$ from the second step of the proof of (3.6). On
the other hand, $\si$ is defined to be the vertical strip by which
the shapes of $T$ and $\bT$ differ. 
Thus, our running example from Figure~9 is mapped under $\Th_2$ to the
pair in Figure~10. There, the vertical strip $\si$ is visualized by bold
lines embedded in the shape of $T$.
\vskip10pt
\vbox{
$$
\(\hbox{\quad }
\Einheit0.5cm
\Pfad(0,-3),22222211111111\endPfad
\Pfad(0,-3),12111221211122\endPfad
\PfadDicke{.5pt}
\Pfad(0,-2),1\endPfad
\Pfad(0,-1),1111\endPfad
\Pfad(0,0),1111\endPfad
\Pfad(0,1),11111\endPfad
\Pfad(0,2),11111111\endPfad
\Pfad(1,-2),22222\endPfad
\Pfad(2,-2),22222\endPfad
\Pfad(3,-2),22222\endPfad
\Pfad(4,0),222\endPfad
\Pfad(5,1),22\endPfad
\Pfad(6,1),22\endPfad
\Pfad(7,1),22\endPfad
\Label\ro{1}(0,2)
\Label\ro{1}(1,2)
\Label\ro{1}(2,2)
\Label\ro{1}(3,2)
\Label\ro{2}(4,2)
\Label\ro{3}(5,2)
\Label\ro{3}(6,2)
\Label\ro{4}(7,2)
\Label\ro{2}(0,1)
\Label\ro{2}(1,1)
\Label\ro{2}(2,1)
\Label\ro{3}(3,1)
\Label\ro{3}(4,1)
\Label\ro{4}(5,1)
\Label\ro{5}(6,1)
\Label\ro{5}(7,1)
\Label\ro{3}(0,0)
\Label\ro{3}(1,0)
\Label\ro{4}(2,0)
\Label\ro{4}(3,0)
\Label\ro{6}(4,0)
\Label\ro{4}(0,-1)
\Label\ro{5}(1,-1)
\Label\ro{5}(2,-1)
\Label\ro{5}(3,-1)
\Label\ro{5}(0,-2)
\Label\ro{6}(1,-2)
\Label\ro{6}(2,-2)
\Label\ro{6}(3,-2)
\Label\ro{6}(0,-3)
\hskip4cm
\hbox{\quad ,\quad }
\PfadDicke{.3pt}
\Pfad(0,-3),22222211111111\endPfad
\Pfad(0,-3),12111221211122\endPfad
\PfadDicke{1.5pt}
\Pfad(3,-2),12212\endPfad
\Pfad(3,-2),22121\endPfad
\PfadDicke{.3pt}
\Pfad(0,-2),1\endPfad
\Pfad(0,-1),1111\endPfad
\Pfad(0,0),1111\endPfad
\Pfad(0,1),11111\endPfad
\Pfad(0,2),11111111\endPfad
\Pfad(1,-2),22222\endPfad
\Pfad(2,-2),22222\endPfad
\Pfad(3,-2),22222\endPfad
\Pfad(4,0),222\endPfad
\Pfad(5,1),22\endPfad
\Pfad(6,1),22\endPfad
\Pfad(7,1),22\endPfad
\hskip4cm
\hbox{\quad }
\)
$$
\centerline{\eightpoint Figure 10}
}
\vskip10pt

By definition of $\Th_2$ (recall
that an element is never inserted into the $(2c)$-th column of some
$\bT_i$), $\si$ is a vertical strip on the rim of the shape of $T$
that avoids the $(2c)$-th column. It is obvious that (5.17) holds.

\medskip
The composition $\Th_2\circ\Th_1$ is by definition the desired
bijection between $(2n+1)$-orth\-og\-on\-al tableaux $S$ of shape $(c^{n-p},(c-1)^p)$ and
pairs $(T,\si)$, where $T$ is an
$n$-tableau whose shape is contained in $((2c)^{n})$,
and where $\si$ is a vertical strip of length $p$
on the rim of $\nu$ that avoids the $(2c)$-th column.
From (5.15) and (5.17) the weight property (5.13) follows immediately.
This completes the proof of (3.8).\quad \quad \qed
\enddemo

\demo {Proof of {\rm (3.10)}}
Here we proceed very similarly to the preceding proof of (3.8).

Using (A.5) and (A.2) in (3.10),
we see that
(3.10) will be proved if we can find a bijection, $\Om$ say, between
$(2n+1)$-orthogonal tableaux $S$ of shape $(c^{n-1},c-p)$ and
pairs $(T,\si)$, where $T$ is an
$n$-tableau whose shape is contained in $((2c)^{n})$,
and where $\si$ is a horizontal strip of length $p$
on the rim of $\nu$ such that the $i$-th cell of the
strip comes before the $(2c-2p+2i)$-th column,
such that
$$(\x^{\pm1})^S=(x_1x_2\cdots x_n)^{-c}\cdot\x^T,\quad \text {if
}(T,\si)=\Om(S),\text { for some }\si.\tag5.18$$

Again, we proceed in two steps.

\medskip
{\it First step}. Let $S$ be a $(2n+1)$-orthogonal tableau of shape
$(c^{n-1},c-p)$. By the definition of $(2n+1)$-orthogonal tableaux in
Section~A4 of the Appendix, $S$ is a $(2n)$-tableau of
shape $\((2c)^{n-1},2c-2p\)$ such that columns $2c-1$,~$2c$, columns
$2c-3$,~$2c-2$,
etc., form $(2n+1)$-orthogonal admissible pairs. An example
with $n=6$, $c=7/2$, $p=3$ is displayed in Figure~11.
\vskip10pt
\vbox{
$$
\Einheit0.5cm
\Pfad(0,0),1211111122222\endPfad
\Pfad(0,0),2222221111111\endPfad
\PfadDicke{1.5pt}
\Pfad(0,1),122111211122\endPfad
\Pfad(0,1),222221111111\endPfad
\PfadDicke{.5pt}
\Pfad(0,1),1\endPfad
\Pfad(0,2),1111111\endPfad
\Pfad(0,3),1111111\endPfad
\Pfad(0,4),1111111\endPfad
\Pfad(0,5),1111111\endPfad
\Pfad(1,1),22222\endPfad
\Pfad(2,1),22222\endPfad
\Pfad(3,1),22222\endPfad
\Pfad(4,1),22222\endPfad
\Pfad(5,1),22222\endPfad
\Pfad(6,1),22222\endPfad
\Label\ro{1}(0,5)
\Label\ro{1}(1,5)
\Label\ro{1}(2,5)
\Label\ro{1}(3,5)
\Label\ro{2}(4,5)
\Label\ro{3}(5,5)
\Label\ro{4}(6,5)
\Label\ro{2}(0,4)
\Label\ro{2}(1,4)
\Label\ro{2}(2,4)
\Label\ro{5}(3,4)
\Label\ro{5}(4,4)
\Label\ro{5}(5,4)
\Label\ro{5}(6,4)
\Label\ro{3}(0,3)
\Label\ro{5}(1,3)
\Label\ro{5}(2,3)
\Label\ro{6}(3,3)
\Label\ro{7}(4,3)
\Label\ro{9}(5,3)
\Label\ro{10}(6,3)
\Label\ro{4}(0,2)
\Label\ro{7}(1,2)
\Label\ro{7}(2,2)
\Label\ro{9}(3,2)
\Label\ro{9}(4,2)
\Label\ro{11}(5,2)
\Label\ro{11}(6,2)
\Label\ro{5}(0,1)
\Label\ro{10}(1,1)
\Label\ro{10}(2,1)
\Label\ro{11}(3,1)
\Label\ro{12}(4,1)
\Label\ro{12}(5,1)
\Label\ro{12}(6,1)
\Label\ro{7}(0,0)
\Label\o{\hskip-0.5cm\oversetbrace \tsize 2c\to{\hbox{\hskip3.5cm}}}(4,6)
\Label\l{\hskip0cm\raise1.4cm\hbox{$\sideset n\and \to
     {\left\{\vbox{\vskip1.6cm}\right.}$}}(0,3)
\Label\u{\hskip0cm\undersetbrace \tsize 2p\to{\hbox{\hskip3cm}}}(4,0)
\hskip3.5cm
$$
\vskip10pt
\centerline{\eightpoint Figure 11}
}
\vskip10pt

We claim that $(2n+1)$-orthogonal tableaux $S$ of shape
$(c^{n-1},c-p)$ are in bijection with pairs
$\big(\bT,(e_1,e_2,\dots,e_p)\big)$, by a bijection $\Om_1$ say,
where $\bT$ is an $n$-tableau contained in
$\((2c)^{n-1},2c-2p\)$,
and where $(e_1,e_2,\dots,e_p)$ is a vector of numbers (we definitely
mean {\it vector\/} here, i.e., the order of the numbers is important) satisfying
$$\spreadmatrixlines{2\jot}\matrix 1\le e_1,e_2,\dots,e_p\le n,\\
 e_l\notin [i_{l}(m),j_{l}(m)]\text { for }1\le l\le p,\ 1\le m\le
s_l,\\
e_l\notin [i_{l}({s_l+1}),n]\text { for }1\le l\le p,
\endmatrix\tag5.19$$
and for all $m$,
$$
\big\vert\{1,2,\dots,m\}\backslash
\{e_l,j_{l}(1),\dots,j_{l}({s_l})\}\big\vert
\le
\big\vert\{1,2,\dots,m\}
\backslash
\{e_{l+1},i_{l+1}(1),\dots,i_{l+1}({t_{l+1}})\}\big\vert,
\tag5.20$$
given that 
$$\matrix i_{l}(1)&j_{l}(1)\\\vdots&\vdots \\ i_{l}({s_l})&j_{l}({s_l})\\
i_{l}({{s_l}+1})\\\vdots\\i_{l}({t_l})\endmatrix$$
are the $(2c-2p+2l-1)$-st and $(2c-2p+2l)$-th column of $\bT$, such that
$$(\x^{\pm1})^S=(x_1x_2\cdots x_n)^{-c}\cdot x_{e_1}\cdots x_{e_p}\cdot
\x^{\bT},\quad \text {if }(\bT,\{e_1,e_2,\dots,e_p\})=\Om_1(S).
\tag5.21$$

Again, the construction of the bijection $\Om_1$ is based on
an analysis of the orthogonal tableaux under
consideration. 
Clearly, the entries $\le n$ in $S$ form an $n$-tableau, $\bT$
say, whose shape is contained
in $\((2c)^{n-1},2c-2p\)$. 
In Figure~11 we have marked the area that
is covered by entries $\le n$ by a bold line. 
The resulting tableau is displayed in the left half of
Figure~12.

If we apply the paragraph containing (5.16) with $p=1$, then we see
that, with $e_l$ being the number that together with its conjugate
$2n+1-e_l$ does not appear in the $(2c-2p+2l-1)$-st or
$(2c-2p+2l)$-th column of $S$, $l=1,2,\dots,p$, the map $\Om_1$ defined by
$$S\to\big(\bT,(e_1,e_2,\dots,e_p)\big)$$
defines the desired bijection. In particular, the fact that also the
entries $>n$ in $S$ form a (skew) tableau is reflected by condition
(5.20). 
For, the entries $>n$ from the $(2c-2p+2l)$-th column
of $S$ are 
$$\{n+1,n+2,\dots,2n\}\backslash
\{2n+1-e_l,2n+1-j_{l}(1),\dots,2n+1-j_{l}({s_l})\},$$
and the entries $>n$ from the 
$(2c-2p+2l+1)$-st column of $S$ are
$$\{n+1,n+2,\dots,2n\}\backslash
\{2n+1-e_{l+1},2n+1-i_{l+1}(1),\dots,2n+1-i_{l+1}({t_{l+1}})\}.$$
That all entries $>n$ from the $(2c-2p+2l)$-th column
of $S$ are less or equal than their right neighbours from the
$(2c-2p+2l+1)$-st column of $S$, is exactly equivalent to requiring
$$\multline
\big\vert\{2n,2n-1,\dots,2n-m\}\backslash
\{2n+1-e_l,2n+1-j_{l}(1),\dots,2n+1-j_{l}({s_l})\}\big\vert\\
\le
\big\vert\{2n,2n-1,\dots,2n-m\}\hskip5cm\\
\backslash
\{2n+1-e_{l+1},2n+1-i_{l+1}(1),\dots,2n+1-i_{l+1}({t_{l+1}})\}\big\vert
\endmultline
$$
for all $m$, which is clearly equivalent to (5.20).

Again, it is easy to
check that the weight property (5.21) holds under this correspondence.
Our running example in Figure~11 is mapped under $\Om_1$ to the pair
in Figure~12.
\vskip10pt
\vbox{
$$
\(\hbox{\quad }
\Einheit0.5cm
\Pfad(0,-2),222221111111\endPfad
\Pfad(0,-2),122111211122\endPfad
\PfadDicke{.5pt}
\Pfad(0,-1),1\endPfad
\Pfad(0,0),1\endPfad
\Pfad(0,1),1111\endPfad
\Pfad(0,2),1111111\endPfad
\Pfad(1,0),222\endPfad
\Pfad(2,0),222\endPfad
\Pfad(3,0),222\endPfad
\Pfad(4,1),22\endPfad
\Pfad(5,1),22\endPfad
\Pfad(6,1),22\endPfad
\Label\ro{1}(0,2)
\Label\ro{1}(1,2)
\Label\ro{1}(2,2)
\Label\ro{1}(3,2)
\Label\ro{2}(4,2)
\Label\ro{3}(5,2)
\Label\ro{4}(6,2)
\Label\ro{2}(0,1)
\Label\ro{2}(1,1)
\Label\ro{2}(2,1)
\Label\ro{5}(3,1)
\Label\ro{5}(4,1)
\Label\ro{5}(5,1)
\Label\ro{5}(6,1)
\Label\ro{3}(0,0)
\Label\ro{5}(1,0)
\Label\ro{5}(2,0)
\Label\ro{6}(3,0)
\Label\ro{4}(0,-1)
\Label\ro{5}(0,-2)
\hskip3.5cm
\quad ,\quad \big(4,3,6\big)
\)
$$
\centerline{\eightpoint Figure 12}
}
\vskip10pt

\medskip
{\it Second step}. In the second step we construct a bijection
$\Om_2$ between pairs $\big(\bT,(e_1,\mathbreak
e_2,\dots,e_p)\big)$ satisfying (5.19) 
as before, and pairs $(T,\si)$, where $T$ is an
$n$-tableau whose shape is contained in $((2c)^{n})$,
and where $\si$ is a is a horizontal strip of length $p$
on the rim of $\nu$ such that the $i$-th cell of the
strip comes before the $(2c-2p+2i)$-th column,
such that
$$\x^T=x_{e_1}\cdots x_{e_p}\cdot \x^{\bT}.\tag5.22$$

Let $\big(\bT,(e_1,e_2,\dots,e_p)\big)$ be such a pair. 
We insert $e_1,e_{2},\dots,e_p$, in this order, into $\bT$,
according to a procedure that is very similar to $\Ph_2$, or $\Th_2$,
which were used in the proofs of (3.6), respectively (3.8). That we have to
modify these procedures is due to the fact that $e_l$ satisfies a
condition, namely (5.19), that depends on the $(2c-2p+2l-1)$-st and
$(2c-2p+2l)$-th column of $\bT$, and not just on the $(2c-1)$-st and
$(2c)$-th column as was the case in (5.4) or (5.14). 
Let $\bT_0:=\bT$. Suppose that,
by inserting $e_1,e_{2},\dots,e_{l}$ we already formed
$\bT_l$. Next we insert $e_{l+1}$ into $\bT_l$ in the following way.
Choose the first row (from top to bottom) of $\bT$ such that
$e_{l+1}$ is less than the entry in
the $(2c-2p+2l-1)$-st column in that row of $\bT_l$. If there is no
such row choose the first row that does not have an entry in the
$(2c-2p+2l-1)$-st column of that row. Then, starting with that row of
$\bT$, ROW-INSERT $e_{l+1}$ into $\bT$, see the definition of $\Ph_2$
in the second step of the proof of (3.6). Thus one obtains the
tableau $\bT_{l+1}$. Finally, set
$T=\Om_2\big((\bT,(e_1,e_2,\dots,e_p))\big):=\bT_p$. 
On the other hand, $\si$ is defined to be the vertical strip (it is
indeed a vertical strip, as will be shown in a moment) by which
the shapes of $T$ and $\bT$ differ. 
Thus, our running example from Figure~12 is mapped under $\Om_2$ to the
pair in Figure~13. There, the vertical strip $\si$ is visualized by bold
lines embedded in the shape of $T$.
\vskip10pt
\vbox{
$$
\(\hbox{\quad }
\Einheit0.5cm
\Pfad(0,-2),222221111111\endPfad
\Pfad(0,-2),121121121122\endPfad
\PfadDicke{.5pt}
\Pfad(0,-1),1\endPfad
\Pfad(0,0),111\endPfad
\Pfad(0,1),11111\endPfad
\Pfad(0,2),1111111\endPfad
\Pfad(1,-1),2222\endPfad
\Pfad(2,-1),2222\endPfad
\Pfad(3,0),222\endPfad
\Pfad(4,0),222\endPfad
\Pfad(5,1),22\endPfad
\Pfad(6,1),22\endPfad
\Label\ro{1}(0,2)
\Label\ro{1}(1,2)
\Label\ro{1}(2,2)
\Label\ro{1}(3,2)
\Label\ro{2}(4,2)
\Label\ro{3}(5,2)
\Label\ro{4}(6,2)
\Label\ro{2}(0,1)
\Label\ro{2}(1,1)
\Label\ro{2}(2,1)
\Label\ro{3}(3,1)
\Label\ro{5}(4,1)
\Label\ro{5}(5,1)
\Label\ro{5}(6,1)
\Label\ro{3}(0,0)
\Label\ro{4}(1,0)
\Label\ro{5}(2,0)
\Label\ro{5}(3,0)
\Label\ro{6}(4,0)
\Label\ro{4}(0,-1)
\Label\ro{5}(1,-1)
\Label\ro{6}(2,-1)
\Label\ro{5}(0,-2)
\hskip3.5cm
\hbox{\quad ,\quad }
\PfadDicke{.3pt}
\Pfad(0,-2),222221111111\endPfad
\Pfad(0,-2),121121121122\endPfad
\PfadDicke{1.5pt}
\Pfad(1,-1),112\endPfad
\Pfad(1,-1),211\endPfad
\Pfad(4,0),12\endPfad
\Pfad(4,0),21\endPfad
\PfadDicke{.3pt}
\Pfad(0,-1),1\endPfad
\Pfad(0,0),111\endPfad
\Pfad(0,1),11111\endPfad
\Pfad(0,2),1111111\endPfad
\Pfad(1,-1),2222\endPfad
\Pfad(2,-1),2222\endPfad
\Pfad(3,0),222\endPfad
\Pfad(4,0),222\endPfad
\Pfad(5,1),22\endPfad
\Pfad(6,1),22\endPfad
\hskip3.5cm
\hbox{\quad }
\)
$$
\centerline{\eightpoint Figure 13}
}
\vskip10pt

It is obvious that (5.22) holds under this correspondence.
Moreover, it is immediate from the definition of $\si$ as a result of the above
insertion procedure, that the $i$-th cell of $\si$
comes before the $(2c-2p+2i)$-th column. However, it is not so
immediate that $\si$ is indeed a
vertical strip. This will be established next. 
We will be done if we are able to 
show that later ``insertion
paths" are strictly to the right of previous ones. It suffices to
consider two successive insertion paths. 

Let the columns $2c-2p+2l-1$, $2c-2p+2l$, $2c-2p+2l+1$, $2c-2p+2l+2$
be given by the four columns in Figure~14 
(ignore for the moment `$e_l\to$' and `$e_{l+1}\to$')
\vskip10pt
\vbox{
$$\matrix &i_{l}(1)&j_{l}(1)\\&\vdots&\vdots \\
&\vdots&\vdots \\&\vdots&\vdots \\
e_l\to\kern-10pt& i_{l}({u})&j_{l}({u})\\
&\vdots&\vdots \\&\vdots&\vdots \\ & i_{l}({s_l})&j_{l}({s_l})\\
&i_{l}({{s_l}+1})\\&\vdots\\&i_{l}({t_l})\endmatrix
\spreadmatrixlines{.3pt}
\matrix &i_{l+1}(1)&j_{l+1}(1)\\&\vdots&\vdots \\ 
e_{l+1}\to\kern-5pt&i_{l+1}({v})&j_{l+1}({v})\\
&\vdots&\vdots \\ &i_{l+1}({s_{l+1}})&j_{l+1}({s_{l+1}})\\
&\vdots\\&i_{l+1}({t_{l+1}})\\
&\vphantom{h}&\vphantom{} \\
&\vphantom{\vdots}&\vphantom{\vdots} \\
&\vphantom{\vdots}&\vphantom{\vdots} \\
&\vphantom{_{l}({{s_l}+1})}
\endmatrix$$
\centerline{\eightpoint Figure 14}
}
\vskip10pt
\noindent
Suppose that the insertion of $e_l$ would start in row $u$, as is
symbolized by $e_l\to i_{l}({u})$ in Figure~14, and that the
subsequent insertion of $e_{l+1}$ would start in row $v$, as is
symbolized by $e_l\to i_{l+1}({v})$ in Figure~14. This would mean
that we have 
$$j_{l}(u-1)<e_l< i_{l}({u})\tag5.23$$ 
and
$$j_{l+1}(v-1)<e_{l+1}< i_{l+1}({v}).\tag5.24$$ 
During insertion of $e_{l+1}$, first the element $e_{l+1}$ bumps
$i_{l+1}({v})$ (which is an element of the $(2c-2p+2l+1)$-st column
of $\bT$) or an element to the left of $i_{l+1}(v)$ in the same row. 
Let us further suppose that until row $w-1$, $w\ge v$, elements of the
$(2c-2p+2l+1)$-st column are bumped, i.e., $i_{l+1}({v})$ bumps
$i_{l+1}({v+1})$, \dots, finally $i_{l+1}({w-2})$ bumps
$i_{l+1}({w-1})$. This would mean 
$$j_{l}({v+1})\le i_{l+1}({v}), \dots, 
j_{l}({w-1})\le i_{l+1}({w-2}).\tag5.25$$
We suppose that then in row $w$ the insertion path jumps to the left of the
$(2c-2p+2l+1)$-st column (note that the case $w=v$ covers the case that
$e_{l+1}$ bumps an element to the left of $i_{l+1}(v)$), which means that
$$i_{l+1}(w-1)< j_l(w),\text {\quad or }e_{l+1}<j_l(v)\text { in case }w=v.\tag5.26$$
We do not care what happens afterwards. 

We claim that
$u\le w$ and that $e_l\le i_{l+1}({w-1})$, respectively $e_l\le
e_{l+1}$ in case $w=v$. This would
imply that in the $u$-th row 
the insertion path caused by $e_{l+1}$ is strictly to the
right of the insertion path caused by $e_{l}$ and therefore has to
stay strictly to the right from thereon, by an elementary property of
ROW-INSERT. And this is what we want to show.

For proving the claim we consider (5.20) with $m=i_{l+1}(w-1)$,
respectively $m=e_{l+1}$ in case $w=v$. For
this choice of $m$, the right-hand side of (5.20) equals
$i_{l+1}(w-1)-w$ since we have
$$i_{l+1}(1)<i_{l+1}(2)<\dots<i_{l+1}(w-1),$$
and by (5.24),
$$e_{l+1}<i_{l+1}({v})<\dots<i_{l+1}({w-1}),$$
respectively equals $e_{l+1}-v$ in case $w=v$ since by (5.24) we have
$$
i_{l+1}(1)<\dots<i_{l+1}(v-1)\le j_{l+1}(v-1)
<e_{l+1}<i_{l+1}(v).
$$
Hence, the left-hand side of (5.20) is bounded above by
$i_{l+1}(w-1)-w$, respectively $e_{l+1}-v$ in case $w=v$.

On the other hand, the left-hand side of (5.20) is at least
$i_{l+1}(w-1)-w$, respectively $e_{l+1}-v$ in case $w=v$, and equal
to the lower bound only if $e_l\le i_{l+1}(w-1)$, respectively
$e_l\le e_{l+1}$ in case $w=v$. For, by (5.25) and (5.26) we have
$$
j_{l}(1)<\dots<j_{l}(w-1)\le i_{l+1}(w-2)
<i_{l+1}(w-1)<j_{l}(w),
$$
and in case $w=v$ we have by (5.24) and (5.26)
$$
j_{l}(1)<\dots<j_{l}(v-1)\le j_{l+1}(v-1)
<e_{l+1}<j_{l}(v).
$$
So the lower bound can only be reached if $e_l\le
i_{l+1}(w-1)$, respectively $e_l\le
e_{l+1}$ in case $w=v$, which is what we wanted to show.

Summarizing, we have shown that indeed
$e_l\le i_{l+1}(w-1)$, respectively $e_l\le e_{l+1}$ in case $w=v$.
Combining  this with
(5.23) and (5.26), we obtain the inequality chain 
$$j_l(u-1)<e_l\le i_{l+1}(w-1)\text { (respectively }e_{l+1})<j_l(w),$$
hence $j_l(u-1)<j_l(w)$. Since columns are strictly
increasing, this immediately implies $u\le w$, as desired.

To show that $\Ps_2$ is a bijection, we have to construct the inverse
mapping. Experienced with three other similar proofs, this is rather
straight-forward. Let $(T,\si)$ be a pair, where $T$ is an
$n$-tableau whose shape is contained in $((2c)^{n})$,
and where $\si$ is a horizontal strip of length $p$
on the rim of $\nu$ such that the $i$-th cell of the
strip comes before the $(2c-2p+2i)$-th column. For $l=1,2,\dots,p$
start a ROW-DELETE with the entry of $T$ that is located in the cell
corresponding to the $l$-th cell of $\si$ (counted from the right), but stop
before an entry in the $(2c-2l+2)$-th column would be bumped. This
procedure reverses the algorithm $\Om_2$, step by step.

\medskip
The composition $\Om_2\circ\Om_1$ is by definition the desired
bijection between $(2n+1)$-orth\-og\-on\-al tableaux $S$ of shape $(c^{n-1},c-p)$ and
pairs $(T,\si)$, where $T$ is an
$n$-tableau whose shape is contained in $((2c)^{n})$,
and where $\si$ is a horizontal strip of length $p$
on the rim of $\nu$ such that the $i$-th cell of the
strip comes before the $(2c-2i+2)$-nd column.
From (5.21) and (5.22) the weight property (5.18) follows immediately.
This completes the proof of (3.10).\quad \quad \qed
\enddemo

\subhead 6. Proof of Theorem~3\endsubhead
In this section we use Littelmann's extension
\cite{\LitPAA} of the Littlewood--Richardson rule
to symplectic and special orthogonal characters. 
This extension is described in Section~A6 of the Appendix. 
We remark that while
there are other rules  
for the decomposition of the product of two
symplectic or orthogonal characters involving ordinary
Littlewood--Richardson coefficients (the {\it Newell--Littlewood
rules}, 
see \cite{\KingAC, Theorem~4.1;
\KoTeAA, Cor.~2.5.3/Prop.~2.5.2; \SunaAE, Theorem~5.3}), these do
not appear to be very helpful for our purposes. Again, this is because they
involve modification rules for characters, these cause alternating signs,
and these in turn cause a lot of cancellations, and all
this is simply not tractable for the applications that we have in mind.

Before we move on to the proofs of (3.13)--(3.20) itself, in 
Proposition~1 we supply
decomposition formulas for the product of a rectangularly shaped
symplectic, respectively special orthogonal, character and an {\it
arbitrarily\/} shaped character of the same type. All these
expansions involve slightly modified Littlewood--Richardson
coefficients,
which even reduce to ordinary Littlewood--Richardson coefficients in a
number of cases, see the Remark after Proposition~1. But there are no
alternating signs here, and hence there is no cancellation. The
formulas (3.13)--(3.20) then follow rather easily from (6.1)--(6.5).
\proclaim{Proposition 1} For any nonnegative integer $c$ and any
partition $\la$ with at most $n$ parts there holds
$$
\sp_{2n}\big((c^n);\x^{\pm1}\big)\cdot
\sp_{2n}\big(\la;\x^{\pm1}\big)=
\sum _{\nu\subseteq \la+(c^n)} ^{}\sp_{2n}(\nu;\x^{\pm1})
\underset \mu\text { even}
\to{\sum _{\mu\subseteq((2c)^n)} ^{}}\bLR_{\la,\mu}^{\nu+(c^n)}(c)
\tag6.1$$
(`$\mu$ even' means that all the rows of $\mu$ are even), 
where $\bLR_{\la,\mu}^{\nu+(c^n)}(c)$ is the number of LR-fillings $F$
of shape $(\nu+(c^n))/\la$ with content $\mu$ and with the additional
property that 
$$\matrix \vbox {\hsize11cm\noindent
if there is an entry $e$ in the $n$-th row of $F$, in column $j$
say (see Section~A1 in the Appendix how columns are counted), 
then $F$ must contain at least $2c-2j+1$ other entries $e$
to the right of column $j$.}\endmatrix
\tag6.2$$

Next, for any nonnegative integer or half-integer $c$ and any
partition or
half-par\-ti\-tion $\la$ with at most $n$ parts there holds
$$
\so_{2n+1}\big((c^n);\x^{\pm1}\big)\cdot
\so_{2n+1}\big(\la;\x^{\pm1}\big)=
\sum _{\nu\subseteq \la+(c^n)} ^{}\so_{2n+1}(\nu;\x^{\pm1})
{\sum _{\mu\subseteq((2c)^n)} ^{}}\bLR_{\la,\mu}^{\nu+(c^n)}(c),
\tag6.3$$
with the understanding that $\nu$ ranges over 
partitions if $\la+(c^n)$ is
a partition and over half-partitions if $\la+(c^n)$ is a half-partition, 
and
where $\bLR_{\la,\mu}^{\nu+(c^n)}(c)$ is defined as before. (Again,
see Section~A1 in the Appendix how columns are counted; 
in particular, in condition (6.2) the column index
$j$ ranges over the integers if $\la$ is a partition and 
over half-integers if $\la$ is a half-partition.)

Finally, for any nonnegative integer or half-integer $c$ and any
$(2n)$-orthogonal partition or
half-partition $\la$ with at most $n$ parts there holds
$$\multline
\so_{2n}\big((c^n);\x^{\pm1}\big)\cdot
\so_{2n}\big(\la;\x^{\pm1}\big)
\\=
\sum _{\nu\subseteq \la+(c^n)} ^{}\so_{2n}(\nu;\x^{\pm1})
\bigg(\kern-10pt
\underset \oddcols\!\big(((2c)^n)/\mu\big)=0\to
{\sum _{\mu\subseteq((2c)^n)} ^{}}\kern-10pt
\tLR_{\la,\mu}^{\nu+(c^n)}(c)\bigg),
\endmultline\tag6.4$$
where $\tLR_{\la,\mu}^{\nu+(c^n)}(c)$ is the number of LR-fillings $F$
of shape $(\nu+(c^n))/\la$ with content $\mu$ and with the additional
property that for $\ell=1,2,\dots,2c$ holds:
$$\matrix \vbox {\hsize11cm\noindent
If the subfilling that arises from $F$ by deleting the
rightmost
$2c-\ell$ entries 1, the rightmost $2c-\ell$ entries 2,
\dots, the rightmost
$2c-\ell$ entries $n$ (if there are less than $2c-\ell$
entries of some size
then delete all of these) has shape $\nu(\ell)/\la$ then
$\nu(\ell)_{n-1}+\nu(\ell)_n\ge\ell$.}
\endmatrix
\tag6.5$$
Again, the sum in (6.4) is understood to range over $n$-orthogonal 
partitions if $\la+(c^n)$ is
an $n$-orthogonal partition and over $n$-orthogonal half-partitions if
$\la+(c^n)$ is an 
$n$-orthogonal half-partition.

\endproclaim

\remark{Remark} It should be noted that in case $\la_n\ge c$ the
condition (6.2) is void, so that the coefficients 
$\bLR_{\la,\mu}^{\nu+(c^n)}(c)$
which appear in (6.1) and (6.3) reduce to the {\it ordinary\/}
Littlewood--Richardson coefficients $\LR_{\la,\mu}^{\nu+(c^n)}$.
Likewise, if $\la_{n-1}+\la_n\ge 2c$ the condition (6.5) is void,
so that the coefficients $\tLR_{\la,\mu}^{\nu+(c^n)}(c)$
which appear in (6.4) reduce to the {\it ordinary\/}
Littlewood--Richardson coefficients $\LR_{\la,\mu}^{\nu+(c^n)}$.
\endremark

\demo{Proof of Proposition~1} For convenience, we start with the
proof of (6.3).

\medskip
{\it Proof of {\rm (6.3)}}. By (A.10) with
$\chi_n(\,.\,)=\so_{2n+1}(\,.\,;\x^{\pm1})$ 
we know that 
$$
\so_{2n+1}\big((c^n);\x^{\pm1}\big)\cdot
\so_{2n+1}\big(\la;\x^{\pm1}\big)=
\sum _{T} ^{}\so_{2n+1}(\la+\con(T);\x^{\pm1}),
\tag6.6$$
where the sum is over all $(2n+1)$-orthogonal tableaux $T$ of shape
$(c^n)$ such that for all $\ell=1,2,\dots,2c$ the vector
$\nu(\ell):=\la+\con(T(\ell))$ is in the Weyl chamber (A.8) of type $B$,
i.e., satisfies
$$\nu(\ell)_1\ge\nu(\ell)_2\ge\dots\ge \nu(\ell)_n\ge0.\tag6.7$$
The content $\con(T)$ of $T$ is defined after (A.3).

By comparing (6.6) with (6.3) we see that (6.3) will be proved once
we construct a bijection, $\Up$ say, between $(2n+1)$-orthogonal
tableaux $T$ of shape $(c^n)$ and content $\rh$ which satisfy (6.7)
for $\ell=1,2,\dots,2c$ and LR-fillings
$F$ of shape $(\la+\rh+(c^n))/\la$ which satisfy property (6.2).

The bijection $\Up$ is defined as follows. Let $\la$ be fixed and let
$T$ be a $(2n+1)$-orthogonal
tableau of shape $(c^n)$ with content $\rh$. By Observation~2 in
Section~A4 of the Appendix,
$(2n+1)$-orthogonal
tableaux of shape $(c^n)$ are nothing else but $(2n)$-tableaux of
shape $((2c)^n)$ where each column contains one of $e$ or $2n+1-e$
for all $e=1,2,\dots,n$. An example with $n=5$ and $c=5/2$ is
displayed in the left half of Figure~15. It satisfies the required
property that $\la+\con(T(\ell))$ is in the Weyl chamber of type $B$,
$\ell=1,2,\dots,2c$ for $\la=(4,4,3,1,1)$, as is easily checked.
\vskip10pt
\vbox{
$$
\Einheit0.5cm
\Pfad(0,-2),2222211111\endPfad
\Pfad(0,-2),1111122222\endPfad
\PfadDicke{.5pt}
\Pfad(0,-1),11111\endPfad
\Pfad(0,0),11111\endPfad
\Pfad(0,1),11111\endPfad
\Pfad(0,2),11111\endPfad
\Pfad(1,-2),22222\endPfad
\Pfad(2,-2),22222\endPfad
\Pfad(3,-2),22222\endPfad
\Pfad(4,-2),22222\endPfad
\Label\ro{1}(0,2)
\Label\ro{1}(1,2)
\Label\ro{1}(2,2)
\Label\ro{4}(3,2)
\Label\ro{6}(4,2)
\Label\ro{2}(0,1)
\Label\ro{3}(1,1)
\Label\ro{4}(2,1)
\Label\ro{6}(3,1)
\Label\ro{7}(4,1)
\Label\ro{3}(0,0)
\Label\ro{4}(1,0)
\Label\ro{5}(2,0)
\Label\ro{8}(3,0)
\Label\ro{8}(4,0)
\Label\ro{5}(0,-1)
\Label\ro{5}(1,-1)
\Label\ro{8}(2,-1)
\Label\ro{9}(3,-1)
\Label\ro{9}(4,-1)
\Label\ro{7}(0,-2)
\Label\ro{9}(1,-2)
\Label\ro{9}(2,-2)
\Label\ro{10}(3,-2)
\Label\ro{10}(4,-2)
\hskip2.5cm
\hbox{\hskip1cm $\longleftrightarrow$\hskip1cm }
\PfadDicke{1pt}
\Pfad(0,-2),122112122\endPfad
\Pfad(0,-2),222221111111\endPfad
\Pfad(0,-2),111122122112\endPfad
\PfadDicke{.5pt}
\Pfad(1,-1),111\endPfad
\Pfad(3,0),1\endPfad
\Pfad(4,1),1\endPfad
\Pfad(4,2),1\endPfad
\Pfad(2,-2),22\endPfad
\Pfad(3,-2),22\endPfad
\Pfad(4,0),2\endPfad
\Pfad(5,2),2\endPfad
\Pfad(6,2),2\endPfad
\Label\ro{1}(4,2)
\Label\ro{1}(5,2)
\Label\ro{1}(6,2)
\Label\ro{2}(4,1)
\Label\ro{2}(3,0)
\Label\ro{3}(4,0)
\Label\ro{1}(1,-1)
\Label\ro{2}(2,-1)
\Label\ro{3}(3,-1)
\Label\ro{3}(1,-2)
\Label\ro{4}(2,-2)
\Label\ro{4}(3,-2)
\hskip3.5cm
$$
\centerline{\eightpoint Figure 15}
}
\vskip10pt

To obtain the image of $T$ under $\Up$, we construct a
sequence $F_0,F_1,\dots,F_{2c}$ of fillings by reading $T$
column-wise, from right to left. The desired filling $F$ will then be
defined to be the last filling, $F_{2c}$. Define $F_0$ to be the only
filling of the shape $\la/\la$ (which is of course the empty
filling). Suppose that we already constructed $F_\ell$. To obtain
$F_{\ell+1}$, we add for $i=1,2,\dots,n$ an entry $e$ to row $i$ of
$F_\ell$ if $i$ is an entry in the $(2c-\ell)$-th column and
the $e$-th row of $T$. As already announced, we define $F$ to be
$F_{2c}$. Thus, with $\la=(4,4,3,1,1)$, our tableau in the
left half of Figure~15 is mapped by $\Up$ to the filling in the right
half of Figure~15.

It is straight-forward from this construction that the mapping $\Up$
can be reversed, step by step. So we shall be done if we show that,
given that $T$ is mapped to $F$ by $\Up$, $T$ is a (2n+1)-orthogonal
tableau of shape $(c^n)$ with content $\rh$ if and only if $F$ is a
LR-filling of shape $(\la+\rh+(c^n))/\la$ satisfying (6.2). We
provide the details only for the forward implication. Since the
arguments for the backward implication are similar, the reader will
have no difficulties to fill in the respective details.

Let $T$ be a (2n+1)-orthogonal
tableau of shape $(c^n)$ with content $\rh$ that is mapped by $\Up$
to the filling $F$. What we have to show is that $F$ is an
$n$-tableau, i.e., that entries are weakly increasing along rows
and strictly increasing along columns, that the LR-condition holds, that
the shape of $F$ is $(\la+\con(T)+(c^n))/\la$, and that $F$
satisfies (6.2). The reader is advised to keep the
example of Figure~15 in mind. It will help to follow the subsequent
arguments.

For the first statement, let $e$ and $f$ be entries in the $i$-th row
of $F$, $e$ being the left neighbour of $f$. Then, by construction of
$\Up$, $e$ was caused by some entry $i$ in the $e$-th row and
$c_e$-th column, say, of $T$, while $f$ 
was caused by some entry $i$ in the $f$-th row and
$c_f$-th column, say, of $T$. Since $f$ is to the right of $e$, $f$
was added ``later", hence $c_f<c_e$. Now, since $T$ is a tableau,
the entry $i$ in column $c_f$ cannot be higher than the entry $i$
in column $c_e$. Therefore we have $e\le f$. This holds for any
left-right neighbours in $F$, so rows are weakly increasing, as
desired. To see that columns of $F$ are strictly increasing, we
consider entries $e$ and $f$ in the same column of $F$, $e$ being the
top neighbour of $f$. Let $e$ be located in the $i$-th row of $F$
(and so $f$ be located in the $(i+1)$-st row of $F$). 
Then, by construction of
$\Up$, $e$ was caused by some entry $i$ in the $e$-th row and
$c_e$-th column, say, of $T$, while $f$ 
was caused by some entry $(i+1)$ in the $f$-th row and
$c_f$-th column, say, of $T$. It is easily seen by induction that 
for all $\ell$ the shape of the partial filling
$F_\ell$ is given by $(\la+\con(T(\ell))+((\ell/2)^n))/\la$. In
particular, because of (6.7), this implies that the ``outer shape" of
$F_\ell$ is always ``well-behaved" in the sense that lower rows
always terminate earlier than higher rows. Therefore the entry $f$ of
$F$ was added ``later" then the entry $e$, which means that the column
$c_f$ of $T$ must be weakly to the left of column $c_e$. We already know that
there is an entry $i$ in column $c_e$ and an entry $(i+1)$ in column
$c_f$. Since column $c_f$ is located weakly to the left of column
$c_e$, the entry $(i+1)$ must be in a lower row than the entry $i$.
As the entry $(i+1)$ is located in the $f$-th row and the entry $i$
is located in the $e$-th row, this means nothing else than $e<f$.
This holds for any
top-bottom neighbours in $F$, so columns are strictly increasing, as
desired.

Now we turn to the LR-condition. We have to show that, while reading
the entries of $F$ row-wise from top to bottom, and in each row from
right to left, at any stage we have 
$$\text {number of $1$'s}\ge\text{number of $2$'s}\ge\text{number of
$3$'s}\ge\cdots.$$
Now, by construction of $\Up$, each entry in $F$ corresponds to some
entry in $T$. Thus, to the above described reading of the entries of
$F$ there corresponds the following reading of the entries of $T$:
First read the entries $1$ in $T$, from left to right, then the
entries $2$, from left to right, then the entries $3$, etc. The
LR-condition is equivalent to saying that at any stage during this
reading of $T$ the number of entries read from the first row is
greater or equal the number of entries read from the second row,
which in turn is greater or equal the number of entries read from the
third row, etc. But this is obviously true because $T$ is a tableau.

What regards the shape, we already noticed that for all $\ell$
the shape of $F_\ell$ is given by $(\la+\con(T(\ell))+((\ell/2)^n))/\la$.
Hence $F=F_{2c}$ has shape $(\la+\con(T)+(c^n))/\la$. 

Finally, we
want to show that $F$ satisfies (6.2). Let $e$ be a fixed entry in
the $n$-th row and $j$-th column of $F$. Write again
$\nu(\ell)=\la+\con(T(\ell))$. By assumption, $\nu(\ell)$ lies in the
Weyl chamber of type $B$, see (6.7). In particular, we have
$\nu(\ell)_n\ge0$. Now, we already observed that 
$$(\la+\con(T(\ell))+((\ell/2)^n))/\la=(\nu(\ell)+((\ell/2)^n))/\la$$
is the shape of $F_\ell$. Hence, the condition $\nu(\ell)_n\ge0$ is
equivalent to saying that $F_\ell$ contains an entry in the $n$-th
row that is located in the $(\ell/2)$-th, respectively 
$(\ell+1)/2$-th, column (depending on whether $\la$ is a partition or
half-partition). (In
passing, we note that the $\ell=2c$ case of this fact implies that
$F=F_{2c}$ contains an entry in the $c$-th, respectively
$(c+1/2)$-th, column. Hence the shape of
$F$ can indeed be written in the form $(\nu+(c^n))/\la$ with $\nu$ a
partition or half-partition.) 
On the other hand, $F_\ell$ is the subfilling of
$F$ that arises by deleting the rightmost $2c-\ell$ entries
$1$, the rightmost $2c-\ell$ entries $2$, \dots, {\it $2c-\ell$ entries
$e$}, etc., from $F$. Now,
suppose that the fixed entry $e$ in the $n$-th row and $j$-th column
is also contained in $F_\ell$, i.e., $j\le (\ell+1)/2$. Then
there must be at least $2c-\ell$ entries $e$ in $F$ to the right of
the $j$-th column. This last property holds for all $\ell$ with $j\le
(\ell+1)/2$. Hence there must be at least $2c-(2j-1)$ 
entries $e$ in $F$ to the right of the $j$-th column. This is exactly
what we wanted to show.

Thus, the proof of (6.3) is complete.

\medskip

{\it Proof of {\rm (6.1)}}. By (A.10) with
$\chi_n(\,.\,)=\sp_{2n}(\,.\,;\x^{\pm1})$ 
we know that 
$$
\sp_{2n}\big((c^n);\x^{\pm1}\big)\cdot
\sp_{2n}\big(\la;\x^{\pm1}\big)=
\sum _{T} ^{}\sp_{2n}(\la+\con(T);\x^{\pm1}),
\tag6.8$$
where the sum is over all $(2n)$-symplectic tableaux $T$ of shape
$(c^n)$ such that for all $\ell=1,2,\dots,2c$ the vector
$\nu(\ell):=\la+\con(T(\ell))$ is in the Weyl chamber (A.8) of type $C$,
i.e., satisfies (6.7). (Recall that the Weyl chambers of types $B$ and
$C$ are the same.) 
By comparing (6.8) with (6.1) we see that (6.1) will be proved once
we construct a bijection between $(2n)$-symplectic
tableaux $T$ of shape $(c^n)$ with content $\rh$ and LR-fillings
$F$ of shape $(\la+\rh+(c^n))/\la$ with even content 
which satisfy property (6.2).

As bijection we can take the mapping $\Up$ from the preceding proof of
(6.3). We only have to observe (see Observation~3 in Section~A3 of
the Appendix)
that $(2n)$-symplectic tableaux $T$ of shape $(c^n)$
are the same as $(2n+1)$-orthogonal tableaux of
shape $(c^n)$ with the additional property that the entries $\le n$ form a
subtableau with only even rows. Suppose that $T$ is mapped by $\Up$ to
$F$. Then the length of the
$i$-th row of this subtableau of $T$ is the same as the number of occurences
of $i$ in the filling $F$. In other words, the shape of the
subtableau equals the content of the corresponding filling $F$.
Since the shape of the subtableau is even, 
the filling must have even content, as desired. 
The final observation is that since (6.7) holds here, too, the
filling must again satisfy (6.2).

\medskip
{\it Proof of {\rm (6.4)}}. By (A.10) with
$\chi_n(\,.\,)=\so_{2n}(\,.\,;\x^{\pm1})$ 
we know that 
$$
\so_{2n}\big((c^n);\x^{\pm1}\big)\cdot
\so_{2n}\big(\la;\x^{\pm1}\big)=
\sum _{T} ^{}\so_{2n}(\la+\con(T);\x^{\pm1}),
\tag6.9$$
where the sum is over all $(2n)$-orthogonal tableaux $T$ of shape
$(c^n)$ such that for all $\ell=1,2,\dots,2c$ the vector
$\nu(\ell):=\la+\con(T(\ell))$ is in the Weyl chamber (A.9) of type $D$.
By comparing (6.9) with (6.4) we see that (6.4) will be proved once
we construct a bijection between $(2n)$-orthogonal
tableaux $T$ of shape $(c^n)$ with content $\rh$ and LR-fillings
$F$ of shape $(\la+\rh+(c^n))/\la$ and content $\mu$ where
all the columns of $\mu$ have the same parity as $n$ and 
where property (6.5) is satisfied.

Again, we can take the mapping $\Up$ from the proof of
(6.3) as the bijection. Here, this is because of the observation 
(see Observation~1 in Section~A5 of the Appendix, with
$\la_{n-1}=\la_n=c$) that
$(2n)$-orthogonal tableaux $T$ of shape $(c^n)$
are the same as $(2n+1)$-orthogonal tableaux of
shape $(c^n)$ with the additional property that the entries $\le n$ form a
subtableau whose shape has only columns of the same parity as $n$.
Again, since the shape of the subtableau equals the content of the
corresponding filling, the content $\mu$ of the filling must have
columns of the same parity as $n$ throughout, as desired. This time
we have to impose (6.5) (instead of (6.2)) since
$\nu(\ell)=\la+\con(T(\ell))$ has to be in the Weyl chamber of type $D$
(and not of type $B$ or $C$).

\medskip
This completes the proof of the Proposition.\quad \quad \qed
\enddemo

Now we are in the position to prove (3.13)--(3.20).

\demo{Proof of (3.13)} We apply (6.1) with $\la=((d+1)^p,d^{n-p})$.
Because of the assumption $c\le d$, the Remark after Proposition~1 applies,
which says that the coefficients $\bLR_{\la,\mu}^{\nu+(c^n)}(c)$
reduce to ordinary Littlewood-Richardson coefficients 
$\LR_{\la,\mu}^{\nu+(c^n)}$. 
Then, obviously, (3.13) is equivalent to the claim
$$\underset \mu\text { even}\to{\sum _{\mu\subseteq((2c)^n)}
^{}}\LR_{((d+1)^p,d^{n-p}),\mu}^{\nu+(c^n)}=\cases
1&((d-c)^n)\subseteq \nu\subseteq ((c+d+1)^n)\\
&\text {and }\oddrows\!\big(\nu/((d-c)^n)\big)=p\\
0&\text {otherwise}.\endcases\tag6.10$$

Let $F$ be a LR-filling of shape $(\nu+(c^n))/((d+1)^p,d^{n-p})$ with
even content. Because of the LR-condition, almost all the entries of
$F$ are uniquely determined. To be precise, except for the entries in
column $d+1$, all the entries in the $i$-th row have to equal $i$
throughout, $i=1,2,\dots,n$, as is exemplified in Figure~16.
\vskip10pt
\vbox{
$$
\Einheit0.5cm
\PfadDicke{1pt}
\Pfad(0,0),11122221222\endPfad
\Pfad(0,0),2222222111111111\endPfad
\Pfad(0,0),1112121212212112\endPfad
\PfadDicke{.5pt}
\Pfad(3,2),1\endPfad
\Pfad(3,3),11\endPfad
\Pfad(4,4),11\endPfad
\Pfad(4,5),111\endPfad
\Pfad(4,6),111\endPfad
\Pfad(4,2),22\endPfad
\Pfad(5,3),2222\endPfad
\Pfad(6,5),22\endPfad
\Pfad(7,6),2\endPfad
\Pfad(8,6),2\endPfad
\Label\ro{1}(4,6)
\Label\ro{1}(5,6)
\Label\ro{1}(6,6)
\Label\ro{1}(7,6)
\Label\ro{1}(8,6)
\Label\ro{2}(4,5)
\Label\ro{2}(5,5)
\Label\ro{2}(6,5)
\Label\ro{3}(4,4)
\Label\ro{3}(5,4)
\Label\ro{4}(4,3)
\Label\ro{4}(5,3)
\Label\ro{5}(4,2)
\Label\ro{*}(3,1)
\Label\ro{*}(3,2)
\Label\ro{*}(3,3)
\Label\o{\hskip0cm\oversetbrace \tsize d+1\to{\hbox{\hskip2cm}}}(2,7)
\Label\l{\hskip0cm\raise2.2cm\hbox{$\sideset n\and \to
     {\left\{\vbox{\vskip1.8cm}\right.}$}}(0,3)
\Label\l{\hskip0.1cm\raise1.15cm\hbox{$\sideset p\and \to
     {\left\{\vbox{\vskip0.8cm}\right.}$}}(4,5)
\Label\l{\hskip0.1cm\raise1.2cm\hbox{$\sideset q\and \to
     {\left\{\vbox{\vskip0.9cm}\right.}$}}(3,2)
\Label\u{\hskip0.5cm\undersetbrace \tsize d\to{\hbox{\hskip1.5cm}}}(1,0)
\hskip4.5cm
$$
\centerline{\eightpoint Figure 16}
}
\vskip10pt
\noindent
However, also the entries in the $(d+1)$-st column of $F$ are
uniquely determined because the content of $F$ should be even.
Namely, $i$ occurs in the $(d+1)$-st column if and only if the number
of the other entries $i$, which is $\nu_i+c-(d+1)$, is odd. In our
example in Figure~16, the entries in the $(d+1)$-st column would have
to be $1,2,5$. Now,
suppose that $F$ contains exactly $q$ entries in the
$(d+1)$-st column, $q\le n-p$ of course. Equivalently,
$\nu=(\nu_1,\dots,\nu_{p+q},d-c,\dots,d-c)$. Then there are exactly
$q$ quantities $\nu_i+c-(d+1)$, $i\le p+q$, that are odd, plus $n-p-q$
quantities $\nu_i+c-(d+1)=-1$ for $i=p+q+1,p+q+2,\dots,n$. Therefore,
the number of odd rows in $\nu/((d+1-c)^n)$ is exactly
$q+(n-p-q)=n-p$, or equivalently, the number of odd rows in $\nu/((d-c)^n)$ 
is exactly $p$.

Summarizing, we have shown that the left-hand side in (6.10) is
different from zero only if $((d-c)^n)\subseteq \nu\subseteq
((d+c+1)^n)$ and $\oddrows\!\big(\nu/((d-c)^n)\big)=p$, the inclusions
being trivial constraints. 
In addition, we have also seen that there is exactly one
LR-filling under those conditions. Thus, (6.10), and thus also (3.13), is
established.\quad \quad \qed
\enddemo

\demo{Proof of (3.14)} Here we apply (6.1) with
$\la=(d^{n-p},(d-1)^{p})$. Now the assumption is $c\ge d$, so the
Remark after Proposition~1 does not apply. We have to show
$$\underset \mu\text { even}\to{\sum _{\mu\subseteq((2c)^n)}
^{}}\bLR_{(d^{n-p},(d-1)^{p}),\mu}^{\nu+(c^n)}(c)=\cases
1&((c-d)^n)\subseteq \nu\subseteq ((c+d)^n)\\
&\text {and }\oddrows\!\big(\nu/((c-d)^n)\big)=p\\
0&\text {otherwise}.\endcases\tag6.11$$
From (6.10) it follows directly that the left-hand side in (6.11) can
only be non-zero if $((d-c-1)^n)\subseteq \nu\subseteq ((c+d)^n)$ and
if $\oddrows\!\big(\nu/((d-1-c)^n)\big)=n-p$. Note that the latter
condition is equivalent to $\oddrows\!\big(\nu/((c-d)^n)\big)=p$ (if
this makes sense, i.e., if $((c-d)^n)\subseteq \nu$). It
also follows from (6.10) that if the left-hand hand side of (6.11)
is non-zero then it can only be 1. So, what remains to see is that it
is non-zero if and only if in addition $((c-d)^n)\subseteq \nu$, or
equivalently, $c-d\le \nu_n$.

We begin with the forward implication. Suppose that there is a
LR-filling $F$ of shape $(\nu+(c^n))/(d^{n-p},(d-1)^p)$ with even
content satisfying (6.2). 
If $c=d$ there is nothing to show. So let $c>d$. Then, by
considering the shape of $F$, we see that there must be an entry in
the $n$-th row and $(d+1)$-st column of $F$. (If there is no entry in
the $(d+1)$-st column of $F$, then $\nu_n+c\le d$. So, $\nu_n\le
d-c<0$, which is impossible since $\nu$ has to be a partition.)
Clearly, this entry
equals $n$ since it is located in the last row of a column of length
$n$. Now, condition (6.2) applied to this entry says that there must
be $2c-2(d+1)+1=2c-2d-1$ more entries $n$ to the right of this
column. All of them necessarily have to be in the $n$-th row of $F$,
hence $\nu_n+c\ge(d+1)+(2c-2d-1)$, or equivalently, $\nu_n\ge c-d$,
as desired.

For the backward implication, assume $\nu_n\ge c-d$. We have to
establish the existence of a LR-filling
of shape $(\nu+(c^n))/(d^{n-p},(d-1)^p)$ with even
content satisfying (6.2). This LR-filling can only be the uniquely
determined filling that was described in the proof of (3.13). 
In fact, the arguments of the first paragraph of this proof and those
in the proof of (3.13)
actually show that this uniquely determined LR-filling 
satisfies all the required properties, except for possibly 
(6.2) for the entry in the $d$-th column. 

We now verify (6.2) for this entry, $e$ say, by distinguishing between too
cases. First let $e\ne n$. Condition (6.2) would require that there
are $2c-2d+1$ more entries $e$ to the right. All these entries
necessarily have to be located in the $e$-th row.  Hence, condition
(6.2) would require $\nu_e+c-d\ge 2c-2d+1$. Now, because of $\nu_n\ge
c-d$, we have $\nu_e+c-d\ge \nu_n+c-d\ge 2c-2d$. So the total number
of $e$'s, which is $\nu_e+c-d+1$, is at least $2c-2d+1$. Since the
content of the filling is even, this number must be even. So it is
actually at least $2c-2d+2$. Hence, $\nu_e+c-d+1\ge 2c-2d+2$, or
equivalently, $\nu_e+c-d\ge 2c-2d+1$, as required.

Now let $e=n$. Condition (6.2) would require that there are $2c-2d+1$
more entries $n$ to the right. Now, the total number of $n$'s (all of
them are located in the $n$-th row) equals $\nu_n+c-d+1$. Because of
$\nu_n\ge c-d$ this number is at least $2c-2d+1$. Again, since the
content of the filling is even, this number must be even. So it is
actually at least $2c-2d+2$. Hence, there are at least $2c-2d+1$ more
$n$'s to the right of the $n$ in the $d$-th column, as required.

This completes the proof of (3.14).\quad \quad \qed
\enddemo

\demo{Proof of (3.15)} We apply (6.4) with
$\la=(d^{n-1},d-p)$. 
Thus, we have to show that
$$\underset \oddcols\!\big(((2c)^n)/\mu\big)=0\to
{\sum _{\mu\subseteq((2c)^n)}
^{}}\tLR_{(d^{n-1},d-p),\mu}^{\nu+(c^n)}(c)=\cases
1&(\vert c-d\vert^{n-1},c-d)\subseteq \nu\subseteq ((c+d)^n)\\
&\text {and }\oddcols\!\big(((c+d)^n)/\nu\big)=p\\
0&\text {otherwise}.\endcases\tag6.12$$

Let $F$ be a LR-filling of shape $(\nu+(c^n))/(d^{n-1},d-p)$ with
content $\mu$, where $\oddcols\!\big(((2c)^n)/\mu\big)=0$. For later
use we note right here 
that the left-hand side of (6.12) can only be non-zero if
$$\nu\subseteq((c+d)^n).\tag6.13$$

Similarly here, because of the LR-condition, the entries in the first $n-1$ rows are
uniquely determined. To be precise, all the entries in the $i$-th row
have to equal $i$, $i=1,2,\dots,n-1$, as is exemplified in Figure~17.
\vskip10pt
\vbox{
$$
\Einheit0.4cm
\PfadDicke{1pt}
\Pfad(0,0),1121111122222\endPfad
\Pfad(0,0),2222221111111111111\endPfad
\Pfad(0,0),1111112112212112112\endPfad
\PfadDicke{.5pt}
\Pfad(7,2),1\endPfad
\Pfad(7,3),1\endPfad
\Pfad(7,4),11\endPfad
\Pfad(7,5),1111\endPfad
\Pfad(3,0),2\endPfad
\Pfad(4,0),2\endPfad
\Pfad(5,0),2\endPfad
\Pfad(6,0),2\endPfad
\Pfad(8,3),222\endPfad
\Pfad(9,4),22\endPfad
\Pfad(10,4),22\endPfad
\Pfad(11,5),2\endPfad
\Pfad(12,5),2\endPfad
\Label\ro{1}(7,5)
\Label\ro{1}(8,5)
\Label\ro{1}(9,5)
\Label\ro{1}(10,5)
\Label\ro{1}(11,5)
\Label\ro{1}(12,5)
\Label\ro{2}(7,4)
\Label\ro{2}(8,4)
\Label\ro{2}(9,4)
\Label\ro{2}(10,4)
\Label\ro{3}(7,3)
\Label\ro{3}(8,3)
\Label\ro{4}(7,2)
\Label\ro{5}(7,1)
\Label\ro{*}(2,0)
\Label\ro{*}(3,0)
\Label\ro{*}(4,0)
\Label\ro{*}(5,0)
\Label\o{\hskip0.4cm\oversetbrace \tsize d\to{\hbox{\hskip2.8cm}}}(3,6)
\Label\l{\hskip0cm\raise1.2cm\hbox{$\sideset n\and \to
     {\left\{\vbox{\vskip1.3cm}\right.}$}}(0,3)
\Label\o{\hskip0.4cm\oversetbrace \tsize p\to{\hbox{\hskip2cm}}}(4,1)
\Label\u{\hskip0cm\undersetbrace \tsize q\to{\hbox{\hskip1.6cm}}}(4,0)
\hskip5.2cm
\hbox{\hskip1cm}
\Einheit0.4cm
\PfadDicke{1pt}
\Pfad(0,0),1121111122222\endPfad
\Pfad(0,0),222222111111111111111\endPfad
\Pfad(0,0),111111112112212112112\endPfad
\PfadDicke{.5pt}
\Pfad(7,1),1\endPfad
\Pfad(7,2),111\endPfad
\Pfad(7,3),111\endPfad
\Pfad(7,4),1111\endPfad
\Pfad(7,5),111111\endPfad
\Pfad(3,0),2\endPfad
\Pfad(4,0),2\endPfad
\Pfad(5,0),2\endPfad
\Pfad(6,0),2\endPfad
\Pfad(7,0),2\endPfad
\Pfad(8,1),22222\endPfad
\Pfad(9,1),22222\endPfad
\Pfad(10,3),222\endPfad
\Pfad(11,4),22\endPfad
\Pfad(12,4),22\endPfad
\Pfad(13,5),2\endPfad
\Pfad(14,5),2\endPfad
\Label\ro{1}(7,5)
\Label\ro{1}(8,5)
\Label\ro{1}(9,5)
\Label\ro{1}(10,5)
\Label\ro{1}(11,5)
\Label\ro{1}(12,5)
\Label\ro{1}(13,5)
\Label\ro{1}(14,5)
\Label\ro{2}(7,4)
\Label\ro{2}(8,4)
\Label\ro{2}(9,4)
\Label\ro{2}(10,4)
\Label\ro{2}(11,4)
\Label\ro{2}(12,4)
\Label\ro{3}(7,3)
\Label\ro{3}(8,3)
\Label\ro{3}(9,3)
\Label\ro{3}(10,3)
\Label\ro{4}(7,2)
\Label\ro{4}(8,2)
\Label\ro{4}(9,2)
\Label\ro{5}(7,1)
\Label\ro{5}(8,1)
\Label\ro{5}(9,1)
\Label\ro{*}(2,0)
\Label\ro{*}(3,0)
\Label\ro{*}(4,0)
\Label\ro{*}(5,0)
\Label\ro{*}(6,0)
\Label\ro{*}(7,0)
\Label\o{\hskip0.4cm\oversetbrace \tsize d\to{\hbox{\hskip2.8cm}}}(3,6)
\Label\l{\hskip0cm\raise1.2cm\hbox{$\sideset n\and \to
     {\left\{\vbox{\vskip1.3cm}\right.}$}}(0,3)
\Label\o{\hskip0.4cm\oversetbrace \tsize p\to{\hbox{\hskip2cm}}}(4,1)
\Label\u{\hskip0cm\undersetbrace \tsize q\to{\hbox{\hskip2.4cm}}}(5,0)
\hskip6cm
$$
\centerline{\eightpoint Figure 17}
}
\vskip10pt
\noindent
However, also the entries in the $n$-th row of $F$ are uniquely
determined because the content $\mu$ of $F$ should satisfy
$\oddcols\!\big(((2c)^n)/\mu\big)=0$. For convenience, let $\tF$
denote the subfilling of $F$ consisting of the first $n-1$ rows of
$F$. Namely, there have to be as many entries $i$ in the $n$-th row of $F$ as
there are columns in $\tF$ of length $i-1$ with parity different from
$n$, $i=1,2,\dots,n$. (On the side we note that therefore all the
entries in the $n$-th row have the same parity as $n$.) In our
examples in Figure~17, the entries in the $n$-th row would have
to be $2,2,4,6$ and $2,2,4,6,6,6$, respectively. 
Now, suppose that $F$ contains exactly $q$
entries in the $n$-th row. Then there are exactly $q$ columns in $\tF$
with parity different from $n$. If $q\le p$ then the number of
columns of $F$ whose parity is different from $n$ equals the
aforementioned $q$ columns plus the $p-q$ empty columns
$d-p+q+1,d-p+q+2,\dots,d$, see the left filling in Figure~17. 
If $q\ge p$ then the
number of columns of $F$ whose parity is different from $n$ equals
the aforementioned $q$ columns minus the $q-p$ columns
$d+1,d+2,\dots,d+q-p$ of length $n$, see the right filling in Figure~17.
Hence, in both cases the number of columns of $F$ whose parity is
different from $n$ equals $p$, or equivalently,
$\oddcols\!\big(((c+d)^n)/\nu\big)=p$. (Recall that by (6.13) the
shape $\nu$ is indeed contained in $((c+d)^n)$.)

Summarizing, we have shown that the left-hand side of (6.12) is
different from zero only if $\nu\subseteq ((c+d)^n)$ and 
$\oddcols\!\big(((c+d)^n)/\nu\big)=p$. We have also shown that if the
left-hand side is non-zero, it can only be 1 since there is at most
one LR-filling. So, what remains is to show that it is non-zero if and
only if in addition $(\vert c-d\vert^{n-1},c-d)\subseteq \nu$.

Again, we begin with the forward implication.
Suppose that there is a
LR-filling $F$ of shape $(\nu+(c^n))/(d^{n-1},d-p)$ with 
$\oddcols\!\big(((2c)^n)/\mu\big)=0$, where $\mu$ is the content of
$F$, satisfying (6.5). Because of $(d^{n-1},d-p)\subseteq \nu+(c^n)$,
we have $((d-c)^{n-1},d-c-p)\subseteq \nu$. Hence, we will be done if
we can show $c-d\le \nu_n$.

We already noted that the first $n-1$ rows of $F$ are uniquely
determined, in particular, the $(n-1)$-st row of $F$ contains only
entries $n-1$. Now we apply condition (6.5) with
$\ell=c+d-\nu_{n-1}$. (Note that in (6.5) 
this choice of $\ell$ has the effect of removing
$2c-\ell=\nu_{n-1}+c-d$ entries of each size. In particular, all the
entries from $(n-1)$-st row are removed.) Since the content $\mu$ of
$F$ satisfies $\oddcols\!\big(((2c)^n)/\mu\big)=0$, the number of
$n$'s must equal the number of $(n-1)$'s in $F$. Therefore the number
of $n$'s, all of which have to be in the $n$-th row, is at least as
large as the number of $(n-1)$'s in the $(n-1)$-st row of $F$, the
latter being $\nu_{n-1}+c-d$. Therefore condition (6.5) with
$\ell=c+d-\nu_{n-1}$ reads
$$\big(c+\nu_{n-1}-(\nu_{n-1}+c-d)\big)+
\big(c+\nu_{n}-(\nu_{n-1}+c-d)\big)\ge c+d-\nu_{n-1}.$$
Simplifying, we obtain $\nu_n\ge c-d$, which is what we wanted.

For the backward implication, assume $\nu_n\ge c-d$. Recall that
there is a uniquely determined LR-filling 
of shape $(\nu+(c^n))/(d^{n-1},d-p)$ with
content $\mu$, where\linebreak
$\oddcols\!\big(((2c)^n)/\mu\big)=0$. What has
to be shown is that it also satisfies (6.5). A careful reading of the
previous paragraph reveals that $\nu_n\ge c-d$ is actually {\it
equivalent\/} to condition (6.5) with $\ell=c+d-\nu_{n-1}$. Now, a
moment's thought will convince the reader that in our particular
situation (the $(n-1)$-st row consists of $(n-1)$'s throughout,
$\ell=c+d-\nu_{n-1}$ in (6.5) therefore empties the complete
$(n-1)$-st row), condition (6.5) holds for {\it all\/} $\ell$ if and
only if it holds for the {\it particular choice\/}
$\ell=c+d-\nu_{n-1}$. Thus, (6.5) is established.

This completes the proof of (3.15).\quad \quad \qed

\enddemo

\demo{Proof of (3.16)} Equation (3.16) is an easy corollary of
equation (3.15) with $p$ replaced by $2d-p$. 
This is due to the observation that the substitution
$x_n\to1/x_n$ in a character $\so_{2n}(\la;\x^{\pm1})$ has the effect
$$\so_{2n}((\la_1,\dots,\la_{n-1},\la_n);\x^{\pm1})\big\vert_{x_n\to1/x_n}
=\so_{2n}((\la_1,\dots,\la_{n-1},-\la_n);\x^{\pm1}),$$
which follows easily from the definition (2.12).\quad \quad \qed
\enddemo

\demo{Proof of (3.17)} We apply (6.3) with $\la=(d^{n-p},(d-1)^p)$.
Note that the only difference between (6.3) and (6.1) is that in (6.1)
the partition $\mu$ is required to be even. Hence, we can use those
arguments from the proofs of (3.13) (here one has to
replace $d$ by $d-1$ and $p$ by $n-p$) and (3.14) that do not rely on
this requirement. 

Again, we have to consider LR-fillings $F$ of shape
$(\nu+(c^n))/(d^{n-p},(d-1)^p)$ that satisfy (6.2), but with
{\it arbitrary\/} content. So, again all entries to the right of the
$d$-th column are uniquely determined, and in the same way, see
Figure~16. But the entries in the $d$-th column are {\it not\/}
unique now. In fact, there is almost complete freedom, the
constraints being that the entries along the columns have to be
strictly increasing, i.e., each entry of a particular size can only
occur once, that the total content
$\mu$ of $F$ has to be a partition of course, 
and that (6.2) must be satisfied. It is easy to see that
the first and second constraint are equivalent to
$(\mu+(d^n))/(\nu+(c^n))$ being a vertical strip of length
$p-m_{d-c-1}(\nu)$ avoiding the $d$-th column and, because of
$\mu\subseteq((2c)^n)$, the $(d+2c+1)$-st column.
Note that the latter ``avoidance" conditions give the first and third
``avoidance" condition in (3.18) when columns are counted with
respect to $\nu$, i.e., when everything is shifted back by $c$.

The inclusion 
$\big(\vert c-d\vert^{n-p},(\max\{d-c-1,c-d\})^{p}\big)
\subseteq\nu$ follows from the trivial inclusion
$((d-c)^{n-p},(d-c-1)^p)\subseteq \nu$, and from the fact that
$c-d\le \nu_n$ if $c\ge d$, which is shown in the same way as in
the proof of (3.14).

Finally, we claim that if $c\ge d$ the vertical strip
$(\mu+(d^n))/(\nu+(c^n))$ has to avoid the $(2c-d+1)$-st column. Note
that this ``avoidance" condition gives the second
``avoidance" condition in (3.18) when columns are counted with
respect to $\nu$, i.e., when everything is shifted back by $c$.
The former ``avoidance'' condition comes
from considering condition (6.2) for the entry in the $n$-th row and
$d$-th column of the filling $F$. See the analogous considerations at
the end of the proof of (3.14). The difference here is that the
content $\mu$ can be arbitrary. Hence the argument in the proof of
(3.14) that the number of $e$'s or $n$'s is even does not apply here.
We leave the details to the
reader.

This finishes the proof of (3.17).\quad \quad \qed
\enddemo

\demo{Sketch of proof of (3.19)} As in the preceding proof of (3.17) we
apply (6.3), this time with $\la=(d^{n-1},d-p)$. The arguments are
very similar here. We have to consider LR-fillings $F$ of shape
$(\nu+(c^n))/(d^{n-1},d-p)$ with {\it arbitrary\/} content. As in the
proof of (3.15), all the entries in the first $n-1$ rows of such a
filling are uniquely determined, and in the same way, see Figure~17.
But the entries in the $n$-th row are not. They can be chosen
arbitrarily as long as the LR-condition and (6.2) are satisfied. The
LR-condition implies that, with $\mu$ denoting the content of the
filling, $(\mu_1+d,\dots,\mu_{n-1}+d)/(\nu_1+c,\dots,\nu_{n-1}+c)$
is a {\it horizontal\/} strip, $\si$ say. Of course, the length of
$\si$, which equals the number of entries $\le n-1$ in the $n$-th row of the
LR-filling, is at most $\nu_n+c-d+p$, the length of the $n$-th row of the
filling. Note that this is one part of the first inequality in
(3.20). At the same time, the length of $\si$ is at most $p$. For,
there cannot be more than $p$ entries $\le n-1$ in the $n$-th row
since the entry in the $(d+1)$-st column (which is a column of length
$n$) and the $n$-th row has to be $n$. Note that this proves the
second part of the second inequality in (3.20). 
Moreover, the LR-condition applied to entries $n-1$ and $n$
implies the first part of the 
second inequality in (3.20). Finally, condition (6.2)
implies the third inequality in (3.20), and, when applied to the first
entry $n$ in the $n$-th row, that the length of $\si$ is at least
$c-d+p-\nu_n$, the latter being the missing part of the first
inequality in (3.20). We leave it to the reader to check this in
detail.

The inclusion $(\vert c-d\vert^{n-1},\max\{c-d,d-c-p\})\subseteq \nu$
follows on the one hand from the trivial inclusion
$((d-c)^{n-1},d-c-p)\subseteq\nu$, and on the other hand from $c-d\le \nu_n$,
which is derived from the first and second inequality in (3.20) as
follows: $\nu_n\ge -\vert\si\vert+c-d+p\ge c-d$. \quad \quad \qed
\enddemo

\subhead 7. Applications to plane partitions and tableaux enumeration\endsubhead
We now apply results from Section~3 to derive some enumeration
results for plane partitions of trapezoidal shape and for tableaux.
Recall \cite{\ProcAE} that the {\it trapezoidal shape\/}
$(N,N-2,\dots,N-2r+2)$ is an array of cells with $r$ rows, 
each row indented by one cell to
the right with respect to the previous row, and $N-2i+2$
cells in row $i$. A plane partition of shape $(N,N-2,\dots,N-2r+2)$
is a filling of the trapezoidal shape $(N,N-2,\dots,N-2r+2)$ with
{\it nonnegative\/} integers (note that we allow 0 as an entry) such
that entries along rows and columns are weakly decreasing.

We begin with an application of (3.3). The second statement in the
theorem below, (7.2),
is a result of Proctor \cite{\ProcAE, Corollary on p.~554}.
\proclaim{Theorem 4}The number of plane partitions of trapezoidal
shape $(N,N-2,\dots,\mathbreak
N-2r+2)$ with entries between $0$ and $c$ and where the
entries on the main diagonal form a partition with exactly $p$
columns of parity different from $r$ (equivalently, 
$\sum _{i=1} ^{r}(-1)^{r-i+1}a_{ii}=p$ if $r$ is even,
respectively $\sum _{i=2} ^{r}(-1)^{r-i+1}a_{ii}=p$ 
if $r$ is odd, 
with $a_{ii}$ denoting the first entry in row $i$), equals
$$\binom cp\frac {\dsize\binom {p+r-1}p} {\dsize\binom{N-r+c}p}
\underset 1\le j\le c\to{\prod _{1\le i\le r} ^{}}\frac {N-i+j}
{i+j-1}.\tag7.1$$
In particular, the number of plane partitions of trapezoidal
shape $(N,N-2,\dots,N-2r+2)$ with entries between $0$ and $c$ equals
$$\underset 1\le j\le c\to{\prod _{1\le i\le r} ^{}}\frac {N+1-i+j}
{i+j-1}.\tag7.2$$
\endproclaim

\demo{Proof}We use (3.3) with $x_i=1$, $i=1,2,\dots,\lfloor
N/2\rfloor$. For this choice of $x_i$'s, the Schur functions reduce
to a number (which is the dimension of the corresponding irreducible
representation of $\GL(N,\C)$) which has a nice closed form (see
\cite{\FuHaAA, Ex.~A.30.(ii); \MacdAC, I, Ex.~4 on p.~45; 
\SunaAE, Theorem~4.4}), and is
therefore easily computed. Thus, the left-hand side of (3.3) turns
into (7.1). On the right-hand side of (3.3), we have a certain sum of
symplectic characters evaluated at $x_i=1$, $i=1,2,\dots,\lfloor
N/2\rfloor$. Now, beside the $(2n)$-symplectic tableaux of
DeConcini and Procesi in Section~A3 of 
the Appendix there are other symplectic
tableaux. Namely, the (even) symplectic characters $\sp_{2n}(\la;\x^{\pm1})$
can also be described by King's \cite{\KingAD} symplectic tableaux of
shape $\la$ (see also \cite{\ProcAK, Theorem~4.2; \SunaAE, Theorem~2.3}), and
the (odd) symplectic characters $\sp_{2n+1}(\la;\x^{\pm1})$
can also be described by Proctor's \cite{\ProcAF} odd symplectic tableaux of
shape $\la$ (see also \cite{\ProcAK, Theorem~4.2 with $z=1$}). 
There is a uniform
definition. Let $N=2n$ or $N=2n+1$. Then a King/Proctor symplectic tableau of
shape $\la$ is an $N$-tableau of shape $\la$ such that the entries in
the $i$-th row are at least $2i-1$ for all $i$. Thus, the right-hand
side of (3.3) can be interpreted as the number of King/Proctor
symplectic tableaux with entries $\le N$ and of some shape $\nu$,
where $\nu$ is contained in $(c^r)$ and
$\oddcols\!\big((c^r)/\nu\big)=p$. These tableaux are now translated
into plane partitions of trapezoidal shape as described in
\cite{\ProcAB, bottom of p.~295}. Namely, given
such a King/Proctor symplectic tableau, replace each entry $e$ by
$2n+1-e$, then interpret each row of the resulting array as a
partition and replace it by its conjugate partition. Next, shift
the $i$-th row by $(i-1)$ cells to the right, $i=1,2,\dots$, to
obtain a plane partition of ``shifted" shape that is contained in the
trapezoidal shape $(N,N-2,\dots,N-2r+2)$. Finally, place a zero in each cell
of the trapezoidal shape that is not yet filled. It is easy to see
that during this transformation the lengths of the rows of a King/Proctor
symplectic tableau become the entries on the main diagonal of the
resulting plane partition of trapezoidal shape. This establishes the
first assertion of Theorem~4.

The number in (7.2) is obtained by summing the numbers in (7.1) over
all $p$, by means of the Vandermonde sum (cf\. \cite{\SlatAC,
(1.7.7)}). This completes the proof of Theorem~4.\quad \quad \qed

\enddemo
The next two theorems give applications of Theorem~2. For 
the proof of these theorems we need a few determinant evaluations
that are listed in the Lemma below.
We remark that the evaluations (7.3), (7.4), (7.5) are basically the
Weyl denominator factorizations of types $C$, $B$, $D$, respectively
(cf\. \cite{\FuHaAA, Lemma~24.3, Ex.~A.52, Ex.~A.62, Ex.~A.66}).
\proclaim{Lemma}The following identities hold true:
$$
\det_{1\le i,j\le n}(x_i^j-x_i^{-j})=(x_1\cdots x_n)^{-n}\prod
_{1\le i<j\le n} ^{}\big((x_i-x_j)(1-x_ix_j)\big)
\ \prod _{i=1} ^{n}(x_i^2-1),\tag7.3$$
\vskip-5pt
$$\multline
\det_{1\le i,j\le n}(x_i^{j-1/2}-x_i^{-(j-1/2)})
\\=(x_1\cdots
x_n)^{-n+1/2}\prod
_{1\le i<j\le n} ^{}\big((x_i-x_j)(1-x_ix_j)\big)
\ \prod _{i=1} ^{n}(x_i-1),
\endmultline\tag7.4$$
\vskip-5pt
$$
\det_{1\le i,j\le n}(x_i^{j-1}+x_i^{-(j-1)})=2\cdot(x_1\cdots x_n)^{-n+1}\prod
_{1\le i<j\le n} ^{}\big((x_i-x_j)(1-x_ix_j)\big)
,\tag7.5$$
\vskip-5pt
$$
\det_{1\le i,j\le n}(x_i^{j}+x_i^{-j})=(x_1\cdots x_n)^{-n}\prod
_{1\le i<j\le n} ^{}\big((x_i-x_j)(1-x_ix_j)\big)
\sum _{k=0} ^{n}e_k(x_1,\dots,x_n)^2,\tag7.6$$
\vskip-5pt
$$\multline
\det_{1\le i,j\le n}(x_i^{j-1/2}+x_i^{-(j-1/2)})\\=(x_1\cdots
x_n)^{-n+1/2}\prod
_{1\le i<j\le n} ^{}\big((x_i-x_j)(1-x_ix_j)\big)
\ \prod _{i=1} ^{n}(x_i+1).
\endmultline\tag7.7$$
\endproclaim

\demo{Proof}Identities (7.3)--(7.5), and (7.7) are readily proved by
the standard argument that proves Vandermonde-type determinant
evaluations. 

For (7.6) there is a little bit of work to do. First, by reversing
the order of columns, and adding some factors, we rewrite the
determinant in (7.6) as
$$(-1)^{\binom n2}\frac {\det\limits_{1\le i,j\le
n}(x_i^{n+1-j}+x_i^{-(n+1-j)})} {\det\limits_{1\le i,j\le
n}(x_i^{n-j}+x_i^{-(n-j)})}{\det\limits_{1\le i,j\le
n}(x_i^{n-j}+x_i^{-(n-j)})}.\tag7.8$$
Next we observe that because of (2.12) the quotient of determinants
is one half of the orthogonal character
$\o_{2n}\big((1,1,\dots,1);\x^{\pm1}\big)$
(recall that $\o_{2n}(\la;\x^{\pm1})$ is the sum of 
$\so_{2n}(\la;\x^{\pm1})$ and $\so_{2n}(\la^-;\x^{\pm1})$ if
$\la_n\ne0$). In
addition, we reverse the order of columns in the single determinant
in (7.8) to obtain for (7.8)
$$\frac {1} {2}\o_{2n}\big((1,1,\dots,1);\x^{\pm1}\big)\,
\det_{1\le i,j\le n}(x_i^{j-1}+x_i^{-(j-1)}).\tag7.9$$
By the ``orthogonal Jacobi--Trudi identity" \cite{\FuHaAA,
Cor.~24.45; \KoTeAA, Theorem~2.3.3, (6)}, the orthogonal character in
(7.9) is nothing else but the elementary symmetric function
$e_n(x_1,x_1^{-1},\dots,x_n,x_n^{-1})$. Moreover, the determinant in
(7.9) can be evaluated by (7.5). Thus, the expression (7.9) becomes
$$\multline
e_n(x_1,x_1^{-1},\dots,x_n,x_n^{-1})\cdot
(x_1\cdots x_n)^{-n+1}\prod
_{1\le i<j\le n} ^{}\big((x_i-x_j)(1-x_ix_j)\big).
\endmultline\tag7.10$$
The elementary symmetric function in (7.10) can be transformed as
follows:
$$\align 
e_n(x_1,x_1^{-1},\dots,x_n,x_n^{-1})&=\sum _{k=0}
^{n}e_k(x_1,\dots,x_n)\cdot e_{n-k}(x_1^{-1},\dots,x_n^{-1})\\
&=(x_1\cdots x_n)^{-1}\sum _{k=0} ^{n}e_k(x_1,\dots,x_n)^2.
\endalign$$
Plugging this into (7.10) completes the proof of (7.6).\quad \quad \qed
\enddemo
As a first application of Theorem~2 we give new proofs of two theorems of the
author \cite{\KratAP, Theorems~21 and 11} by means of (3.6). These
are refinements of the Bender--Knuth and MacMahon (ex-)conjectures,
see \cite{\KratAP, Sections~3.3, 4.3} for more information and
references. In order
to be able to formulate the Theorem, 
we have to introduce a few $q$-notations. 
We write
$[\al]_q:=1-q^\al$, $[n]_q!:=[1]_q[2]_q\dotsb [n]_q$, $[0]_q!:=1$,
and
$$\bmatrix n\\k\endbmatrix_q=\cases \dfrac{[n]_q\cdot
[n-1]_q\dotsb[n-k+1]_q}
{[k]_q!}&k\ge0\\
0&k<0\endcases\ .$$
The base $q$ in $[\al]_q$, $[n]_q!$, 
and $\[\smallmatrix n\\k\endsmallmatrix\]_q$ will in
most cases be omitted. Only if the base is different from $q$ it
will be explicitly stated.

\proclaim{Theorem 5}Let $n(T)$ denote the sum of all the entries of a
tableau $T$. The generating function $\sum q^{n(T)}$ 
for tableaux $T$ with $p$
odd rows, with at most $c$ columns, and with entries between $1$ and $n$
is given by 
$$q^{\binom {p+1}2}\frac {[2r]} {[2r+p]}\bmatrix n\\p\endbmatrix
\frac {\dsize \bmatrix n+2r\\n\endbmatrix} {\dsize\bmatrix 
n+2r+p\\n\endbmatrix}\prod _{1\le i\le j\le n} ^{}\frac {[2r+i+j]}
{[i+j]}\quad \quad \text {if\/ }c=2r\tag7.11$$
and
$$q^{\binom{p+1}2}\bmatrix n\\p\endbmatrix\prod _{1\le i\le j\le n} ^{}
\frac {[2r+i+j]} {[i+j]}\quad \quad \text {if\/ }c=2r+1.\tag7.12$$

The generating function $\sum q^{n(T)}$ for
tableaux $T$ with $p$ odd rows, with at most $c$ columns,
and with only odd entries which lie between $1$ and $2n-1$, is given by
$$\multline q^{p^2} \frac {[2r+2p]_{q^2}\,[r]_{q^2}}
{[2r+p]_{q^2}\,[r+p]_{q^2}}\bmatrix n\\p\endbmatrix_{q^2}
\frac {\bmatrix n+2r\\n\endbmatrix_{q^2}\,} {\bmatrix
n+2r+p\\n\endbmatrix_{q^2}}\\
\times\prod _{i=1} ^{n}\frac {[r+i]_{q^2}} {[i]_{q^2}}
\prod _{1\le i<j\le n} ^{}\frac {[2r+i+j]_{q^2}} {[i+j]_{q^2}}\quad \quad 
\text {if }c=2r
\endmultline\tag7.13$$
and
$$
q^{p^2}\bmatrix n\\p\endbmatrix_{q^2}
\prod _{i=1} ^{n}\frac {[r+i]_{q^2}} {[i]_{q^2}}
\prod _{1\le i<j\le n} ^{}\frac {[2r+i+j]_{q^2}} {[i+j]_{q^2}}\quad \quad 
\text {if }c=2r+1.
\tag7.14$$

\endproclaim

\demo{Proof} For proving (7.11), replace $x_i$ by $q^i$,
$i=1,2,\dots,n$, in (3.6). By (2.13), the left-hand side of (3.6) can
be written as a quotient of two determinants which for this choice of
$x_i$'s factor by means of (7.3). On the other hand, by (A.2) the
right-hand side of (3.6) with $x_i=q^i$, $i=1,2,\dots,n$, 
clearly equals the generating
function for the tableaux in the first assertion of Theorem~5.

The expression (7.12) is obtained without much effort from the $p=0$ case
of (7.11) as is described in \cite{\KratAP, Proof of (4.3.2b)}.

The proof of (7.13) proceeds in the same way as the proof of (7.11). 
To point out the
differences, instead of substituting $q^i$ now substitute $q^{2i-1}$
for $x_i$, $i=1,2,\dots,n$, and use (7.4) instead of (7.3). 
The expression (7.14) is obtained without much effort from the $p=0$ case
of (7.13) as is described in \cite{\KratAP, Proof of (3.3.3b)}.
\quad \quad \qed
\enddemo

Our second application of Theorem~2 leads to different refinements of
the Bender--Knuth and MacMahon (ex-)conjectures which are new. What
we do is to use (3.7) to obtain a ``column-analogue"
of Theorem~5.
Now, the tableaux generating functions in question do
not factor completely in terms of cyclotomic polynomials, but it is
still possible to write the results in a reasonably compact form.

\proclaim{Theorem 6}Let $n(T)$ denote the sum of all the entries of a
tableau $T$. The generating function $\sum q^{n(T)}$ for tableaux 
$T$ with at most $c$ columns, $p$ of which being odd, 
and with entries between $1$ and $n$,
is given by 
$$\multline q^p\frac {\bmatrix p+n-1\\p\endbmatrix \bmatrix
c-p+n-1\\c-p\endbmatrix} {2\bmatrix c+n-1\\c\endbmatrix} \prod _{1\le
i<j\le n} ^{}\frac {[c+i+j-2]} {[i+j]}\\
\times\bigg(\sum _{k=0} ^{n}q^{k(k+c-1)}\(q^k\bmatrix n-1\\k\endbmatrix
+q^{-p}\bmatrix n-1\\k-1\endbmatrix\)^2+
(1-q^{c-2p})\prod _{i=1} ^{n-1}(1-q^{c+2i})\bigg).\\
\endmultline\tag7.15$$

The generating function $\sum q^{n(T)}$ for
tableaux $T$ with at most $c$ columns, $p$ of which being odd,
and with only odd entries which lie between $1$ and $2n-1$, is given by
$$\multline q^p\frac {\bmatrix p+n-1\\p\endbmatrix_{q^2} \bmatrix
c-p+n-1\\c-p\endbmatrix_{q^2}} {\bmatrix c+n-1\\c\endbmatrix_{q^2}} 
\prod _{i=1} ^{n}\frac {1} {1+q^{2(i-1)}}\prod _{1\le
i<j\le n} ^{}\frac {[c+i+j-2]_{q^2}} {[i+j-2]_{q^2}}\\
\times\bigg((1+q^{c-2p})\prod _{i=1}
^{n-1}(1+q^{c+2i})+
(1-q^{c-2p})\prod _{i=1} ^{n-1}(1-q^{c+2i})\bigg).
\endmultline\tag7.16$$

\endproclaim

\demo{Proof}Since it is simpler, we start with the proof of (7.16).
Replace $c$ by $c/2$ and 
substitute $q^{2i-1}$ for $x_i$ in (3.7), $i=1,2,\dots,n$. By (2.12), 
the left-hand side of (3.7) can
be written as a quotient of a sum of two determinants by another
determinant. The determinants, with this choice of $x_i$'s, 
factor by means of (7.4) and (7.7), respectively. 
After simplification we arrive at the expression
$$\multline q^p\frac {\bmatrix p+n-1\\p\endbmatrix_{q^2} \bmatrix
c-p+n-1\\c-p\endbmatrix_{q^2}} {\bmatrix c+n-1\\c\endbmatrix_{q^2}} 
\prod _{i=1} ^{n}\frac {1} {1+q^{2(i-1)}}\prod _{1\le
i<j\le n} ^{}\frac {[c+i+j-2]_{q^2}} {[i+j-2]_{q^2}}\\
\times\bigg((1+q^{c-2p})\prod _{i=1}
^{n-1}(1+q^{c+2i})+(-1)^n
(1-q^{c-2p})\prod _{i=1} ^{n-1}(1-q^{c+2i})\bigg).
\endmultline\tag7.17$$
Note that the only difference from (7.16) is the term $(-1)^n$ in the
last line of (7.17).

On the other hand, by (A.2) the
right-hand side of (3.7) with $c$ replaced by $c/2$ and
with $x_i=q^{2i-1}$, $i=1,2,\dots,n$, equals the generating function for 
tableaux with at most $c$ columns, $p$ of which having parity
different from $n$,
and with only odd entries which lie between $1$ and $2n-1$.
Now we distinguish between $n$ even or odd. If $n$ is even then the
right-hand side of (3.7) with these replacements
equals the generating function for the tableaux
in the second assertion of Theorem~6, and the expressions (7.17) and
(7.16)
agree. If $n$ is odd then the
right-hand side of (3.7), with the above replacements
and with $p$ replaced by $c-p$, 
equals the generating function for the tableaux
in the second assertion of Theorem~6, and, as a few manipulations show, 
the expressions (7.17), with $p$ replaced by $c-p$, and (7.16)
agree.

For the proof of (7.15) we proceed in the same manner, only that
things are more complicated here. 
Now we replace $c$ by $c/2$ and substitute $q^i$ for $x_i$,
$i=1,2,\dots,n$, in (3.7). Again, by (2.12), 
the left-hand side of (3.7) can
be written as a quotient of a sum of two determinants by another
determinant. The determinant in the denominator, with this choice of
$x_i$'s, can be evaluated by
means of (7.5). The second determinant in the numerator, 
with this choice of $x_i$'s, can be evaluated by
means of (7.3). The first determinant in the numerator is evaluated
by means of (7.6). After simplification we arrive at the expression
$$\multline q^p\frac {\bmatrix p+n-1\\p\endbmatrix \bmatrix
c-p+n-1\\c-p\endbmatrix} {2\bmatrix c+n-1\\c\endbmatrix} \prod _{1\le
i<j\le n} ^{}\frac {[c+i+j-2]} {[i+j]}\\
\times\bigg(\sum _{k=0} ^{n}q^{k(k+c-1)}\(q^k\bmatrix n-1\\k\endbmatrix
+q^{-p}\bmatrix n-1\\k-1\endbmatrix\)^2+(-1)^n
(1-q^{c-2p})\prod _{i=1} ^{n-1}(1-q^{c+2i})\bigg).\\
\endmultline\tag7.18$$
Again, note that the only difference from (7.15) is the term $(-1)^n$ in the
last line of (7.18).

On the other hand, by (A.2) the
right-hand side of (3.7) with $c$ replaced by $c/2$ and
with $x_i=q^{i}$, $i=1,2,\dots,n$, equals the generating function for 
tableaux with at most $c$ columns, $p$ of which having parity
different from $n$,
and with entries between $1$ and $n$.
Again, we distinguish between $n$ even or odd. If $n$ is even then the
right-hand side of (3.7) with these replacements
equals the generating function for the tableaux
in the first assertion of Theorem~6, and the expressions (7.18) and
(7.15)
agree. If $n$ is odd then the
right-hand side of (3.7), with the above replacements
and with $p$ is replaced by $c-p$, 
equals the generating function for the tableaux
in the first assertion of Theorem~6, and, as a few manipulations show, 
the expressions (7.18), with $p$ replaced by $c-p$, and (7.15)
agree. In particular, in the last step it is used that
$$\multline
 \sum _{k=0} ^{n}q^{k(k+c-1)}\(q^k\bmatrix n-1\\k\endbmatrix
+q^{-c+p}\bmatrix n-1\\k-1\endbmatrix\)^2
\\=
q^{-c+2p}\sum _{k=0} ^{n}q^{k(k+c-1)}\(q^k\bmatrix n-1\\k\endbmatrix
+q^{-p}\bmatrix n-1\\k-1\endbmatrix\)^2,
\endmultline$$
which is verified by expanding the squares.

This finishes the proof of Theorem~6.\quad \quad \qed
\enddemo

Since the results in Theorem~6 are completely combinatorial
in nature, it would of course be desirable to find a combinatorial proof
of Theorem~6. The papers \cite{\KratAP, \KratAW} (see also
\cite{\KratAQ}) contain 
combinatorial proofs for the assertions in Theorem~5, 
in particular avoiding representation
theory. But these are already quite difficult. The fact that the
expressions in (7.15) and (7.16) are more complicated than those
in (7.11)--(7.14) supports the suspicion that it will be even much harder
to find combinatorial proofs for the assertions in Theorem~6.

\medskip
We conclude by providing an application of Theorem~3, (3.13), (3.14).

\proclaim{Theorem 7}The number of plane partitions of trapezoidal
shape $(2n,2n-2,\dots,\mathbreak 
2)$ with entries between $0$ and $N$ where the
entries on the main diagonal are at least $M$ and exactly $p$ of them
have parity different from $N$ is
$$\multline
\binom np\frac {\dsize\binom{N-M+n}{n+1}} {\dsize\binom{N-M+n+p}{n+1}}
\prod _{1\le i\le j\le n} ^{}\frac {(N+M+i+j)(N-M+i+j)} {(i+j)^2}\\
\quad \quad \text {if }N+M\text { is even}
\endmultline\tag7.19$$
and 
$$\multline
\binom np\frac {\dsize\binom{N+M+n+1}{n+1}}
{\dsize\binom{N+M+1+n+p}{n+1}}
\prod _{1\le i\le j\le n} ^{}\frac {(N+M+1+i+j)(N-M-1+i+j)} {(i+j)^2}\\
\quad \quad \text {if }N+M\text { is odd}.
\endmultline\tag7.20$$
In particular, the number of plane partitions of trapezoidal
shape $(2n,2n-2,\dots,2)$ with entries between $0$ and $N$ where the
entries on the main diagonal are at least $M$ equals
$$2^n
\prod _{1\le i\le j\le n} ^{}\frac {(N+M+i+j)(N-M-1+i+j)} {(i+j)^2}.
\tag7.21$$

\endproclaim

\demo{Proof}We set $c=(N+M)/2$, $d=(N-M)/2$, $x_i=1$, $i=1,2,\dots,n$,
in 
(3.14) and we set $c=(N-M-1)/2$, $d=(N+M-1)/2$, $x_i=1$, $i=1,2,\dots,n$, 
and replace $p$ by $n-p$ in
(3.13). For this choice of $x_i$'s, 
a symplectic character reduces to a number (which is the dimension
of the corresponding irreducible representation of $\Sp(2n,\C)$)
which has a nice closed form (see \cite{\FuHaAA, Ex.~24.20; 
\SunaAE, Theorem~4.5.(1)}), and is therefore easily computed. 
Thus, the left-hand sides of (3.14) and (3.13) turn
into (7.19) and (7.20), respectively. 
On the right-hand sides of (3.14) and (3.13), we have certain sums of
symplectic characters evaluated at $x_i=1$, $i=1,2,\dots,n$. In the
same way as in the proof of Theorem~4, these sums can be interpreted
as the cardinalities of certain sets of plane partitions of
trapezoidal shape. These are exactly the plane partitions in the
first assertion of Theorem~7, they are counted by the specialized
right-hand side of (3.14) if $N+M$ is even, and by the specialized
right-hand side of (3.13) if $N+M$ is odd. This proves the first
assertion of Theorem~7.

To obtain (7.21), we sum the respective expressions in (7.19) and
(7.20) by means of Kummer's very well-poised $_2F_1[-1]$ summation
(cf\. \cite{\SlatAC, (2.3.2.9), Appendix (III.5)}). Thus, Theorem~7
is proved.\quad \quad \qed
\enddemo
It would be easy to generalize the above theorem to trace generating
functions of the type appearing in \cite{\ProcAB, Theorem~1; \KratAO,
sec.~5} by replacing $x_i$ by $q^i$ or $q^{2i-1}$ in (3.13) and
(3.14). We omit the details here.

\head Appendix \endhead

\subhead A1. Partitions and their diagrams\endsubhead
In Section~2 we already defined a {\it partition\/} to be a sequence
$\la=(\la_1,\la_2,\dots,\la_r)$ of integers such that 
$\la_1\ge\la_2\ge\dots\ge\la_r\ge0$. A partition can be viewed
geometrically, in terms of its {\it Ferrers diagram}.
The {\it Ferrers
diagram\/} of a partition $\lambda=(\la_1,\la_2,\dots,\la_r)$
is an array of cells with $r$ left-justified rows and $\lambda_i$
cells in row $i$. Figure~A.a shows the Ferrers
diagram corresponding to $(4,3,3,1)$. We identify partitions with
their Ferrers diagram. For example, if we say ``the first row of the
partition $\la$" then we mean ``the first row of the Ferrers diagram
of the partition $\la$." The {\it conjugate} of a partition $\la$ is the
partition $\la'=(\la'_1,\la'_2,\dots,\la'_{\la_1})$, where $\la'_j$ is
the length of the $j$-th column in the Ferrers diagram of $\la$.

\vskip10pt
\vbox{
$$
\PfadDicke{1pt}
\Pfad(0,0),12112212\endPfad
\Pfad(0,0),22221111\endPfad
\PfadDicke{.5pt}
\Pfad(0,1),1\endPfad
\Pfad(0,2),111\endPfad
\Pfad(0,3),111\endPfad
\Pfad(1,1),222\endPfad
\Pfad(2,1),222\endPfad
\Pfad(3,3),2\endPfad
\hskip2cm
\hbox{\hskip1cm}
\PfadDicke{1pt}
\Pfad(0,0),12112212\endPfad
\Pfad(0,0),22121211\endPfad
\PfadDicke{.5pt}
\Pfad(0,1),1\endPfad
\Pfad(1,2),11\endPfad
\Pfad(2,3),1\endPfad
\Pfad(1,1),2\endPfad
\Pfad(2,1),22\endPfad
\Pfad(3,3),2\endPfad
\hskip2cm
\hbox{\hskip1cm}
\PfadDicke{1pt}
\Pfad(0,0),12112212\endPfad
\Pfad(0,0),22121211\endPfad
\SPfad(0,2),2211\endSPfad
\PfadDicke{.5pt}
\Pfad(0,1),1\endPfad
\Pfad(1,2),11\endPfad
\Pfad(2,3),1\endPfad
\Pfad(1,1),2\endPfad
\Pfad(2,1),22\endPfad
\Pfad(3,3),2\endPfad
\Label\ro{*}(1,2)
\hskip2cm
\hbox{\hskip1cm}
\PfadDicke{1pt}
\Pfad(1,0),12112212\endPfad
\Pfad(1,0),22121211\endPfad
\SPfad(0,0),2222111\endSPfad
\SPfad(0,0),1\endSPfad
\PfadDicke{.5pt}
\Pfad(1,1),1\endPfad
\Pfad(2,2),11\endPfad
\Pfad(3,3),1\endPfad
\Pfad(2,1),2\endPfad
\Pfad(3,1),22\endPfad
\Pfad(4,3),2\endPfad
\Label\ro{*}(2,2)
\hskip2.5cm
$$
\centerline{\eightpoint 
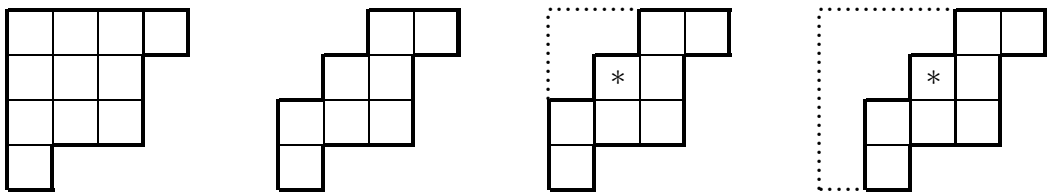
a. Ferrers diagram\quad 
b. skew Ferrers diagrams\quad 
c. $(4,3,3,1)/(2,1)$\quad 
d. $(5,4,4,2)/(3,2,1,1)$}
\vskip7pt
\centerline{\eightpoint Figure A}
}
\vskip10pt

Given two partitions, half-partitions, or orthogonal (half-)partitions 
$\la=(\la_1,\la_2,\mathbreak\dots)$ and
$\mu=(\mu_1,\mu_2,\dots)$, we write $\la+\mu$ for the vector sum
$(\la_1+\mu_1,\la_2+\mu_2,\dots)$. We write $\mu\subseteq \la$ if
$\mu_i\le \la_i$
for all $i$. Given two partitions $\la,\mu$ with $\mu\subseteq \la$,
the {\it skew Ferrers diagram\/} $\la/\mu$ consists of
all cells that are contained in (the Ferrers diagram of) $\la$ 
but not in (the Ferrers diagram of) $\mu$. Figure~A.b
shows the skew Ferrers diagram $(4,3,3,1)/(2,1)$. Of course, it also
shows the skew Ferrers diagram $(5,4,4,2)/(3,2,1,1)$. So, the notation
for skew Ferrers diagrams is not unique. We might even allow vectors
containing negative coordinates for denoting skew Ferrers diagrams,
e.g., $(3,2,2,0)/(1,0,-1,-1)$ for the same skew Ferrers diagram in
Figure~A.b, or vectors with half-integer coordinates, e.g.,
$(9/2,7/2,7/2,3/2)/(5/2,3/2,1/2,1/2)$ for the same skew Ferrers
diagram. We do all this freely in the text, without further notice.
However, when we count columns of skew Ferrers diagrams $\la/\mu$ we
do distinguish between the different notations for the same skew Ferrers
diagram. Each column gets the number that it has as a column of
$\la$. Thus, the first cell in the second row of the skew Ferrers diagram
in Figure~A.b is located in the second column of $(4,3,3,1)/(2,1)$,
see Figure~A.c, in the third column of $(5,4,4,2)/(3,2,1,1)$, see
Figure~A.d, in the first column of $(3,2,2,0)/(1,0,-1,-1)$, in the
$5/2$-th column of $(9/2,7/2,7/2,3/2)/(5/2,3/2,1/2,1/2)$, etc.
A {\it horizontal strip} is a skew Ferrers diagram with no more than
one cell in each of its columns.
A {\it vertical strip} is a skew Ferrers diagram with no more than
one cell in each of its rows.

\subhead A2. $n$-tableaux (ordinary tableaux)\endsubhead
Let $\la,\mu$ be partitions with $\mu\subseteq\la$. An {\it
$n$-tableau of shape $\la$}, respectively {\it of shape $\la/\mu$}, is a
filling of the cells of $\la$, respectively $\la/\mu$, with integers
between $1$ and $n$ such that entries along rows are weakly
increasing and entries along columns are strictly increasing. If we
just say {\it tableau\/} instead of {\it $n$-tableau\/} then we mean
the same but without requiring that the entries are bounded above by
$n$.
Figure~4 shows a $6$-tableau ($7$-tableau, \dots) of shape $(8,8,5,3,2)$.

The {\it weight\/} $\x^T$ for an $n$-tableau $T$ is defined by 
$$\x^T:=x_1^{\#(\text {1's in $T$})}
x_2^{\#(\text {2's in $T$})}\cdots x_n^{\#(\text {$n$'s in $T$})}.
\tag A.1$$
The vector $({\#(\text {1's in $T$})},{\#(\text {2's in
$T$})},\dots,{\#(\text {$n$'s in $T$})})$ of exponents in (A.1) is
called the {\it content\/} of $T$ and is denoted by $\con(T)$. 
For example, the weight and
content of the tableau in Figure~4 are
$x_1^4x_2^4x_3^5x_4^6x_5^5x_6^2$ and $(4,4,5,6,5,2)$, respectively.

As is well-known (see \cite{\MacdAC, (5.12) with $\mu=\boldkey0$;
\SagaAL, Def.~4.4.1}), the irreducible general linear character
(Schur function) $s_n(\la;\x)$ equals the generating function for
all $n$-tableaux of shape $\la$,
$$s_n(\la;\x)=\underset \text {of shape }\la\to
{\sum _{T\text { an $n$-tableau }} ^{}}\x^T,
\tag A.2$$
with $\x^T$ as defined in (A.1).

\subhead A3. $(2n)$-symplectic tableaux\endsubhead
Let $\la=(\la_1,\la_2,\dots,\la_n)$  be a partition. A
{\it $(2n)$-symplectic tableau of shape $\la$} is a $(2n)$-tableau of
shape $2\la=(2\la_1,2\la_2,\dots,2\la_n)$ such that columns 1,~2,
columns 3,~4, \dots, columns $2\la_1-1$,~$2\la_1$ form {\it
$(2n)$-symplectic admissible pairs\/}.
\definition{Definition 1}  A pair $(C,D)$ of two columns of the length
$k\le n$ is called a {\it $(2n)$-symplectic admissible pair\/} if the
following conditions are satisfied:
\roster
\item "(a)"Entries in $C$ and $D$ are between $1$ and $2n$ and in
strictly increasing order.
\item "(b)"If $e$ is in $D$ then $2n+1-e$ is not in $D$. The same
holds for $C$.
\item "(c)"$C$ arises from $D$ by a (possibly empty) sequence of
operations $O$ of the following type: The operation $O$ to be
described applies only to columns $E$ and integers $e_1,e_2$ with
$e_2< e_1\le n$, where $e_1\in E$, $2n+1-e_2\in E$, and for all $t$
between $e_2$ and $e_1$ either $t$ or $2n+1-t$ belongs to $E$. The
operation $O$ itself consists of forming the new column $O(E)$ out of
$E$ by replacing $e_1$ by $e_2$ and $2n+1-e_2$ by $2n+1-e_1$ and
rearranging the new set of entries in strictly increasing order.
\endroster
\enddefinition

For example, Figure~2 displays a $12$-symplectic tableau of shape
$(4,4,4,4,3,3)$. There, the next-to-last column arises from the
last by one operation, as described in item (c) of the Definition,
with $e_1=4$, $e_2=3$.
\remark{Remark} Our description of $(2n)$-symplectic admissible pairs
follows Lakshmibai \cite{\LaksAB}. The description in \cite{\LitPAA,
Appendix~A.2} is equivalent. Note that Littelmann's {\it rows\/} are our
{\it columns}.
\endremark

\medskip
There are a few observations that are immediate from the Definition.

\remark{Observation 1} $C\le D$, meaning that if the entries of $C$
are $c_1,c_2,\dots,c_k$ (from top to bottom) and those of $D$ are
$d_1,d_2,\dots,d_k$ (from top to bottom), then $c_i\le d_i$ for all $i$.
\endremark
\remark{Observation 2} If $(C,D)$ is a $(2n)$-symplectic admissible
pair then the number of entries $\le n$ in $C$ is the same as that in
$D$.
\endremark
\remark{Observation 3} If $C$ and $D$ are columns of length $n$ which
satisfy (a), (b) in the Definition, $C\le D$, and have the same
number of entries $\le n$, then $(C,D)$ is a $(2n)$-symplectic admissible
pair, i.e., condition (c) is satisfied automatically in this
situation.
This is due to the fact that a column
of length $n$ must contain either $t$ or $2n+1-t$, for all
$t=1,2,\dots,n$. Hence, the most obvious way to transform $D$ into $C$
is a legal sequence of operations according to (c): Namely, apply the
operation $O$ described in (c) first with $e_1$ the topmost entry of
$C$ and $e_2$ the topmost entry of $D$ (unless $e_1=e_2$), then
apply $O$ with $e_1$ the next-to-the top entry of
$C$ and $e_2$ the next-to-the-top entry of $D$, etc. 
\endremark

\medskip
The {\it weight\/} $(\x^{\pm1})^S$ for a $(2n)$-symplectic tableau $S$ is defined by 
$$\multline
(\x^{\pm1})^S:=x_1^{\frac {1} {2}\big(\#(\text {1's in $S$})-\#(\text
{$(2n)$'s in $S$})\big)}
x_2^{\frac {1} {2}\big(\#(\text {2's in $S$})-\#(\text {$(2n-1)$'s in
$S$})\big)}
\\\cdots \ 
x_n^{\frac {1} {2}\big(\#(\text {$n$'s in $S$})-\#(\text {$(n+1)$'s in
$S$})\big)}.
\endmultline\tag A.3$$
Again, the vector 
$$\multline
\frac {1} {2}\big({\#(\text {1's in $S$})}-\#(\text
{$(2n)$'s in 
$S$}),\ {\#(\text {2's in
$S$})}-\#(\text {$(2n-1)$'s in $S$}),\\
\dots,\ 
{\#(\text {$n$'s in $S$})}-\#(\text {$(n+1)$'s in $S$})\big)
\endmultline$$ 
of exponents in (A.3) is
called the {\it content\/} of $S$ and is denoted by $\con(S)$. 
For example, the weight and
content of the tableau in Figure~2 are
$x_3x_4^2x_5x_6^{-2}$ and $(0,0,1,2,1,-2)$, respectively.

It is a theorem (see \cite{\DeCoAA, \DePrAA, \LaShAC})
that the irreducible symplectic character
$\sp_{2n}(\la;\x^{\pm1})$ equals the generating function for
all $(2n)$-symplectic tableaux of shape $\la$,
$$\sp_{2n}(\la;\x^{\pm1})=\underset \text {of shape }\la\to
{\sum _{S\text { a $(2n)$-symplectic tableau}} ^{}}(\x^{\pm1})^S,
\tag A.4$$
with $(\x^{\pm1})^S$ as defined in (A.3).

\subhead A4. $(2n+1)$-orthogonal tableaux\endsubhead
Let $\la=(\la_1,\la_2,\dots,\la_n)$  be a partition or half-partition. A
{\it $(2n+1)$-orthogonal tableau of shape $\la$} is a $(2n)$-tableau of
shape $2\la=(2\la_1,2\la_2,\dots,2\la_n)$ such that columns
$2\la_1-1$,~$2\la_1$, columns $2\la_1-3$,~$2\la_1-2$, \dots, form {\it
$(2n+1)$-orthogonal admissible pairs\/}.
\definition{Definition 2}  A pair $(C,D)$ of two columns of the length
$k\le n$ is called a {\it $(2n+1)$-orthogonal admissible pair\/} if the
following conditions are satisfied:
\roster
\item "(a)"Entries in $C$ and $D$ are between $1$ and $2n$ and in
strictly increasing order.
\item "(b)"If $e$ is in $D$ then $2n+1-e$ is not in $D$. The same
holds for $C$.
\item "(c)"$C$ arises from $D$ by a (possibly empty) sequence of
operations that can be either operations of the type that are described in
Definition~1.(c), or operations $O$ of the following type: 
The operation $O$ to be
described applies only to columns $E$ and an entry $e$ of $E$,
$e> n$, where for all $t$
between $n+1$ and $e$, $n+1$ included,
either $t$ or $2n+1-t$ belongs to $E$. The
operation $O$ itself consists of forming the new column $O(E)$ out of
$E$ by replacing $e$ by $2n+1-e$ and
rearranging the new set of entries in strictly increasing order.
\endroster
\enddefinition

For example, Figure~11 displays a $13$-orthogonal tableau of shape
$(7/2,7/2,7/2,\mathbreak 7/2,7/2,1/2)$. 
There, the $4$-th column arises from the
$5$-th by one operation as in item (c) of the Definition above, with
$e=7$, and by one operation as in Definition~1.(c), 
with $e_1=2$, $e_2=1$.
\remark{Remark} Again, our description of $(2n+1)$-orthogonal admissible pairs
follows Lakshmibai \cite{\LaksAB}. The description in \cite{\LitPAA,
Appendix, {\it Spin$_{2m+1}$}-standard Young tab\-leaux} 
is equivalent. Again, note that Littelmann's {\it rows\/} are our
{\it columns}.
\endremark

\medskip
We make similar observations here, also immediate from the Definition.

\remark{Observation 1} $C\le D$, meaning that if the entries of $C$
are $c_1,c_2,\dots,c_k$ (from top to bottom) and those of $D$ are
$d_1,d_2,\dots,d_k$ (from top to bottom), then $c_i\le d_i$ for all $i$.
\endremark
\remark{Observation 2} If $C$ and $D$ are columns of length $n$ which
satisfy (a), (b) in the Definition, and $C\le D$, 
then $(C,D)$ is a $(2n+1)$-orthogonal admissible
pair, i.e., condition (c) is satisfied automatically in this
situation. 
\endremark

\medskip
The {\it weight\/} $(\x^{\pm1})^S$ for a $(2n+1)$-orthogonal tableau $S$ is
again defined by (A.3).
Also here, the vector of exponents in (A.3) is
called the {\it content\/} of $S$ and is denoted by $\con(S)$. 

It is a theorem (see \cite{\LaShAC})
that the irreducible orthogonal character
$\so_{2n+1}(\la;\x^{\pm1})$ equals the generating function for
all $(2n+1)$-orthogonal tableaux of shape $\la$,
$$\so_{2n+1}(\la;\x^{\pm1})=\underset \text {of shape }\la\to
{\sum _{S\text { a $(2n+1)$-orthogonal tableau}} ^{}}(\x^{\pm1})^S,
\tag A.5$$
with $(\x^{\pm1})^S$ as defined in (A.3).

\subhead A5. $(2n)$-orthogonal tableaux\endsubhead
The definition of $(2n)$-orthogonal tableaux is the most intricate
one. We provide the full definition for the sake of completeness.
However, in the current paper we actually need only a special case,
the one that is discussed in Observation~1 below.
The reader who is only interested in this paper's applications of
$(2n)$-orthogonal tableaux can safely skip the full definition and
move on directly to Observation~1. 

We basically reproduce Littelmann's description \cite{\LitPAA,
Appendix, A.3},
with a little modification in the description of
$(2n)$-orthogonal admissible pairs, where we follow \cite{\LaksAB}.
Note again that Littelmann's {\it rows\/} are our
{\it columns}.

Let $\la=(\la_1,\la_2,\dots,\la_n)$  be an $n$-orthogonal 
partition or half-partition. A
{\it $(2n)$-orthogonal tableau of shape $\la$} is a triple
$(S_3,S_2,S_1)$ of $(2n)$-tableaux of respective shapes
$((\la_{n-1}+\la_n)^n)$, $((\la_{n-1}-\la_n)^n)$, and
$(2\la_1-2\la_{n-1},2\la_2-2\la_{n-1},\dots,2\la_{n-2}-2\la_{n-1})$,
subject to conditions (0)--(4) below. Note that when $S_3$, $S_2$, and
$S_1$ are glued together, in that order, an array of shape
$(2\la_1,2\la_2,\dots,2\la_{n-2},2\la_{n-1},2\la_{n-1})$ is obtained.

(0) If $e$ is an element of a column of $S_3$, $S_2$, or $S_1$ then
$2n+1-e$ is not an element of the column.

(1) Columns 1,~2, columns 3,~4, \dots, columns
$2\la_{1}-2\la_{n-1}-1$,~$2\la_{1}-2\la_{n-1}$ of $S_1$ form
$(2n)$-orthogonal admissible pairs (see the Definition below). In
addition, for $i=1,2,\dots,\la_{1}-\la_{n-1}-1$, let the entries of
the $(2i)$-th column of $S_1$ be (from top to bottom)
$k_1,k_2,\dots,k_s$ and let the entries of
the $(2i+1)$-st column be (again, from
top to bottom) $l_1,l_2,\dots,l_t$, $s\le t$. Then for all sequences
$1\le j_1<\dots<j_q\le s$ with
$$n+1-q\le k_{j_1}<\dots<k_{j_q}\le n+q$$
and 
$$n+1-q\le l_{j_1}<\dots<l_{j_q}\le n+q$$
one has $k_{j_1}+\dots+k_{j_q} \equiv l_{j_1}+\dots+l_{j_q}$ mod 2.
(Note that this condition is empty if neither $n$ nor $n+1$ is an
entry in one of the columns.) 

(2) The number of entries $>n$ in each column of $S_2$ is odd. In
addition, denote by $R$ the column (of length $n$) 
that arises from the leftmost
column of $S_1$ by adding all integers $e$, $n+1\le e\le 2n$, that
together with their ``conjugate" $2n+1-e$ do not already appear in
the leftmost column of $S_1$, arranging everything in strictly increasing
order, and replacing the smallest added element, $f$ say,
by $2n+1-f$ in case that the number of entries $>n$
would be odd otherwise. Denote by $x$ the element of $R$ that is
closest to $n+1/2$. To be precise,
let $k_1,k_2,\dots,k_s$ be the entries (from top to bottom)
of the leftmost column of $S_1$. Then $R$ consists of the entries
$k_1,\dots,k_s,l_1,\dots,l_{n-s-1},x$ with the following properties:
\vskip10pt
{\leftskip20pt\rightskip20pt
\noindent 
$2n\ge l_1>\dots>l_{n-s-1}>n$, $l_{n-s-1}>x$, $l_{n-s-1}>2n+1-x$,
$l_i\ne k_j$ and $x\ne k_j$ for all $1\le i\le n-s-1$, $1\le j\le s$,
and if $e\in R$ then $2n+1-e\notin R$. Furthermore, if the number of
integers strictly greater than $n$ in $R\backslash\{x\}$ is odd then
$x>n$ else $x\le n$. 

}
\vskip10pt
\noindent
Let $R'$ be the column with entries $(R\backslash \{x\})\cup\{2n+1-x\}$
in increasing order.
Then the concatenation $S_2\cup R'$ forms a $(2n)$-tableau. By the
concatenation $S_2\cup R'$ we mean that $R'$ is glued from the right 
to $S_2$ such that the topmost entries in each column are
aligned in one row.

(3) The number of entries $>n$ in each column of $S_3$ is even.
Denote by $S_2'$ the tableau of shape $((\la_{n-1}-\la_{n}+1)^n)$ 
obtained from $S_2$ as follows. The rightmost column of $S_2'$ is
$R$. Now assume that $1\le i\le \la_{n-1}-\la_{n}$ and the
$(\la_{n-1}-\la_n+2-i)$-th column of $S_2'$ has already been defined. 
The $(\la_{n-1}-\la_n+1-i)$-th column of $S_2'$ consists of the 
entries of the $(\la_{n-1}-\la_n+1-i)$-th
column of $S_2$ with one entry $e$ replaced by $2n+1-e$, arranged
in strictly increasing order. The entry $e$ to be chosen is the
smallest possible such that rows will be weakly increasing (so that
$S_2'$ indeed becomes a tableau). Then the concatenation $S_3\cup
S_2'$ is a $(2n)$-tableau. By the concatenation $S_3\cup S_2'$ we mean
that $S_3$ is glued to the left to $S_2'$ such that the topmost 
entries in each column are aligned in one row.

\definition{Definition 3} A pair $(C,D)$ of two columns of the length
$k\le n$ is called a {\it $(2n)$-orthogonal admissible pair\/} if the
following conditions are satisfied:
\roster
\item "(a)"Entries in $C$ and $D$ are between $1$ and $2n$ and in
strictly increasing order.
\item "(b)"If $e$ is in $D$ then $2n+1-e$ is not in $D$. The same
holds for $C$.
\item "(c)"$C$ arises from $D$ by a (possibly empty) sequence of
operations that can be either operations of the type that are described in
Definition~1.(c), or operations $O$ of the following type: 
The operation $O$ to be
described applies only to columns $E$ and entries $e_1,e_2$ of $E$,
$n<e_1<e_2$, where for all $t$
between $n+1$ and $e_2$, $n+1$ included,
either $t$ or $2n+1-t$ belongs to $E$. The
operation $O$ itself consists of forming the new column $O(E)$ out of
$E$ by replacing $e_1$ by $2n+1-e_1$ and $e_2$ by $2n+1-e_2$, and
rearranging the new set of entries in strictly increasing order.
\endroster
\enddefinition

\remark{Observation 1} $(2n)$-orthogonal tableaux of shape
$(\la_{n-1},\dots,\la_{n-1},\la_n)$ are just pairs $(S_3,S_2)$ of
$(2n)$-tableaux of respective shapes $((\la_{n-1}+\la_n)^n)$ and 
$((\la_{n-1}-\la_n)^n)$ such that each column of $S_3$ or $S_2$ does
not contain $2n+1-e$ if it contains $e$, such that the number of
entries $>n$ in $S_2$ is odd and each entry in the first row of $S_2$
is $\le n$, such that the number of
entries $>n$ in $S_3$ is even, and where the concatenation
$S_3\cup S_2'$ is a $(2n)$-tableau, 
with $S_2'$ the tableau that arises from $S_2$ by
replacing the topmost element, $e_i$ say, in column $i$
of $S_2$ by its ``conjugate" $2n+1-e_i$, for all
$i=1,2,\dots,\la_{n-1}-\la_{n}$, and
by rearranging the columns in increasing order. For, in case of the
above particular shape the tableau $S_1$ is empty, hence $R$ equals
$\{n+1,n+2,\dots,2n\}$ or $\{n,n+2,\dots,2n\}$, depending on
whether $n$ is even or odd. Thus, $R$ does not restrict $S_2$ in item
(2) above, except that for odd $n$ it forces all the entries in the
first row to be at most $n$. If $n$ is even then all the entries in
the first row have to be $\le n$ as well. For, in each column the
number of entries $>n$ is odd, with $n$ being even. This implies that 
the number of entries $\le n$ has to be odd as well, in particular, it
has to be at least $1$. Finally, in item (3) above, when forming
$S_2'$ out of $S_2$  
always the smallest element in each column can be replaced by its
``conjugate".
\endremark

The {\it weight\/} $(\x^{\pm1})^S$ for a $(2n)$-orthogonal tableau $S$ is
again defined by (A.3).
Also here, the vector of exponents in (A.3) is
called the {\it content\/} of $S$ and is denoted by $\con(S)$. 

It is a theorem (see \cite{\LaShAC})
that the irreducible orthogonal character
$\so_{2n}(\la;\x^{\pm1})$ equals the generating function for
all $(2n)$-orthogonal tableaux of shape $\la$,
$$\so_{2n}(\la;\x^{\pm1})=\underset \text {of shape }\la\to
{\sum _{S\text { a $(2n)$-orthogonal tableau}} ^{}}(\x^{\pm1})^S,
\tag A.6$$
with $(\x^{\pm1})^S$ as defined in (A.3).

\subhead A6. Littelmann's decomposition rule and the
Littlewood--Richardson rule\endsubhead

Littelmann's rule \cite{\LitPAA} for the decomposition of the product
of two general linear, two symplectic, or two special orthogonal
characters can be stated uniformly. (In fact, it is also valid for
the simple, simply connected algebraic groups of type $G_2$ and
$E_6$.) What is needed in the formulation of Littelmann's theorem is
the notion of ({\it dominant}) {\it Weyl chamber\/} associated to each of
the characters.

The ({\it dominant}) {\it Weyl chamber of type $A$\/}, which is associated
to $s_n(\,.\,;\x)$, is the set of points (cf\. \cite{\FuHaAA, p.~215})
$$\{(x_1,x_2,\dots,x_n):x_1\ge x_2\ge \dots\ge x_n\}.\tag A.7$$

The ({\it dominant}) {\it Weyl chamber of type $C$\/}, which is associated
to $\sp_{2n}(\,.\,;\x^{\pm1})$, is the set of points (cf\. \cite{\FuHaAA, p.~243,
(16.5)})
$$\{(x_1,x_2,\dots,x_n):x_1\ge x_2\ge \dots\ge x_n\ge0\}.\tag A.8$$

The ({\it dominant}) {\it Weyl chamber of type $B$\/}, which is associated
to $\so_{2n+1}(\,.\,;\x^{\pm1})$, is the same set (A.8) of points (cf\. \cite{\FuHaAA,
p.~272}).

\NoBlackBoxes
Finally, the ({\it dominant}) {\it Weyl chamber of type $D$\/}, which is associated
to\linebreak 
$\so_{2n}(\,.\,;\x^{\pm1})$, is the set of points (cf\. \cite{\FuHaAA, p.~272})
$$\{(x_1,x_2,\dots,x_n):x_1\ge x_2\ge \dots\ge \vert x_n\vert\}.\tag A.9$$

\BlackBoxes
Furthermore, if $T$ is a tableau, then by $T(\ell)$ we denote the
tableau that consists of the {\it last\/} $\ell$ columns of $T$.

Now we are in the position to state Littelmann's theorem.

\proclaim{Theorem (Littelmann \cite{\LitPAA, Theorem.~(a), p.~346})}
Let $\chi_n(\,.\,)$ be any of the characters $s_n(\,.\,;\x)$,
$\sp_{2n}(\,.\,;\x^{\pm1})$, 
$\so_{2n+1}(\,.\,;\x^{\pm1})$, or $\so_{2n}(\,.\,;\x^{\pm1})$. Then
$$\chi_n(\la)\cdot \chi_n(\mu)=\sum _{T} ^{}\chi_n(\la+\con(T)),
\tag A.10$$
where the sum is over all corresponding tableaux (that is,
$n$-tableaux in case $\chi_n(\,.\,)=s_n(\,.\,;\x)$, $(2n)$-symplectic
tableaux in case $\chi_n(\,.\,)=\sp_{2n}(\,.\,;\x^{\pm1})$, etc.) of shape $\mu$
such that $\la+\con(T(\ell))$ is in the Weyl chamber of the
corresponding type for all $\ell$.
\endproclaim
In case $\chi_n(\,.\,)=s_n(\,.\,;\x)$ this rule translates to the classical
Littlewood--Richardson rule (cf\. \cite{\MacdAC, I, sec.~9; \SagaAL,
(4.26)~+~Theorem~4.9.4}), which we want to describe next.
\definition{Definition 4}  A {\it Littlewood--Richardson filling\/}
({\it LR-filling}) {\it of shape $\nu/\la$ and content $\mu$\/} is an
(ordinary) tableau $F$ of shape $\nu/\la$ and content $\mu$ where the {\it 
Littlewood--Richardson condition\/} ({\it LR-condition}) is
satisfied. The latter condition is the following
\vskip10pt
{\leftskip20pt\rightskip20pt
\noindent
Read the entries of $F$ row-wise from top to bottom and in each row
from right to left. Then at any stage during the reading the number
of $1$'s is greater or equal the number of $2$'s, which in turn is
greater or equal the number of $3$'s, etc.

}
\vskip10pt
\noindent
Furthermore, define the {\it Littlewood--Richardson coefficient\/}
$\LR_{\la,\mu}^{\nu}$ by
$$\LR_{\la,\mu}^{\nu}=\text {number of LR-fillings of shape $\nu/\la$
and content $\mu$}.\tag A.11$$
\enddefinition
Then the {\it Littlewood--Richardson rule\/} reads as follows,
$$s_n(\la;\x)\cdot s_n(\mu;\x)=\sum _{\nu}
^{}\LR_{\la,\mu}^{\nu}\,s_n(\nu;\x).\tag A.12$$

The translation from (A.10) with $\chi_n(\,.\,)=s_n(\,.\,;\x)$ to (A.12)
proceeds as follows. It is basically the same idea as the one we use
in the proof of Proposition~1. What we do is to construct, for any
fixed $\la,\mu,\nu$, a bijection between $n$-tableaux $T$ of shape
$\mu$, with $\nu=\la+\con(T)$, and where $\la+\con(T(\ell))$ is in
the Weyl chamber (A.7) of type $A$ for all $\ell$, and LR-fillings $F$
of shape $\nu/\la$ and content $\mu$. Clearly, this would establish
the equivalence of (A.10) with $\chi_n(\,.\,)=s_n(\,.\,;\x)$ and (A.12). 

Given $T$ as above, we construct a
sequence $F_0,F_1,\dots,F_{\mu_1}$ of fillings by reading $T$
column-wise, from right to left. The desired filling $F$ will then be
defined to be the last filling, $F_{\mu_1}$. Define $F_0$ to be the only
filling of the shape $\la/\la$ (which is of course the empty
filling). Suppose that we already constructed $F_\ell$. To obtain
$F_{\ell+1}$, we add for $i=1,2,\dots,n$ an entry $e$ to row $i$ of
$F_\ell$ if $i$ is an entry occuring in the $\ell$-th last column and
the $e$-th row of $T$. As already announced, we define $F$ to be
$F_{\mu_1}$. 

\subhead A7. Littlewood's branching rules\endsubhead
Here we quote the branching rules for the restriction of Schur
functions to symplectic or orthogonal characters, which we use in
Section~4.
\proclaim{Theorem (Littlewood)} There holds 
(see \cite{\LittAA, App., p.~295; \KingAC, (3.8b); 
\SunaAE, Theorem~3.13}
$$s(\la;\x)=\sum _{\nu} ^{}\sp(\nu;\x)\sum _{\mu,\ \mu'\text { even}}
^{}\LR_{\mu,\nu}^{\la}\tag A.13$$
(`$\mu'$ even' means that all the columns of $\mu$ are even),
and (see \cite{\LittAA, p.~240, (II); \KingAC, (3.8a); 
\SunaAE, Theorem~3.16})
$$s(\la;\x)=\sum _{\nu} ^{}\o(\nu;\x)\sum _{\mu,\ \mu\text { even}}
^{}\LR_{\mu,\nu}^{\la}\tag A.14$$
(`$\mu$ even' means that all the rows of $\mu$ are even), with 
$\LR_{\mu,\nu}^{\la}$ the Littlewood--Richardson coefficients as
defined in (A.11).

\endproclaim

\enddocument